%% file: main_arxiv.tex
\documentclass{amsart}
\input{sections/preamble}
\input{notation}

\theoremstyle{remark}
\newtheorem{remark}{Remark}[section]

\newcommand{\RA}[1]{{\color{black}#1}}
\newcommand{\RB}[1]{{\color{black}#1}}
\newcommand{\RALL}[1]{{\color{black}#1}}
\begin{document}

\title[DATA-DRIVEN EFR for NSE]{Data-driven 
Optimization for the Evolve-Filter-Relax regularization 
of convection-dominated flows}

\author[A. Ivagnes, M.Strazzullo, M. Girfoglio, T. Iliescu, G. Rozza]{Anna Ivagnes$^{a,1}$, Maria Strazzullo$^{a,2}$, Michele Girfoglio$^{3}$, Traian Iliescu$^{4}$, Gianluigi Rozza$^{a,1}$}
\date{\today}
\address{$^1$ International School of Advanced Studies, SISSA, Mathematics Area, Via Bonomea, 265, 34136, Trieste, Italy}
\address{$^2$ Politecnico di Torino, Department of Mathematical Sciences ``Giuseppe Luigi Lagrange'', Corso Duca degli Abruzzi, 24, 10129, Torino, Italy}
\address{$^3$ Department of Engineering, University of Palermo, Via delle Scienze, Ed. 7, Palermo, 90128, Italy}
\address{$^4$ Department of Mathematics, Virginia Tech, Blacksburg, VA 24061, USA}
\address{$^{\MakeLowercase{a}}$ INdAM research group GNCS member}

\maketitle

\begin{abstract}
    \input{sections/abstract}
\end{abstract}

\section{Introduction}
\label{sec:intro}
\input{sections/introduction}

\section{The EFR algorithm}
\label{sec:methods}
\input{sections/methods-EFR}

\section{The Data-Driven optimized EFR algorithms}
\label{sec:methods-opt}
\input{sections/methods-Opt-intro}

\section{Numerical results}
\label{sec:results}
\input{sections/results-intro}

\subsection{Numerical results for $\chi$-Opt-EFR}
\label{sec:results-chi-opt}
\input{sections/results-Opt-EFR}

\subsection{Numerical results for $\delta$-Opt-EF}
\label{sec:results-delta-opt}
\input{sections/results-Opt-EF}

\subsection{Numerical results for $\delta \chi$--Opt-EFR}
\label{sec:results-double-opt}
\input{sections/results-D-Opt-EFR}

\section{Conclusions and outlook}
\label{sec:conclusions}
\input{sections/conclusions}

\input{sections/extra-final-sections}

\bibliography{biblio}
\bibliographystyle{abbrv}

\appendix

\input{sections/supplementary}

\end{document}

%% file: sections/preamble.tex
\usepackage[foot]{amsaddr}
\usepackage{comment}
\usepackage[utf8]{inputenc}
\usepackage{graphicx}
\usepackage{verbatim}
\usepackage{subfig}
\usepackage{amsmath}
\usepackage{mathtools}
\usepackage{siunitx}
\usepackage[dvipsnames, table]{xcolor}
\usepackage[top=1in, bottom=1.25in, left=1.25in, right=1.25in]{geometry}
\usepackage[hidelinks]{hyperref}
\usepackage{algorithm}
\usepackage{algpseudocode}

\usepackage{ulem}

\newcommand{\highlight}[1]{\colorbox{gray!30}{\parbox{\dimexpr\linewidth-2\fboxsep\relax}{#1}}}

\usepackage{booktabs}
\usepackage{multirow}
\usepackage{multicol}
\usepackage{pifont}

%% file: notation.tex
\def\PP{{{\rm l}\kern - .15em {\rm P} }}
\def\PN2{{\PP_{N}-\PP_{N-2}}}

\newcommand{\bu}{\boldsymbol{u}}

\newcommand{\bw}{\boldsymbol{w}}

\newcommand{\bx}{\boldsymbol{x}}

\newcommand{\optEFRloc}{\ensuremath{\chi\text{-Opt-EFR}_{\text{local}}^{\bu}}}
\newcommand{\optEFRlocGRAD}{\ensuremath{\chi\text{-Opt-EFR}_{\text{local}}^{\bu, \nabla \bu}}}
\newcommand{\optEFRglob}{\ensuremath{\chi\text{-Opt-EFR}_{\text{global}}^{\bu}}}
\newcommand{\optEFRglobGRAD}{\ensuremath{\chi\text{-Opt-EFR}_{\text{global}}^{\bu, \nabla \bu}}}
\newcommand{\optEFRglobPRESS}{\ensuremath{\chi\text{-Opt-EFR}_{\text{global}}^{\bu, \nabla \bu, p}}}

\newcommand{\optEFglob}{\ensuremath{\delta\text{-Opt-EF}_{\text{global}}^{\bu}}}
\newcommand{\optEFglobGRAD}{\ensuremath{\delta\text{-Opt-EF}_{\text{global}}^{\bu, \nabla \bu}}}
\newcommand{\optEFglobPRESS}{\ensuremath{\delta\text{-Opt-EF}_{\text{global}}^{\bu, \nabla \bu, p}}}

\newcommand{\DoptEFRglobGRAD}{\ensuremath{\delta \chi\text{-Opt-EFR}_{\text{global}}^{\bu, \nabla \bu}}}
\newcommand{\DoptEFRglobPRESS}{\ensuremath{\delta \chi\text{-Opt-EFR}_{\text{global}}^{\bu, \nabla \bu, p}}}

\newcommand{\lossLOC}{\ensuremath{\mathcal{L}_{\text{local}}}}
\newcommand{\lossLOCu}{\ensuremath{\mathcal{L}_{\text{local}}^{\bu}}}
\newcommand{\lossLOCgrad}{\ensuremath{\mathcal{L}_{\text{local}}^{\nabla \bu}}}
\newcommand{\lossGLOB}{\ensuremath{\mathcal{L}_{\text{global}}}}
\newcommand{\lossGLOBu}{\ensuremath{\mathcal{L}_{\text{global}}^{\bu}}}
\newcommand{\lossGLOBgrad}{\ensuremath{\mathcal{L}_{\text{global}}^{\nabla \bu}}}
\newcommand{\lossGLOBp}{\ensuremath{\mathcal{L}_{\text{global}}^{p}}}

\newcommand{\deleted}[1]{{}}


%% file: sections/abstract.tex
  Numerical stabilization techniques are often employed in under-resolved simulations of convection-dominated flows to improve accuracy and mitigate spurious oscillations. Specifically, the evolve--filter--relax (EFR) algorithm is a framework which consists in evolving the solution, applying a filtering step to remove high-frequency noise, and relaxing through a convex combination of filtered and original solutions. The stability and accuracy of the EFR solution strongly depend on two parameters, the filter radius $\delta$ and the relaxation parameter $\chi$. 
    Standard choices for these parameters are usually fixed in time, and related to the full order model setting, 
    i.e., the grid size for $\delta$ and the time step for $\chi$.  
   \RA{The key novelties with respect to the standard EFR approach are:}
    (i) time-dependent parameters \RA{$\delta(t)$ and $\chi(t)$}, and (ii) data-driven adaptive optimization of the parameters in time, considering a fully-resolved simulation as reference.
    In particular, we propose three different classes of optimized-EFR \RA{(\textbf{Opt-EFR})} strategies, aiming to optimize one or both parameters.
    The 
    new \RA{Opt-EFR} strategies are tested 
    in the under-resolved simulation of a turbulent flow past a cylinder at $Re=\num{1000}$.
    \RA{The \RA{Opt-EFR} proved to be more accurate than standard approaches by up to 99$\%$,} 
   while maintaining a similar computational time.
    In particular, \RA{the key new finding of our analysis is that such accuracy can be obtained only
     if the optimized objective function includes: (i) a \emph{global} metric (as the kinetic energy), and (ii) \emph{spatial gradients}' information.}

%% file: sections/introduction.tex

%



It is well known that, in the convection-dominated (e.g., turbulent), under-resolved regime (i.e., when the number of degrees of freedom is not high enough to accurately represent the complex underlying dynamics), classical numerical discretizations of the Navier-Stokes equations (NSE) generally yield spurious numerical oscillations~\cite{berselli2006mathematics,layton2012approximate}.
To alleviate this inaccurate behavior, these discretizations are generally equipped with numerical stabilization, large eddy simulation \RA{\cite{choi2025performance, li2025optimizing}}, or 
Reynolds-Averaged Navier--Stokes (RANS) models~\cite{berselli2006mathematics,layton2012approximate,rebollo2014mathematical}\RA{\cite{guo2025effects, ahmadi2024predicting, martinez2025rans, zou2025study}}.
Regularized models (see \cite{layton2012approximate} for a review) are a class of numerical stabilization strategies that leverage spatial filtering to increase the numerical stability of the underlying model and eliminate the spurious numerical oscillations.
One of the most popular regularized model is the {\it evolve-filter-relax (EFR)} algorithm, which consists of the following three steps:
(i) {\it Evolve}: \, Use the NSE discretization to advance the velocity 
approximation at the current time step to an intermediate approximation at the next time step.
(ii) {\it Filter}: \, Use a spatial filter to smooth the intermediate approximation in step (i).
The goal of this step is to eliminate the spurious oscillations from the intermediate velocity approximation in step (i). 
(iii) {\it Relax}: \, Set the velocity approximation at the next time step to be a convex combination of the filtered velocity in step (ii) and the intermediate velocity approximation in step (i).
The goal of this step is to prevent an overdiffusive velocity approximation, which could occur in step (ii) if, e.g., the filter radius is too large.
The EFR strategy has been successfully used in conjunction with classical numerical discretization techniques, ranging from finite element (FE), to spectral element, to finite volume methods, to reduced order modeling (see, e.g., \cite{boyd2001chebyshev, canuto1988some,ervin2012numerical,fischer2008nek5000,fischer2001filter,germano1986differential, germano1986differential2,girfoglio2019finite, girfoglio2021fluid,xie2018evolve,mullen1999filtering,NekROM,olshanskii2013connection,pasquetti2002comments, rezaian2023predictive,takhirov2018computationally,tsai2025time,strazzullo2023new,strazzullo2022consistency})
in challenging applications, ranging from nuclear engineering design, to the quasi-geostrophic equations describing the large scale ocean circulation, to cardiovascular modeling (see, e.g, 
\cite{girfoglio2019finite, girfoglio2023novel, besabe2025linear, xuthesis, quaini2024bridging, xu2020backflow,tsai2023parametric,layton2012approximate}).
The EFR's popularity is due to its accuracy and low computational cost.
We emphasize, however, that probably the main reasons for the EFR's popularity in practical applications are its {\it ease of implementation},
{\it modularity}, and {\it minimally intrusive} 
character.
Indeed, the EFR algorithm is embarrassingly easy to implement:
Given access to the stiffness matrix of the underlying numerical discretization, the spatial filter (and, thus, the EFR algorithm) can be implemented in a matter of minutes.
Furthermore, the EFR algorithm is modular and 
minimally intrusive: 
Given a {\it legacy} NSE code, one can very easily add a modular spatial filter subroutine and a modular relaxation step without having to modify the legacy code. 


The EFR algorithm has been proven to be a popular and effective strategy in important applications.
Its further development, however, currently faces a major roadblock:
The EFR algorithm involves two parameters (i.e., the filter radius and the relaxation parameter) that need to be carefully calibrated in order to obtain accurate approximations.
This parameter calibration is critical for the EFR's success:
If the filter radius is too large, the EFR model is overdiffusive and the results are inaccurate (smeared) \cite{strazzullo2024variational}.
On the other hand, if the filter radius is too small, there is not enough dissipation and the EFR model can actually blow up.
Similarly, a relaxation parameter that is too small or too large can yield grossly inaccurate EFR results.
Thus, finding the optimal filter radius and relaxation parameter is critical for the EFR's success.
We note that the sensitivity of the EFR results to the model parameters have been well documented (see, e.g., the numerical investigations in \cite{bertagna2016deconvolution,girfoglio2019finite,strazzullo2022consistency, strazzullo2024variational,tsai2025time}).
We also note that the importance of the EFR's parameter calibration is not a surprise since parameter sensitivity is a central issue for most numerical stabilizations \cite{berselli2006mathematics,layton2012approximate,ZoccolanStrazzulloRozza20242,ZoccolanStrazzulloRozza2024}.

Recognizing the importance of finding the optimal filter radius and relaxation parameter in the EFR model, several approaches have been proposed early on.
For example, in \cite{girfoglio2019finite, strazzullo2022consistency}, the EFR filter radius was chosen to be the mesh size of the underlying FE discretization.
Another choice was proposed in \cite{bertagna2016deconvolution}, where the EFR filter radius was chosen to be the Kolmogorov lengthscale.
Similarly, several choices for the EFR relaxation parameter have been proposed.
For example, choosing the relaxation parameter proportional to the time step is a popular choice in the EFR literature \cite{bertagna2016deconvolution,strazzullo2022consistency}.
We note, however, that heuristic arguments have been used to motivate higher values for the relaxation parameter \cite{bertagna2016deconvolution,strazzullo2023new,strazzullo2022consistency}.

We emphasize that the vast majority of the current choices for the EFR filter radius and relaxation parameter yield {\it constant} values \RB{in both space and time} for these two parameters (an important exception is the nonlinear filtering proposed in \cite{girfoglio2019finite}.)
Using constant values for the parameters imposes significant limitations on the EFR's performance.
For example, if the flow transitions from a laminar regime to a turbulent regime, the EFR's filter radius and/or time relaxation parameter should increase to adapt to this transition.
We emphasize that this is simply not possible with the current parameter choices.

\RA{The key novelties of this contribution are:
\begin{itemize}
    \item[] (i) the introduction of time-dependent EFR parameters;
    \item[] (ii) computing the optimal parameters at different time instances with an adaptive {\it data-driven optimization} strategy, leveraging the available reference data;
    \item[] (iii) the inclusion of \emph{global} quantities into the objective function, such as the kinetic energy; and
    \item[] (iv) the inclusion of \emph{spatial gradients} into the objective function.
\end{itemize}
}

The proposed strategies are tested in the marginally-resolved simulation of a flow past a cylinder at a Reynolds number $Re=1000$.

We emphasize that our new strategy is fundamentally different from the current approaches: 
\begin{itemize}
    \item This is the first time that data-driven modeling is leveraged to determine the EFR parameters.
    \item Furthermore, the new data-driven parameters vary in time, which is in stark contrast to the vast majority of current methods used to determine the EFR parameters.
\end{itemize}

The rest of the paper is organized as follows:
In Section~\ref{sec:methods}, we present the EFR strategy and discuss the critical role played by the EFR parameters, the filter radius and relaxation parameter.
In Section~\ref{sec:methods-opt}, we introduce \RA{the \textbf{Opt-EFR} approaches,} novel data-driven optimization algorithms for the EFR strategy.
Specifically, we introduce three \RA{Opt-EFR} strategies:
(i) In the $\chi$-Opt-EFR approach, we fix the filter radius and optimize the relaxation parameter.
(ii) In the $\delta$-Opt-EFR approach, we fix the relaxation parameter, and optimize the filter radius.
(iii) Finally, in the $\delta\chi$-Opt-EFR approach, we optimize both the relaxation parameter and the filter radius. 
In Section~\ref{sec:results}, we perform a numerical investigation of the three new data-driven EFR strategies in the numerical simulation of incompressible flow past a circular cylinder at Reynolds number $Re=1000$ on a relatively coarse mesh.
For comparison purposes, we also consider the standard EFR approach.
As benchmark for our numerical investigation, we use a DNS on a fine mesh.
Finally, in Section~\ref{sec:conclusions} we draw conclusions and outline future research directions.

%% file: sections/methods-EFR.tex
This Section is dedicated to a review of the state of the art related to the methodology we investigate.
In particular, we briefly recall the discretized Navier--Stokes formulation and the Evolve-Filter-Relax (EFR) approach in Section 
\ref{subsec:NSE}. For more details, we refer the reader to \cite{wells2017evolve, strazzullo2022consistency, girfoglio2021fluid,bertagna2016deconvolution}. Moreover, in 
Section \ref{subsec:role-parameters}, we highlight the bottleneck of the EFR formulation, and the motivations to enhance it by optimization.

\subsection{The Navier--Stokes Equations and the EFR algorithm}
\label{subsec:NSE}
As a mathematical model, we consider the \textbf{incompressible Navier--Stokes equations (NSE)}.
Given a fixed domain $\Omega \subset \mathbb R^{D}$, with $D=2,3$, we consider the motion of an incompressible fluid having velocity $\bu \doteq \bu(\bx, t) \in \mathbb U$ and pressure $p \doteq p(\bx, t) \in \mathbb Q$ whose dynamics is represented by the incompressible NSE:
\begin{equation}
\label{eq:NSE}
\begin{cases}
 \displaystyle \frac{\partial \bu}{\partial t} + (\bu \cdot \nabla) \bu - \nu \Delta \bu + \nabla p = 0 & \text{in }  \Omega \times (t_0, T), \\
\nabla \cdot \bu = 0  & \text{in }  \Omega \times (t_0, T), \\
\bu = \bu_D & \text{on } \partial \Omega_D \times (t_0, T), \\
\displaystyle -p \ \boldsymbol n + \nu \frac{\partial \boldsymbol \bu}{\partial \boldsymbol n} = 0  & \text{on }  \partial \Omega_N \times (t_0, T), \\
\end{cases}
\end{equation}
endowed with the initial condition $\bu = \bu_0$ in $\Omega \times \{ t_0 \}$, where $\partial \Omega_D \cup \partial \Omega_N = \partial \Omega$, $\partial \Omega_D \cap \partial \Omega_N = \emptyset$, $\nu$ is the kinematic viscosity, and $\mathbb U$ and $\mathbb Q$ are suitable Hilbert function spaces. The functions $\bu_D$ and $\bu_0$ are given. 
The flow regime is defined by the adimensional \textbf{Reynolds number}
\begin{equation}
\label{eq:Re}
Re \doteq \frac{UL}{\nu},
\end{equation}
where $U$ and $L$ represent the characteristic velocity and length scales of the system, respectively. 
When the Reynolds number is large, the inertial forces 
dominate the viscous forces: this setting is generally referred to as the convection-dominated regime.

In the numerical tests, we employ the \textbf{Finite Element Method (FEM)} and a backward differentiation formula of order 1 (BDF1) for space and time discretization, respectively. 
In under-resolved or marginally-resolved numerical simulations, namely when employing a coarse grid to run the FEM simulations, standard spatial discretization techniques yield spurious numerical oscillation in convection-dominated regimes.

To alleviate this behavior, we equip the simulation with the \textbf{Evolve-Filter-Relax (EFR)} algorithm.
\RB{To introduce the notation, we consider  a time step $\Delta t$. 
Let $t^n = t_0 + n\Delta t$ for $n = 0, \dots, N_T$, and $T = t_0 + N_T\Delta t$. }
We denote the semi-discrete FEM velocity and pressure \RB{at time $t^n$} with $\bu\RB{^n} \in \mathbb U^{N_h^{\bu}}$ and $p\RB{^n} \in \mathbb Q^{N_h^p}$, respectively, where ${N_h^{\bu}}$ and ${N_h^p}$ are 
their degrees of freedom, respectively.
The EFR algorithm at time $t^{n+1}$ reads:

\begin{eqnarray*}      
&	\text{\bf (I)}& \text{\emph{Evolve}:} \quad 
\begin{cases}
\begin{split}
     &\frac{\bw^{n + 1} - \bu^n}{\Delta t} + (\bw^{n+1} \cdot \nabla) \bw^{n+1} +\\
    &-\nu \Delta \bw^{n+1} + \nabla p^{n+1} = 0 
\end{split} & \text{in } \Omega \times \{t_{n+1}\}, \vspace{2mm}\\
\nabla \cdot \bw^{n+1} = 0 & \text{in } \Omega \times \{t_{n+1}\}, \vspace{1mm}\\
\bw^{n+1} = \bu_D^{n+1} & \text{on } \partial \Omega_D \times \{t_{n+1}\}, \vspace{2mm}\\
\displaystyle -p^{n+1} \cdot \boldsymbol n + \frac{\partial \boldsymbol \bw^{n+1}}{\partial \boldsymbol n} = 0  & \text{on } \partial \Omega_N \times \{t_{n+1}\}. \\
\end{cases}
            \label{eqn:ef-rom-1}\nonumber \\[0.3cm]
            &	\text{\bf (II)} &\text{\emph{Filter:}} \quad
\begin{cases} 
        	 -2 \delta^2 \, \Delta \overline{\bw}^{n+1} +  \overline{\bw}^{n+1} = \bw^{n+1}& \text{in } \Omega \times \{t_{n+1}\}, \vspace{2mm}\\
\overline{\bw}^{n+1} = \bu^{n+1}_D &\text{on } \partial \Omega_D \times \{t_{n+1}\}, \vspace{2mm}\\
\displaystyle \frac{\partial \overline{\bw}^{n+1}}{\partial \boldsymbol n} = 0  & \text{on } \partial \Omega_N \times \{t_{n+1}\}.
\end{cases}
	\label{eqn:ef-rom-2} \nonumber \\[0.3cm]
            &	\text{\bf (III)} &\text{\emph{Relax:}} \qquad 
        	   \bu^{n+1}
            = (1 - \chi) \, \bw^{n+1}
            + \chi \, \overline{\bw}^{n+1} \, ,
            \label{eqn:ef-rom-3}\nonumber
\end{eqnarray*}

for $\chi \in [0,1]$ relaxation parameter. We denote by $\bw$ and $\overline \bw$ the evolved and the filtered velocity, respectively. 
Solving the first step (I) with $\bw^{n+1} = \bu^{n+1}$ equals solving the standard NSE.
In step (II), we employ a \emph{differential filter} (DF) filter radius $\delta$ 
(i.e., the spatial lengthscale over the spatial filter takes action). 
The success of DFs 
is due to their mathematical robustness \cite{berselli2006mathematics}.
Indeed, the DF leverages an elliptic operator and eliminates the small scales (i.e., high frequencies) from the input velocity, smoothing it and acting as a spatial filter. 
Step (III) is a relaxation step: a convex combination of the evolved and filtered velocity is performed at a new time step. In particular, the extreme cases are:
\begin{itemize}
    \item[(i)] $\chi=1$: we only retain the filtered solution $\bu^{n+1}=\overline{\bw}^{n+1}$;
    \item[(ii)] $\chi=0$: the final solution is exactly the solution of the incompressble NSE (I), i.e., no filtering contribution is included and $\bu^{n+1}=\bw^{n+1}$.
\end{itemize}
In other words, the relaxation parameter $\chi$ is used to regulate the amount of the filter contribution \cite{ervin2012numerical, fischer2001filter,mullen1999filtering}, 
to increase the accuracy of possible over-diffusive results of Step (II), as shown for example in \cite{bertagna2016deconvolution}. 
{We remark that we call {\it noEFR simulation} the numerical solution of the discretized Navier--Stokes equations with no regularization applied, following the notation of \cite{strazzullo2022consistency}.}

\subsection{On the role of the EFR parameters}
\label{subsec:role-parameters}
A natural question arises in this setting: \emph{what are the optimal choices for $\delta$ and $\chi$?} In the literature some scales arguments have been proposed (we refer the reader to \cite{strazzullo2022consistency,bertagna2016deconvolution}  and the references therein for more details on the topic):
\begin{itemize}
    \item[$\circ$] Common choices for the filter radius $\delta$ are the minimum mesh-scale diameter $h_{min}$ or the Kolmogorov  length scale $\eta =L\cdot Re^{\frac{3}{4}}$.
    Another choice is the energy-based lengthscale~\cite{mou2023energy}.
    \item[$\circ$] For the relaxation, the scaling $\chi \sim \Delta t$ is commonly used \cite{ervin2012numerical}. However, its choice highly depends on the specific test case taken into account and higher values have been employed following heuristic arguments in \cite{bertagna2016deconvolution,girfoglio2019finite, strazzullo2022consistency}.
\end{itemize}
While the choice of the filter radius is related to physical quantities, the choice of the relaxation parameter is based on experience and it might completely change with respect to the test case considered. In this context, many adjustments and improvements can be made.
For example, time-dependent flow features can be better captured with a time-dependent parameter $\chi(t)$, or the differential filter may have local action \cite{girfoglio2019finite}. Clearly, {\it The main bottleneck of the EFR procedure relies upon tuning the filter radius and the relaxation parameters.}

Indeed, in the standard EFR algorithm, the parameters' values need to be carefully tuned in a pre-processing stage to improve the accuracy of the solution:
\begin{itemize}
    \item[(i)] If the values of $\delta$ and/or $\chi$ are \textbf{too large}, it may happen that the solution is \textbf{over-diffusive} and it does not manifest any turbulent behavior.
    \item[(ii)] If the values of $\delta$ and/or $\chi$ are \textbf{too small}, the filter action may be not adequate enough, leading to \textbf{spurious oscillations}, as happens in the noEFR scenario.
\end{itemize}
 Moreover, as mentioned above, the flow dynamics may not be properly captured considering parameters constant in time, especially when dealing with turbulent flows or with test cases including a laminar-to-turbulent transition. For instance, we may want larger filter parameters (hence, higher $\delta$ and/or $\chi$) in the laminar-to-turbulent transition, as well as in fully turbulent regimes.
 
 For this purpose, having time-dependent parameters would help in better capturing the flow dynamics, avoiding, for instance, overdiffusion phenomena.
 In this manuscript, \textbf{for the first time, we introduce a data-driven strategy} to compute the time-dependent parameters $\delta(t)$ and $\chi(t)$, which allow us to better capture the flow dynamics in unsteady turbulent regimes.
In this framework, 
we aim at investigating different \textbf{
data-driven optimization algorithms}, with the ultimate goal of developing 
\emph{an 
adaptive  
strategy} to optimize $\delta(t)$ and $\chi(t)$.

%% file: sections/methods-Opt-intro.tex
In this section, we aim at obtaining \emph{optimal values} for the EFR parameters, enhancing standard approaches by adaptive optimization algorithms. 
In this way, we increase the accuracy of the reg-simulation with respect to a Direct Numerical Simulation (DNS) simulation taken as a reference.

Our investigation goes beyond constant parametric values and focuses on three specific \RA{Opt-EFR} strategies\RA{. Specifically, calling $\boldsymbol{\mu}^n=\boldsymbol{\mu}(t^n)$ the parameter(s) optimized at each time step $t^n$, the three approaches are:
} 

\begin{itemize}
    \item[\textbf{1.}] $\chi$-\textbf{Opt-EFR}, aimed to find the optimal $\chi(t)$ for fixed $\delta=\eta$\RA{, hence $\boldsymbol{\mu}^n=\chi^n$;}
    \item[\textbf{2.}] $\delta$-\textbf{Opt-EF}, aimed to find the optimal $\delta(t)$, with no relaxation, namely $\chi=1$\RA{, hence $\boldsymbol{\mu}^n=\delta^n$;}
    \item[\textbf{3.}] $\delta \chi$-\textbf{Opt-EFR}, aimed to find the optimal $\delta(t)$ and $\chi(t)$\RA{, hence $\boldsymbol{\mu}^n=(\delta^n, \chi^n)$.}
\end{itemize}

The presented optimization procedures are based on a data-driven algorithm trained taking as reference some pre-computed data. In our specific setting, the training data are represented by the reference velocity and pressure fields of a simulation in a resolved regime.\\
We refer to those fields as $\bu_{ref} \in \mathbb{U}^{M_h^{\bu}}$ and $p_{ref} \in \mathbb{Q}^{M_h^p}$, where $M_h^{\bu} \gg N_h^{\bu}$ and $M_h^p \gg N_h^p$. Such data are provided by a simulation over 
a refined mesh yielding $M_h^{\bu}+ M_h^{p}$ degrees of freedom (dofs). The regularization, instead, is applied to a coarser mesh with respect to the DNS one, with $N_h^{\bu} + N_h^{p}$ dofs.
Moreover, it is important to remark 
the \textbf{generality} and the \textbf{modularity} of the proposed algorithms, which can be applied to any experimental setting.
The motivations behind the choices of the test cases of this contribution are postponed to the numerical sections.

\RA{A schematic comparison between the standard EFR and the Opt-EFR methods is displayed in Figure \ref{fig:opt-efr-flowchart}. The standard EFR is characterized by a standard time stepping procedure, as described in Section \ref{subsec:NSE} and showed in Figure \ref{fig:opt-efr-flowchart} (A).
On the other hand, the Opt-EFR algorithms perform an additional optimization stage to find the optimal $\delta(t)$ and/or $\chi(t)$, as highlighted in Figure \ref{fig:opt-efr-flowchart} (B).

\begin{figure*}[htpb!]
\begin{center}
    \subfloat[Standard EFR]{
    \includegraphics[width=0.9\linewidth]{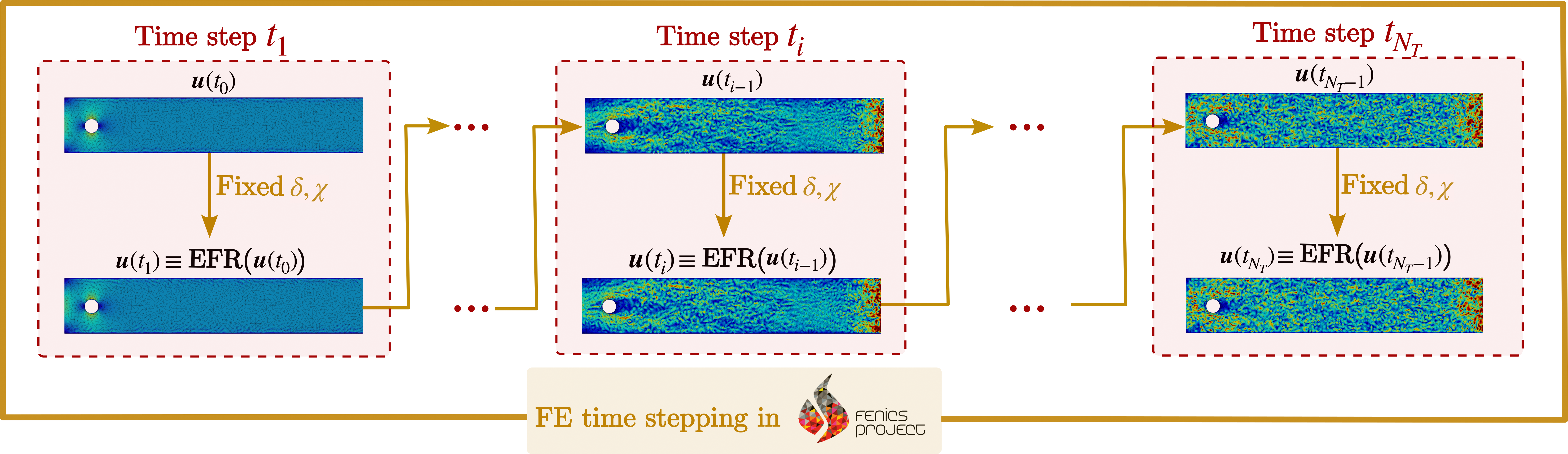}}\\
    \subfloat[Opt-EFR]{
    \hspace{7mm}\includegraphics[width=\linewidth]{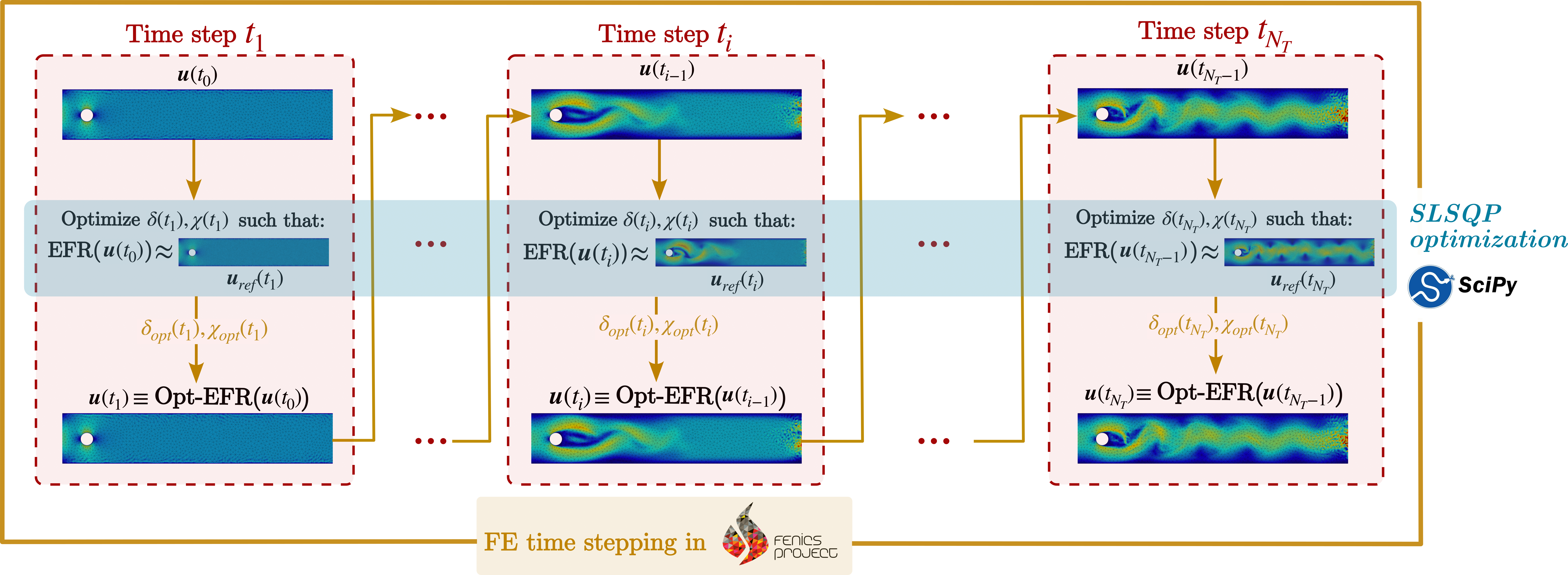}}
    \caption{\RA{Flowchart of the Opt-EFR and standard EFR approaches.}}
    \label{fig:opt-efr-flowchart}
\end{center}
\end{figure*}
}

\subsection{The objective functions}

\RA{In the optimization algorithms, we will consider either \emph{local} or \emph{global} contributions. Specifically, the \emph{local} approach refers to a pointwise evaluation and minimization of the discrepancy between the reference and the predicted quantity. On the other hand, the \emph{global} approach consists of the evaluation and minimization of the discrepancy of the global energy.

We refer to \emph{local} quantities if we minimize the mean squared error (MSE) between the DNS reference solution (projected on the coarse mesh) and the updated EFR solution.

The local objective function can then be written as follows:

\begin{equation}
    \label{eq:local-loss-general}
\lossLOC(\boldsymbol{\mu}^{n}) = w_{\bu} \, \lossLOCu(\boldsymbol{\mu}^{n}) + w_{\nabla \bu} \, \lossLOCgrad(\boldsymbol{\mu}^{n}),
\end{equation}
where the individual contributions are:
\begin{equation}
\label{eq:local-contributions-general}
\begin{split}
    &\lossLOCu(\boldsymbol{\mu}^{n}) = \text{MSE}\left( \bu^{n}(\boldsymbol{\mu}^n), \bu_{ref}^{n} \right), \\
   &\lossLOCgrad(\boldsymbol{\mu}^{n}) = \text{MSE}\left( \nabla \bu^{n}(\boldsymbol{\mu}^n), \nabla \bu_{ref}^{n} \right),
   \end{split}
\end{equation}
with $w_{\bu}$ and $w_{\nabla \bu}$ scalar quantities.
}
\RALL{As pointed out in the results in Section \ref{sec:results-chi-opt}, the local contributions may not be appropriate metrics to optimize in this test case. In particular, the Opt-EFR results and the DNS may differ a lot in terms of qualitative and quantitative behaviors, even if the local MSE is relatively small.}
\RA{For this reason, we investigate the following  global objective function:

\begin{equation}
    \label{eq:global-loss-general}
\lossGLOB(\boldsymbol{\mu}^{n}) = w_{\bu} \lossGLOBu(\boldsymbol{\mu}^{n}) + w_{\nabla \bu} \lossGLOBgrad(\boldsymbol{\mu}^{n}) + w_{p} \lossGLOBp(\boldsymbol{\mu}^{n}),
\end{equation}

where the global contributions are:
\begin{equation}
\label{eq:global-contributions-general}
\begin{split}
    &\lossGLOBu(\boldsymbol{\mu}^n) = \left| \dfrac{\|\bu^{n}(\boldsymbol{\mu}^n)\|^2_{L^2(\Omega)} - \|\bu_{ref}^{n}\|^2_{L^2(\Omega)}}{\|\bu_{ref}^{n}\|^2_{L^2(\Omega)}} \right|, \\
   & \lossGLOBgrad(\boldsymbol{\mu}^{n}) =\left| \dfrac{\| \nabla \bu^{n}(\boldsymbol{\mu}^n)\|^2_{L^2(\Omega)} - \| \nabla \bu_{ref}^{n}\|^2_{L^2(\Omega)}}{\| \nabla \bu_{ref}^{n}\|^2_{L^2(\Omega)}} \right|,\\
   &\lossGLOBp(\boldsymbol{\mu}^{n}) = \left| \dfrac{\|p^{n+1}(\boldsymbol{\mu}^n)\|^2_{L^2(\Omega)} - \|p_{ref}^{n+1}\|^2_{L^2(\Omega)}}{\|p_{ref}^{n+1}\|^2_{L^2(\Omega)}} \right|,
\end{split}
\end{equation}
for $w_p$ scalar quantity. In particular, we set $w_{\bu}\neq 0$ and $w_{\nabla \bu}=w_{p}=0$ in algorithms of type $(\cdot)^{\bu}_{\text{global}}$, 
$w_{\bu}\neq 0, w_{\nabla \bu}\neq 0$, and  $w_{p}=0$ in $(\cdot)^{\bu, \nabla \bu}_{\text{global}}$ algorithms, while all the contributions are activated in $(\cdot)^{\bu, \nabla \bu, p}_{\text{global}}$.

An overview of the methods investigated in this work, depending on the different loss functional considered, is provided in Table \ref{tab:acronyms-all}.

\begin{table*}[h]
\caption{Acronyms and main features of the Opt-EFR algorithms. Gray cells indicate that the algorithm does not present that feature.}
\label{tab:acronyms-all}
\centering
{
\begin{tabular}{cccccc}
\\ \hline
\multirow{2}{*}{\textbf{Acronym}} & \multicolumn{2}{c}{\textbf{Objective function type}} & \multicolumn{3}{c}{\textbf{Contributions in objective function}}
\\ \cline{2-3}\cline{4-6}
& local & global & $\bu$ discrepancy & $\nabla \bu$ discrepancy & $p$ discrepancy \\ \hline
Opt-EFR$^{\bu}_{\text{local}}$ & $\checkmark$ & \cellcolor{gray!50} & $\checkmark$ & \cellcolor{gray!50}& \cellcolor{gray!50}\\ \hline
Opt-EFR$^{\bu, \nabla \bu}_{\text{local}}$ & $\checkmark$ & \cellcolor{gray!50} & $\checkmark$ & $\checkmark$ & \cellcolor{gray!50}\\ \hline
Opt-EFR$^{\bu}_{\text{global}}$ & \cellcolor{gray!50} & $\checkmark$ & $\checkmark$ & \cellcolor{gray!50} & \cellcolor{gray!50}\\ \hline
Opt-EFR$^{\bu, \nabla \bu}_{\text{global}}$ & \cellcolor{gray!50} & $\checkmark$ & $\checkmark$ &$\checkmark$& \cellcolor{gray!50}\\ \hline
Opt-EFR$^{\bu, \nabla \bu, p}_{\text{global}}$ & \cellcolor{gray!50} & $\checkmark$ & $\checkmark$ &$\checkmark$ &$\checkmark$  \\ \hline
\end{tabular}}
\end{table*}

It is worth highlighting that the expression of \lossGLOBp{} in \eqref{eq:global-contributions-general} penalizes the pressure obtained as a solution of the evolve step at time instance $t^{n+1}$. Therefore, the optimization algorithms penalizing the pressure contribution, namely $(\cdot)^{\bu, \nabla \bu, p}_{\text{global}}$, are computationally more expensive than the others, since a further evolve step should be performed for each evaluation.

Indeed, in Steps (II) and (III) at time $t^{n}$ apply only to the velocity field. Thus, the pressure contribution in Equations \eqref{eq:global-contributions-general} needs to be computed considering the pressure at time $t^{n+1}$, which is the first pressure field affected by the optimization of $\boldsymbol{\mu}$ at step $t^n$.
We investigate the effect of the pressure contribution only in the \emph{global} approach because the numerical results for the \emph{local} method are not accurate (as showed in Section \ref{sec:results}).

For additional details on the three classes $\chi$-Opt-EFR, $\delta$-Opt-EF, and $\delta \chi$-Opt-EFR, we refer the reader to the supplementary Section \ref{sec:supp-methodology}.

\remark{It is worth highlighting that the scalar weights $w_{\bu}$, $w_{\nabla \bu}$, and $w_p$ are taken exactly equal to 1 when they are different from $0$, namely when the corresponding terms are activated in the objective functions \eqref{eq:local-loss-general} and \eqref{eq:global-loss-general}.
Other choices can be made, but finding the optimal choices for the loss weights was beyond the goal of the manuscript.}

\subsection{The optimization process}

\begin{algorithm*}[htpb!]
\caption{Pseudocode for algorithms Opt-EFR optimizing both parameters. Row 9 (underlined) is considered only for algorithms of type Opt-EFR$^{\bu, \nabla \bu, p}_{\text{global}}$. The optimization part is highlighted in gray.}
\label{alg:opt-all}
\begin{algorithmic}[1]
\State $\bu_{ref}^n$, $n=0, \dots, N_T - 1$; $N_h^u$; $\bu^0$, $\Delta t_{opt}$, $\boldsymbol{\mu}^0$, $\mathrm{A}$ range of parameter(s). \Comment{$\text{\it Inputs needed}$}
\For{$n \in [0, \dots, N_T - 1]$}
\State (I) $\rightarrow \bw^{n+1}$ \Comment{$\text{\it Perform Evolve step}$}
\State $\bu_{ref}^{n+1} \gets \bu_{ref}^{n+1}\,  | \, \mathbb{U}^{N_h^u}$ \Comment{$\text{\it Project the reference field into the coarse mesh}$}
\State $\boldsymbol{\mu}^{n+1} \gets \boldsymbol{\mu}^n$  \Comment{$\text{\it Initialize parameter(s)}$}
\highlight{\If {$\bmod(n, \Delta t_{opt}) = 0$} \Comment{$\text{\it Perform optimization every }\Delta t_{opt}$} 
\State $iter=0$
\Repeat 
        \State (II)+(III) $\rightarrow {\bu}^{n+1}$ \Comment{$\text{\it Filter and Relax step}$}
        \State \underline{(I) $\rightarrow p^{n+2}$} \Comment{$\text{\it \underline{(Second) Evolve step}}$}
        \State Optimize $\mathcal{L}(\boldsymbol{\mu}^{n+1})$, subject to $\boldsymbol{\mu}^{n+1} \in \mathrm{A}$ 
        \State $iter \gets iter + 1$
        \State Output current solution: $\boldsymbol{\mu}^{n+1}$
    \Until{Convergence or maximum iterations reached}
    \EndIf}
    \State Solve (II)+(III) with the optimized $\boldsymbol{\mu}^{n+1} \rightarrow$ Update $\bu^{n+1}$  \Comment{$\text{\it Perform Filter and Relax}$}
\EndFor 
\end{algorithmic}
\end{algorithm*}

The algorithms in Table \ref{tab:acronyms-all} perform at each time instance Step (I), and look for an optimal time-dependent $\delta(t)$ and/or $\chi(t)$ to be employed in Step (II) and/or (III).

The process is summarized in Algorithm \ref{alg:opt-all}, where $y | \mathbb Y$ denotes the projection of the field $y$ onto the space $\mathbb Y$.

Algorithm \ref{alg:opt-all} presents the general pseudocode that can be applied for all the algorithms presented in Table \ref{tab:acronyms-all}, by replacing the objective function $\mathcal{L}(\chi^{n+1})$ with the specific functional of the proposed algorithm. 

Moreover, we employ a \emph{sequential least squares programming} (SLSQP) algorithm \cite{
philip2010sequential,wilson1963simplicial}, which is an iterative method for constrained nonlinear optimization.

The only considered constraint is $\boldsymbol{\mu}^{n+1} \in \mathrm{A}$, where $\mathrm{A}$ is the admissibility range of the parameter(s).

In particular, we consider:
\begin{enumerate}
    \item $\chi^{n+1} \in [0, 1]$ in $\chi$-Opt-EFR algorithms, namely $\mathrm{A} = [0, 1]$;
    \item $\delta^{n+1} \in \RALL{[0.019 \eta, 1.9\eta]}$ in $\delta$-Opt-EF algorithms, namely $\mathrm{A} = \RALL{[0.019 \eta, 1.9\eta]}$; and
    \item $\delta^{n+1} \in \RALL{[0.019 \eta, 1.9\eta]}$ and $\chi^{n+1} \in [0, 1]$ in $\delta\chi$-Opt-EFR algorithms, namely $\mathrm{A} = \RALL{[0.019 \eta, 1.9\eta]}\times [0,1]$;
\end{enumerate}

By adding the constraint $\delta^{n+1} \in \RALL{[0.019 \eta, 1.9\eta]}$, the optimization problem seeks the optimal parameter within the range spanning an order of magnitude lower to one higher than the Kolmogorov length scale. This window has been established after a preliminary sensitivity analysis on the parameter
$\delta$, as it provides more physical results. The unconstrained optimization may indeed converge to values that are too large and
far from the standard literature choices.

The underlined row (row 10 in the pseudocode \ref{alg:opt-all})
is considered only in the Opt-EFR$^{\bu, \nabla \bu, p}_{\text{global}}$ algorithms, since it performs an additional evolve step in order to evaluate the pressure at the following time instance $p^{n+2}$. 
As an additional consideration, we consider as the initial value of the parameter(s) the common choice in the literature \cite{strazzullo2022consistency}, namely $\delta^0=\eta$ for the filter radius, and $\chi_0=K \Delta t$, with $K=5$, for the relaxation parameter.

\remark{
It is worth specifying that the parameters are not optimized every time step. Instead, the time window between two consecutive optimizations is $\Delta t_{opt}>\Delta t$. Specifically, $\Delta t_{opt} = k \Delta t$, with $k \in \{10,20, \cdots, 90,100\}$. 

Every $\Delta t_{opt}$, the algorithm takes into account the updated solution of the system and therefore needs to be trained \emph{on the fly}. The parameter is fixed to its latest optimal value in the time window between two subsequent optimization steps.
In the numerical results, we will also evaluate how the optimization step influences the accuracy and efficiency of the solution.}

\paragraph*{Software details}
The EFR method is implemented using the finite-element Python-based open-source library FEniCS~\cite{alnaes2015fenics}. We perform the optimization algorithms using the open-source Python-based library Scipy~\cite{virtanen2020scipy}.
}

%% file: sections/results-intro.tex
This section 
shows the numerical results obtained by applying the 
methods proposed in Sections \ref{sec:methods} and \ref{sec:methods-opt}.
We consider the test case of the incompressible flow past a circular cylinder at Reynolds number $Re=1000$ constant in time.

We consider as reference in the optimization algorithms a DNS 
on a fine mesh, and our goal is to minimize the discrepancy between the optimized-EFR fields and the reference DNS fields.

\begin{figure*}[htpb!]
    \centering
    \subfloat[DNS on fine mesh - velocity]{\includegraphics[width=0.5\textwidth, trim={3cm 7cm 3cm 7cm}, clip]{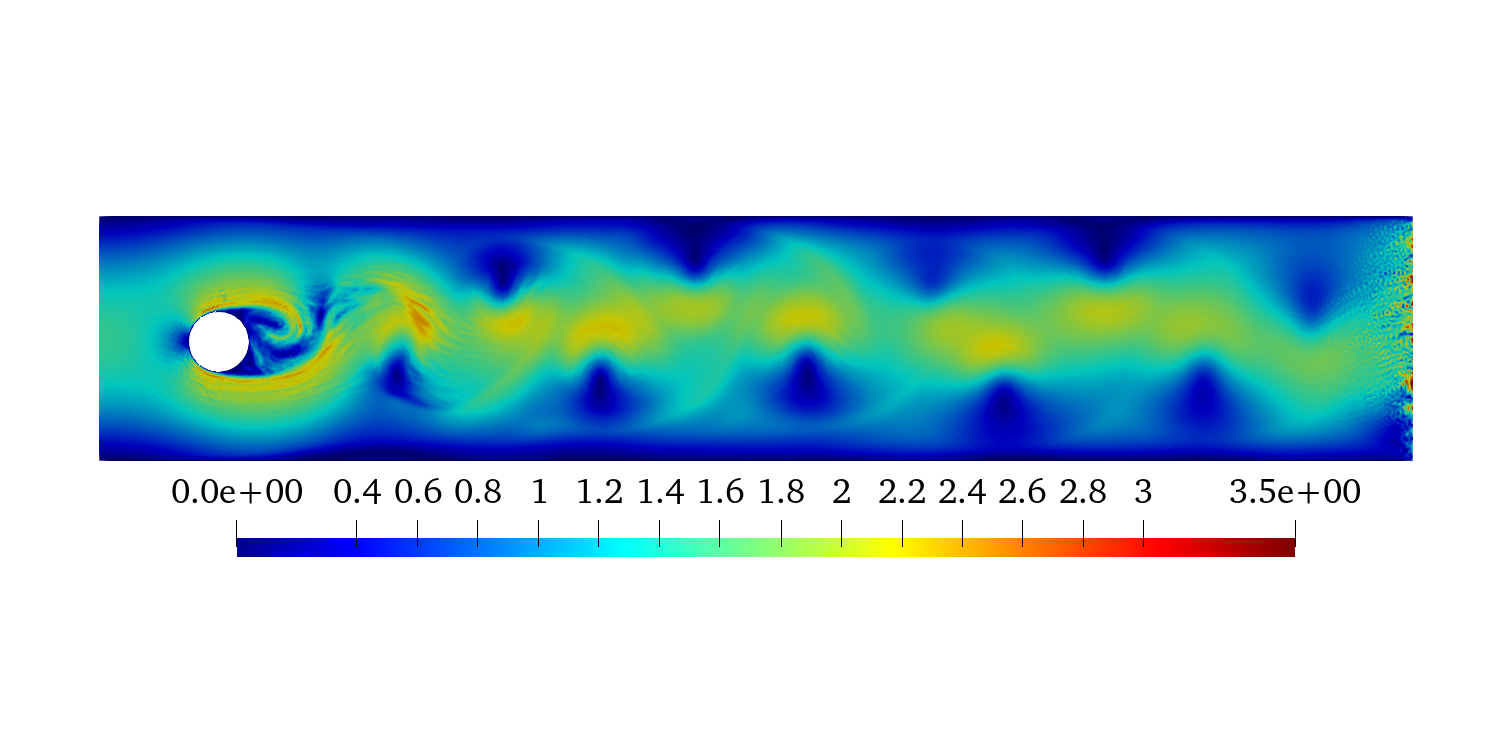}}    \subfloat[DNS on fine mesh - pressure]{\includegraphics[width=0.5\textwidth, trim={3cm 7cm 3cm 7cm}, clip]{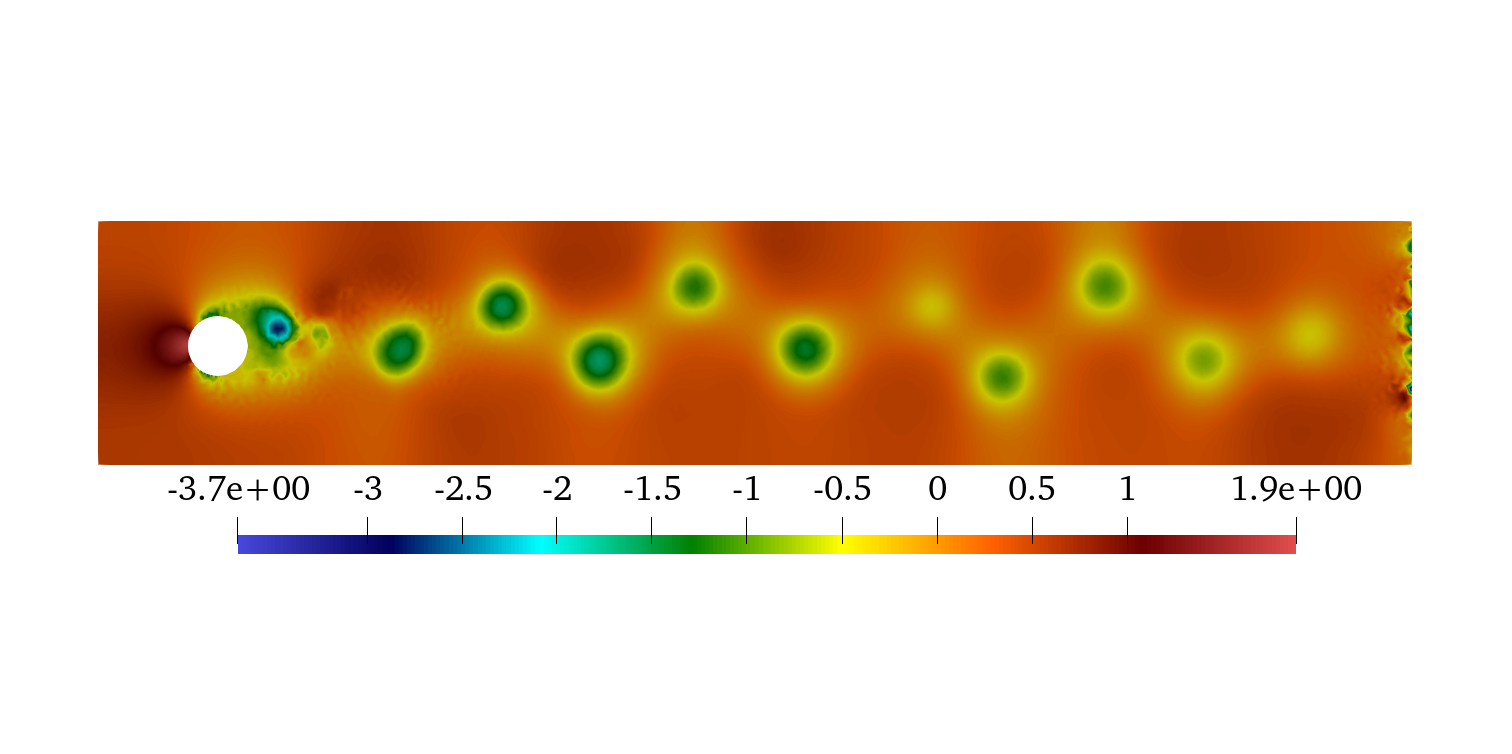}}\\    
    \subfloat[noEFR - velocity]{\includegraphics[width=0.5\textwidth, trim={3cm 6.5cm 3cm 7cm}, clip]{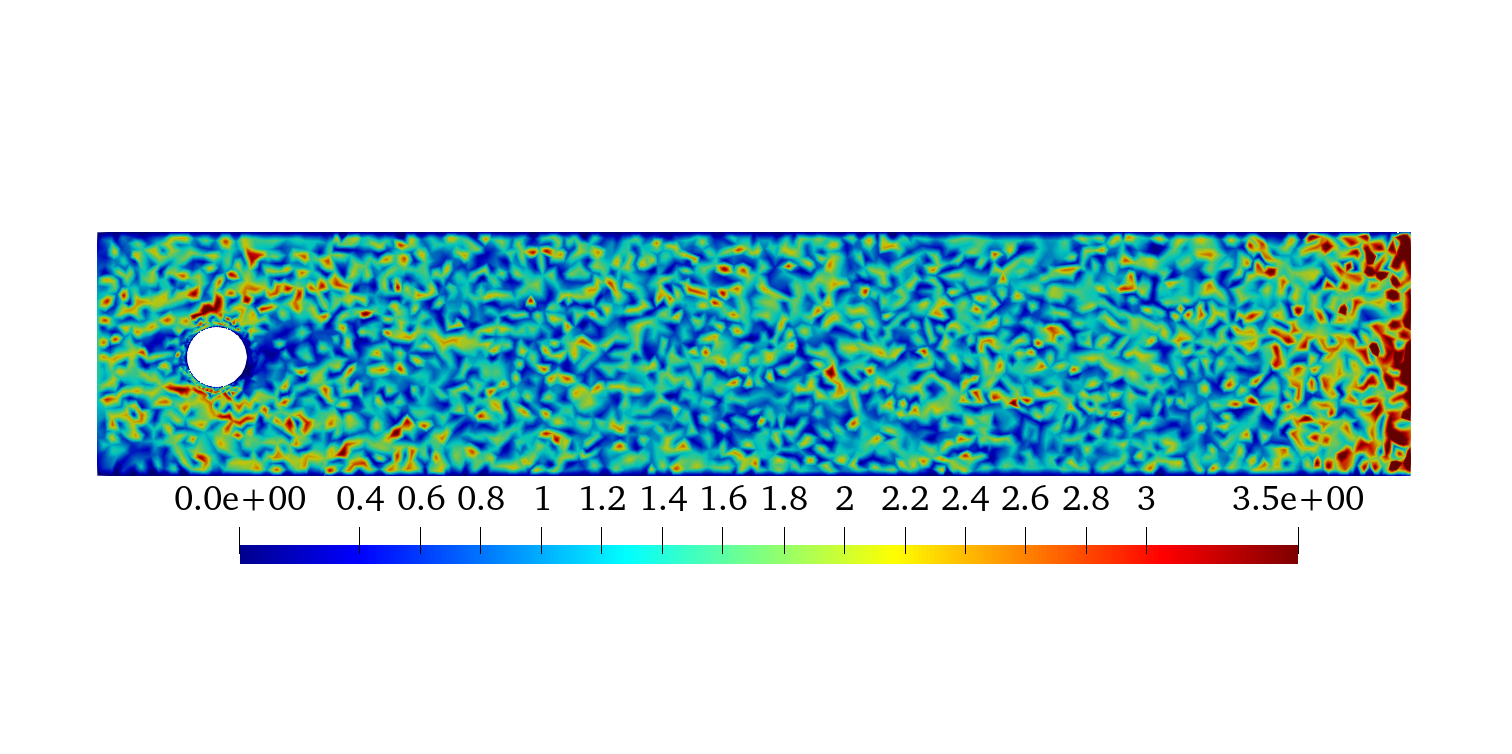}}    \subfloat[noEFR - pressure]{\includegraphics[width=0.5\textwidth, trim={3cm 6.5cm 3cm 7cm}, clip]{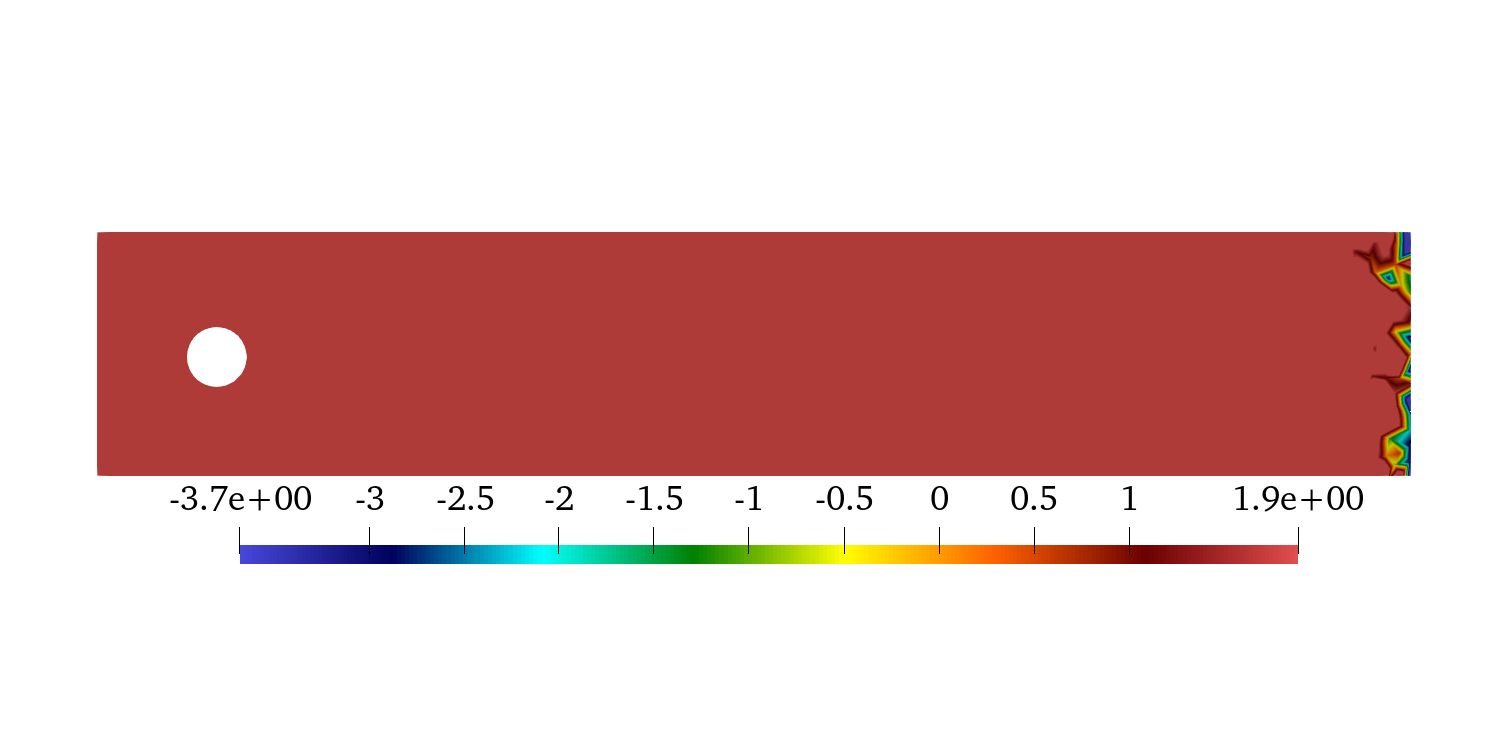}}
    \\    
    \subfloat[Standard EFR ($\delta=\eta$, $\chi=5 \Delta t$) - velocity]{\includegraphics[width=0.5\textwidth, trim={3cm 6.5cm 3cm 7cm}, clip]{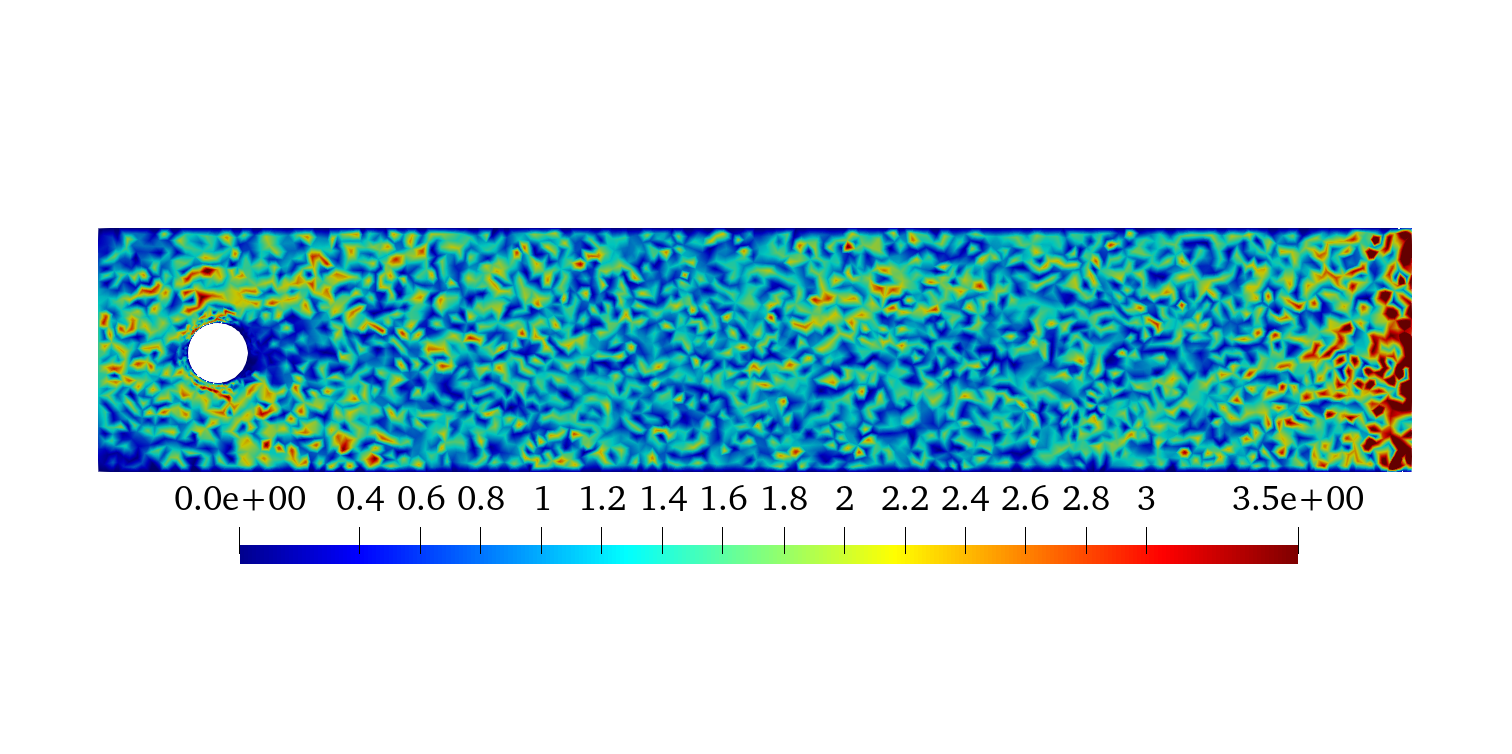}}    \subfloat[Standard EFR ($\delta=\eta$, $\chi=5 \Delta t$) - pressure]{\includegraphics[width=0.5\textwidth, trim={3cm 6.5cm 3cm 7cm}, clip]{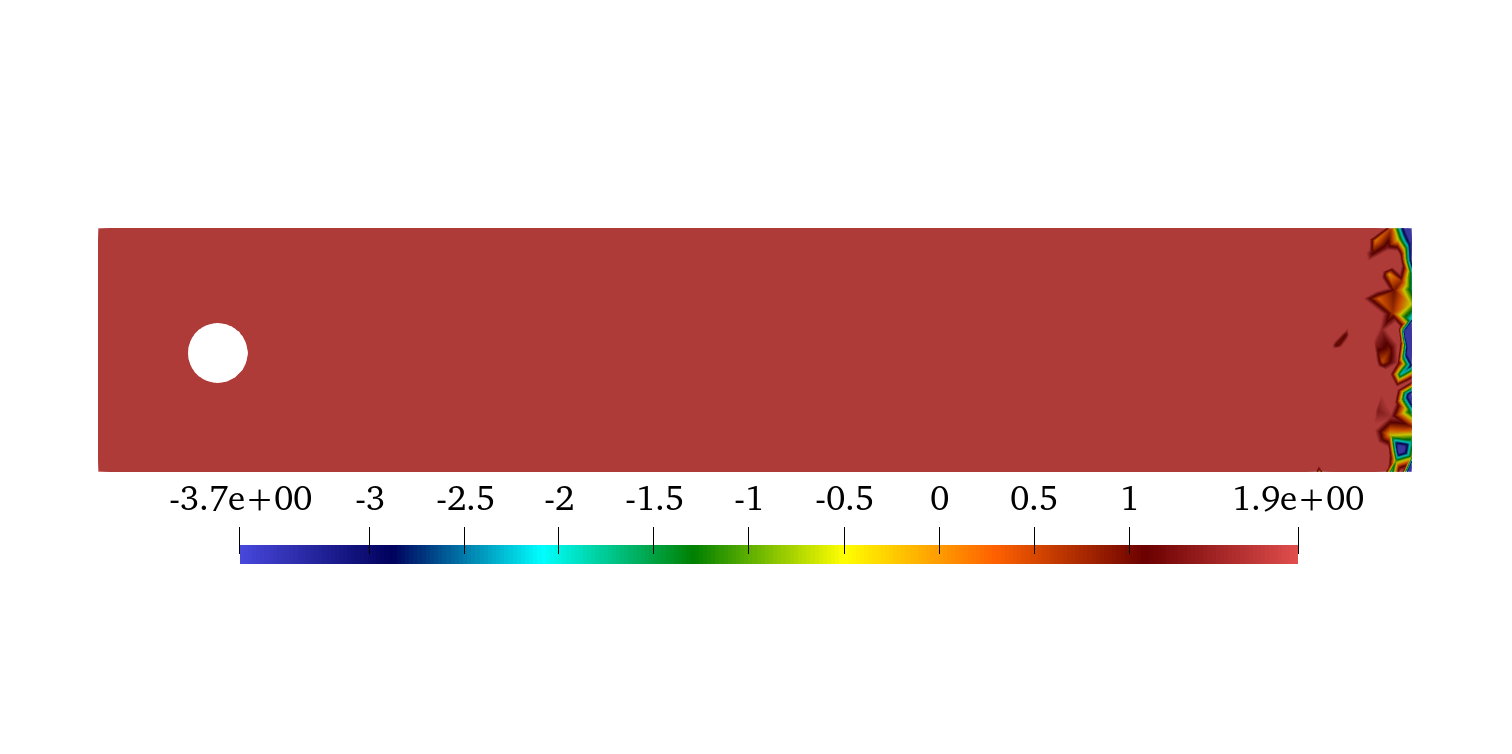}}
    \\    
    \subfloat[Standard EF ($\delta=\eta$) - velocity]{\includegraphics[width=0.5\textwidth, trim={3cm 6.5cm 3cm 7cm}, clip]{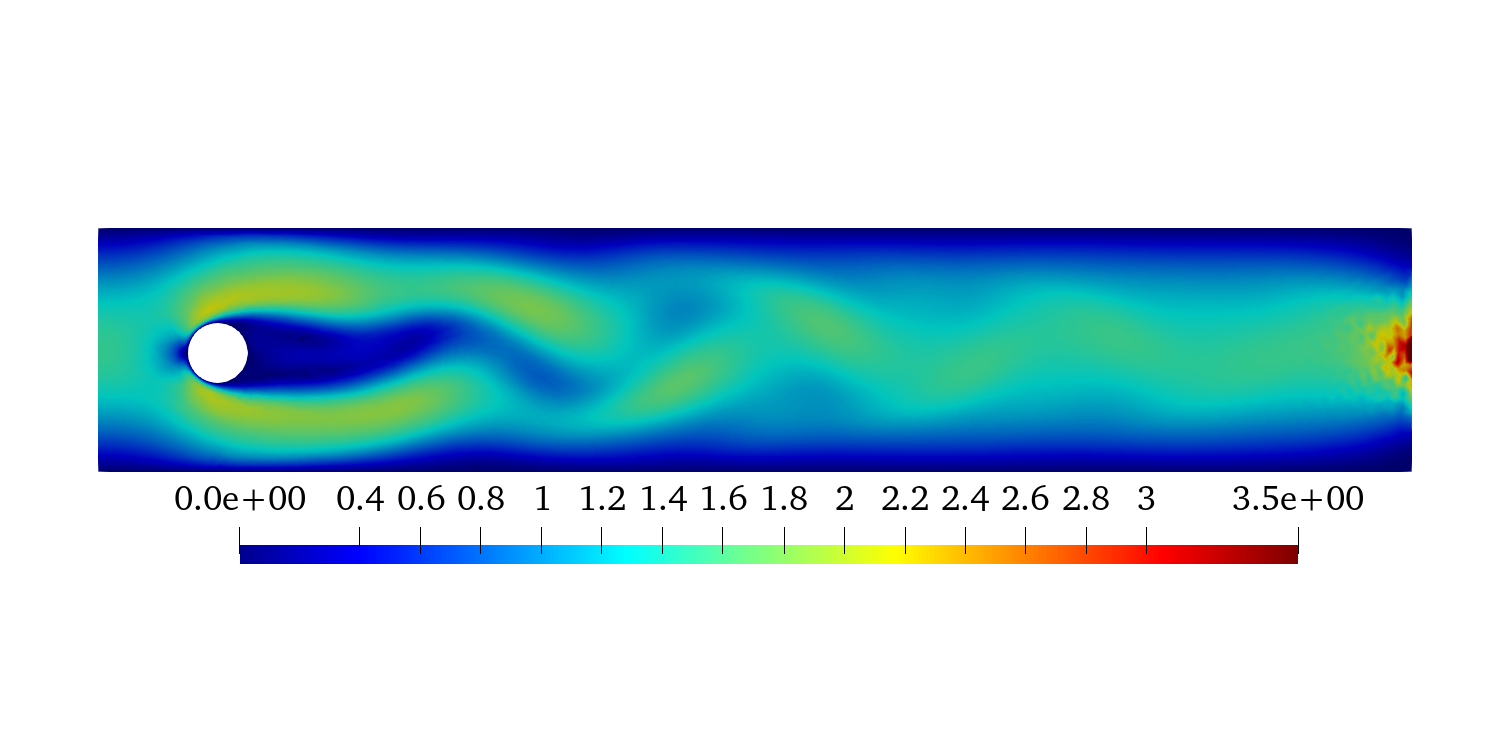}}   
    \subfloat[Standard EF ($\delta=\eta$) - pressure]{\includegraphics[width=0.5\textwidth, trim={3cm 6.5cm 3cm 7cm}, clip]{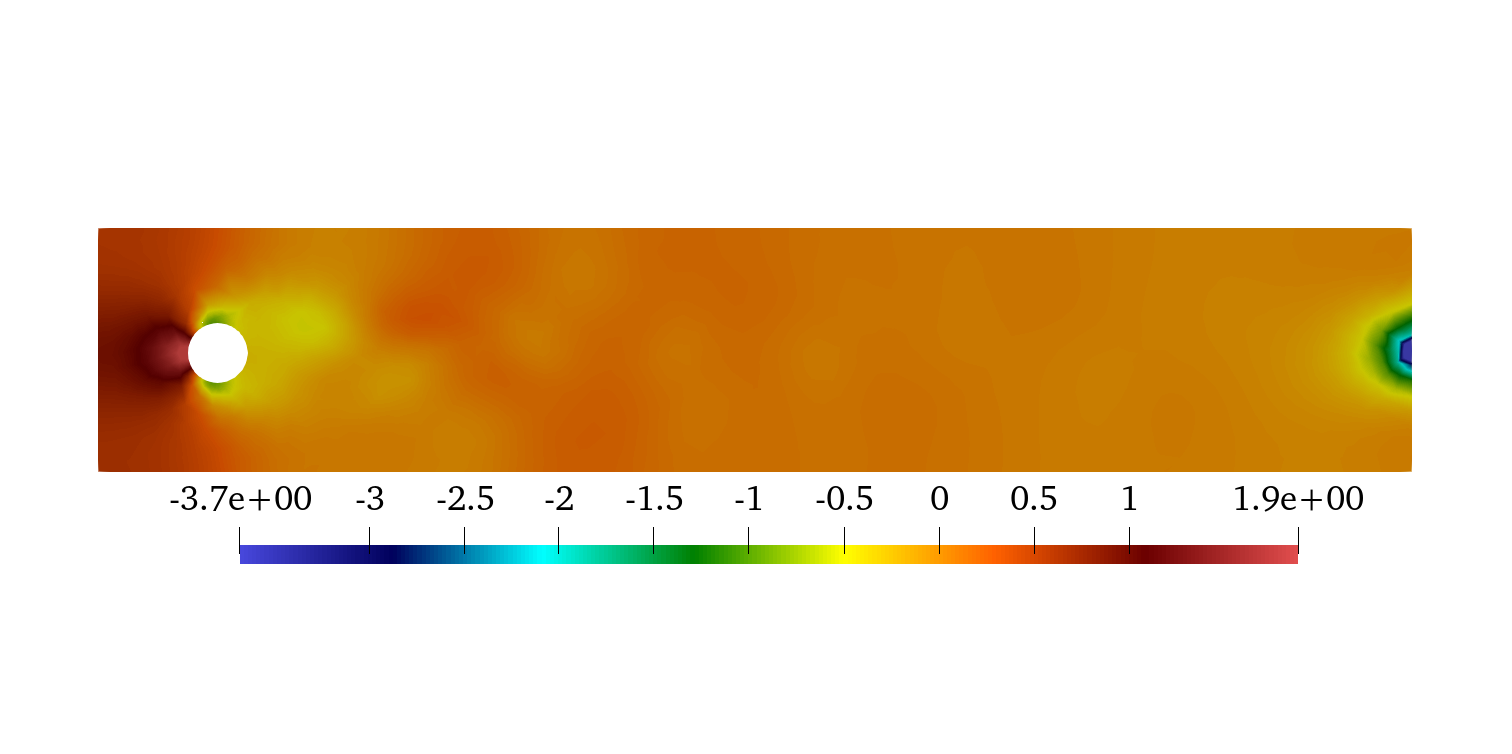}}
    \caption{From top to bottom: velocity (left plots) and pressure (right plots) fields at final $t=4$ for DNS 
    on the fine mesh (the reference solution), and noEFR method on the coarse mesh (i.e., the solution to improve with regularization). The results of standard EFR and EF approaches on the coarse mesh are also 
    displayed.}
    \label{fig:noefr}
\end{figure*}

{
\remark{In this project, we call \emph{DNS} the reference simulation considered to train the proposed data-driven algorithm. However, our reference cannot be considered a \emph{fully-resolved} simulation in its classical sense. Indeed, in Figure \ref{fig:noefr} plots (A) and (B) show that the velocity fields still exhibit spurious numerical oscillations.}}

\begin{figure*}[h!]
    \centering
    \subfloat[Coarse mesh]{\includegraphics[width=0.49\textwidth]{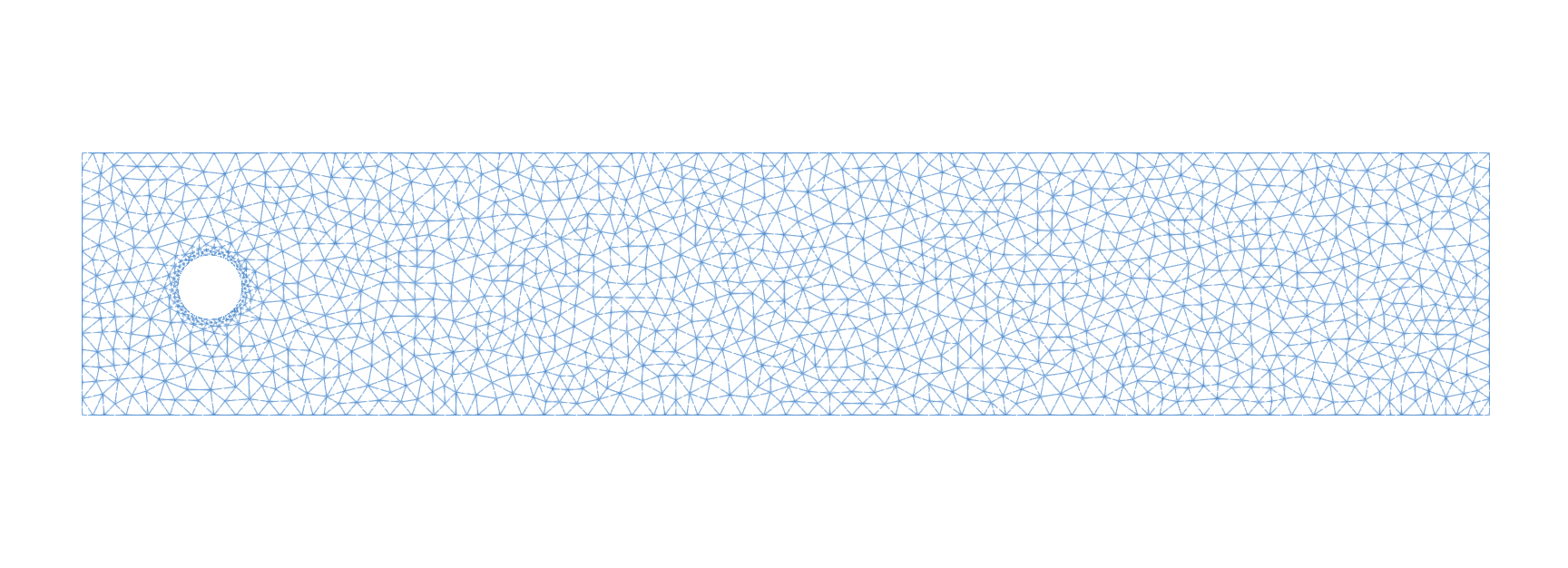}\label{mesh-coarse}}
    \subfloat[Fine mesh]{\includegraphics[width=0.49\textwidth]{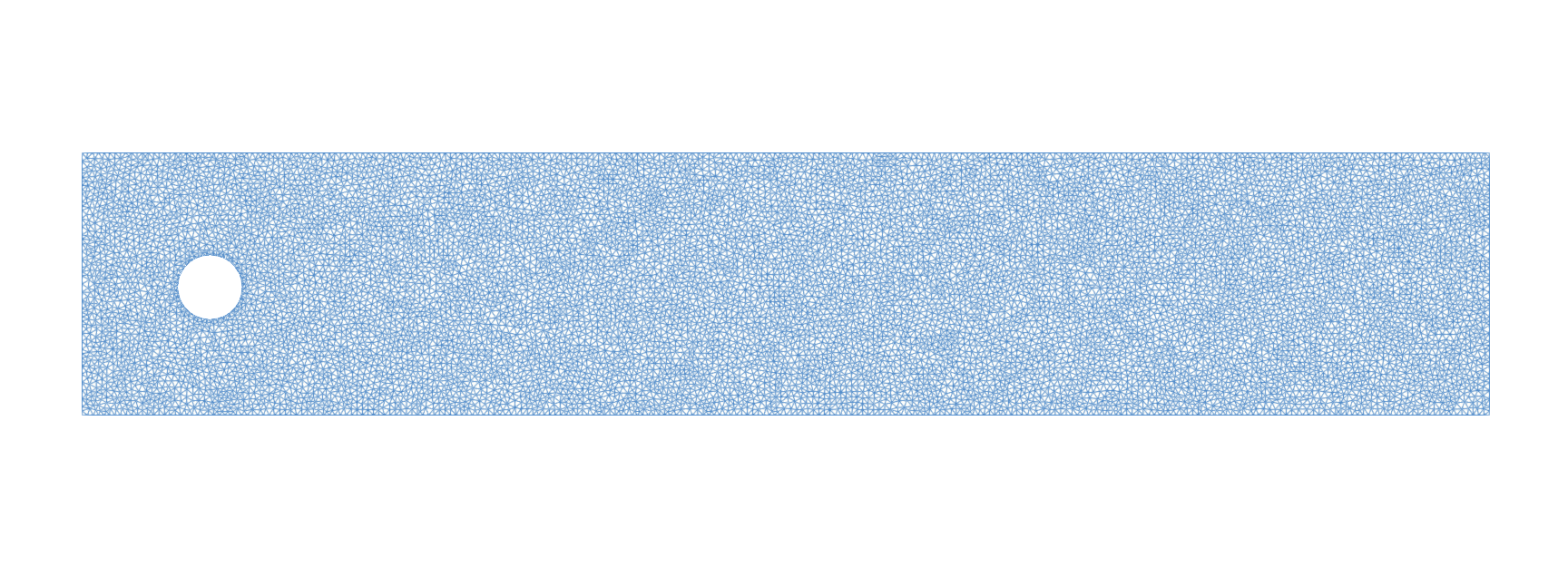}\label{mesh-ref}}
    \caption{Meshes taken into account to run the simulations.}
    \label{fig:meshes}
\end{figure*}

The two meshes taken into account are displayed in Figure \ref{fig:meshes}.
The coarse mesh is characterized by $N^{\bu}_h=12430$ dofs for the velocity field and $N^p_h=1623$ for the pressure.
On the other hand, the fine mesh is characterized by $M^{\bu}_h=100676$ and $M^p_h=12760$.
We perform the simulations in the time interval $[0, 4]$, and the adaptive optimization algorithms described in Section \ref{sec:methods-opt} are tested in the reconstructive regime.

The results obtained by the EFR simulations are compared with the DNS results using the following error metrics:
\begin{itemize}
    \item the $L^2$ relative errors 
    \begin{equation}
        E_{L^2}^u(t_n)=\dfrac{\|\bu^n - \bu_{ref}^n\|_{L^2(\Omega)}}{\| \bu_{ref}^n\|_{L^2(\Omega)}}, \quad
    E_{L^2}^p(t_n)=\dfrac{\|p^n - p_{ref}^n\|_{L^2(\Omega)}}{\| p_{ref}^n\|_{L^2(\Omega)}},
    \label{eq:errs-L2}
    \end{equation}
    \item the velocity relative error in the $H^1$ norm:
    \begin{equation}
E_{H^1}^u(t_n) = \dfrac{\|\bu^n - \bu^n_{ref}\|_{H^1(\Omega)}}{\|\bu_{ref}^n\|_{H^1(\Omega)}}.
    \label{eq:errs-h1}
\end{equation}
\end{itemize}

\begin{figure}[h!]
    \centering
\includegraphics[width=0.5\textwidth]{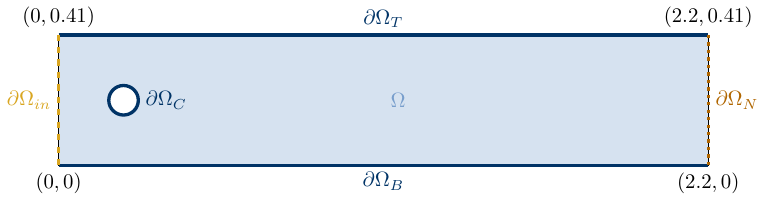}
    \caption{The domain considered in the test cases.}
    \label{fig:domain}
\end{figure}

The domain is $\Omega \doteq \{(x, y) \in [0, 2.2] \times [0, 0.41] \text{ such that } (x-0.2)^2 +(y-0.2)^2 \geq 0.05^2\}$, and it is represented in Figure \ref{fig:domain}. We consider no-slip boundary conditions on $\partial \Omega_{wall}=\partial \Omega_B \cup \partial \Omega_T \cup \partial \Omega_C$ (blue solid line in Figure \ref{fig:domain}), and 
an inlet velocity $\mathbf{u}_{in}$ on $\partial \Omega_{in}$ (gold dashed line in Figure \ref{fig:domain}), given by:
\begin{equation}
    \bu_{in }\doteq \left( \dfrac{6}{0.41^2}y(0.41-y), 0
    \right)\,.
    \label{u_inlet}
\end{equation}

On the outlet of the domain $\partial \Omega_N$ (brown dotted line in Figure \ref{fig:domain}), we employ standard free-flow conditions.  The initial condition is $\bu_0=(0, 0)$. The viscosity is set to $\nu=\num{1e-4}$ and, consequently, the Reynolds number is $Re=\num{1000}$.

For the test case described above, not applying regularization yields inaccurate field representation. We recall that, for the sake of the presentation, when no regularization is applied, we say that noEFR is performed, following the acronyms introduced in  \cite{strazzullo2022consistency}.

Some preliminary results motivating our analysis are available in Figure \ref{fig:noefr}, which yields the following conclusions: 
\begin{itemize}
    \item The noEFR simulations present spurious oscillations, 
as can be seen in plots (C) and (D).
This motivates the use of the EFR algorithm, see, also, \cite{strazzullo2023new, strazzullo2024variational}. 
    \item The challenges related to this specific test case are clear even from the DNS simulations that, despite the use of a finer mesh, still present some numerical oscillations near the outlet domain, see plots (A) and (B). 
    \item Moreover, the EFR method with standard choices $\delta=\eta=\num{5.62e-4}$ and $\chi=5\Delta t$ yields
poor results, as well. Indeed:
\begin{itemize}
    \item The standard EF approach is overdiffusive and does not capture the vortex shedding in the region near the cylinder, as can be seen from plots (G) and (H).
    \item The standard EFR method yields spurious oscillations, meaning that the relaxation parameter $\chi$ is too small and should be carefully tuned, as shown in plots (E) and (F).
\end{itemize}
\end{itemize}

The above statements motivate the following analysis of parameter optimization, whose main goal is to improve the results of the noEFR 
and standard EFR approaches with adaptive optimization. In particular:
\begin{itemize}
    \item Section \ref{sec:results-chi-opt} shows the results for $\chi$-Opt-EFR;
    \item Section \ref{sec:results-delta-opt} shows the results for $\delta$-Opt-EF;
    \item Section \ref{sec:results-double-opt} shows the results for $\delta \chi$-Opt-EFR.
\end{itemize}

%% file: sections/results-Opt-EFR.tex
\begin{figure*}[htpb!]
    \centering    
    \subfloat[]{\includegraphics[width=0.9\textwidth]{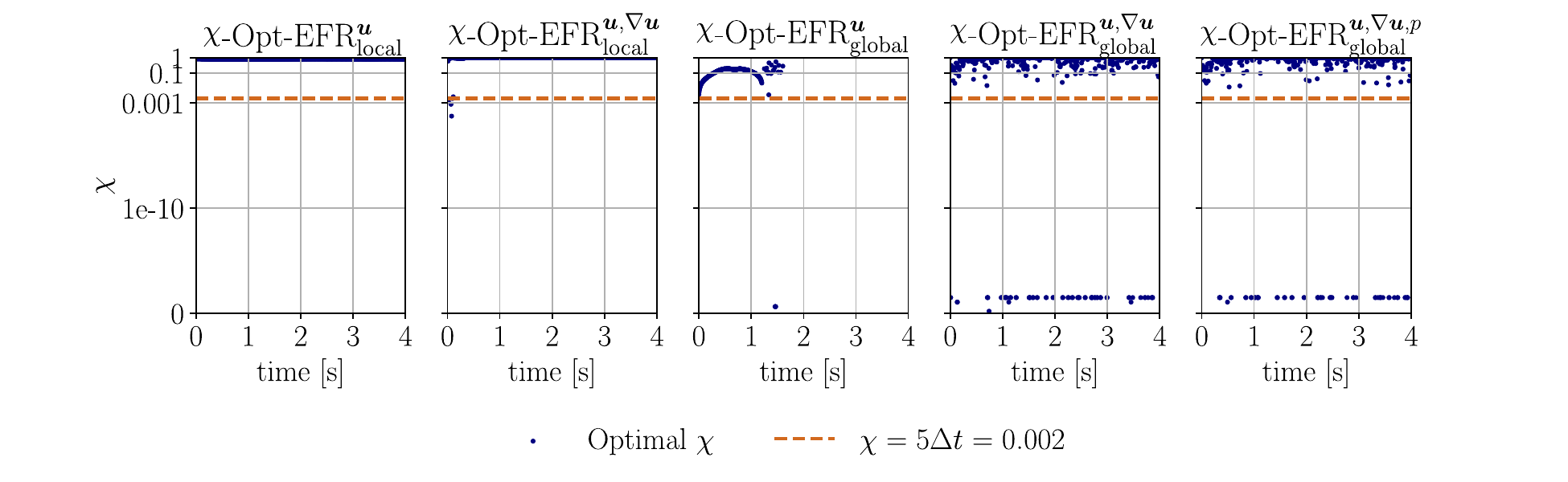}}\\
    \subfloat[]{\includegraphics[width=0.9\textwidth]{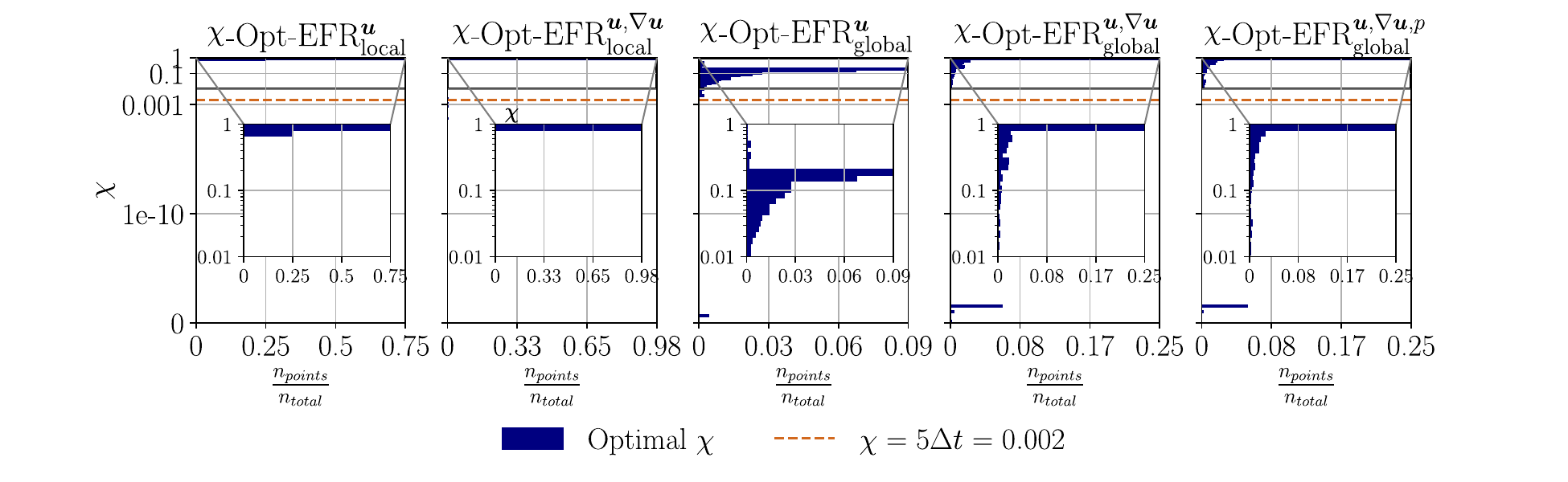}}
    \caption{Optimal value of $\chi(t)$ in the $\chi$-Opt-EFR algorithms (A), and count of observations of the parameter values (B).}
    \label{fig:chi}
\end{figure*}

First, 
in Figure \ref{fig:chi} (A) we include the $\chi(t)$ trends obtained with the different local or global objective functions taken into account.

When minimizing the \textbf{local} MSE discrepancy, the optimal $\chi(t)$ is close to 1 for each optimization time step. It means that the algorithm activates for almost 
$100 \%$ of the filtered velocity, and hence we expect an over-diffusize behavior, as happens in the standard EF approach.

\begin{figure*}[htpb!]
    \centering
     \subfloat[DNS - velocity ($t=4$)]{\includegraphics[width=0.5\textwidth, trim={3cm 7cm 3cm 7cm}, clip]{images/updated_Re1000/u_reference.png}}
    \\
    \subfloat[\optEFRloc{} - velocity ($t=4$)]{\includegraphics[width=0.5\textwidth, trim={3cm 7cm 3cm 7cm}, clip]{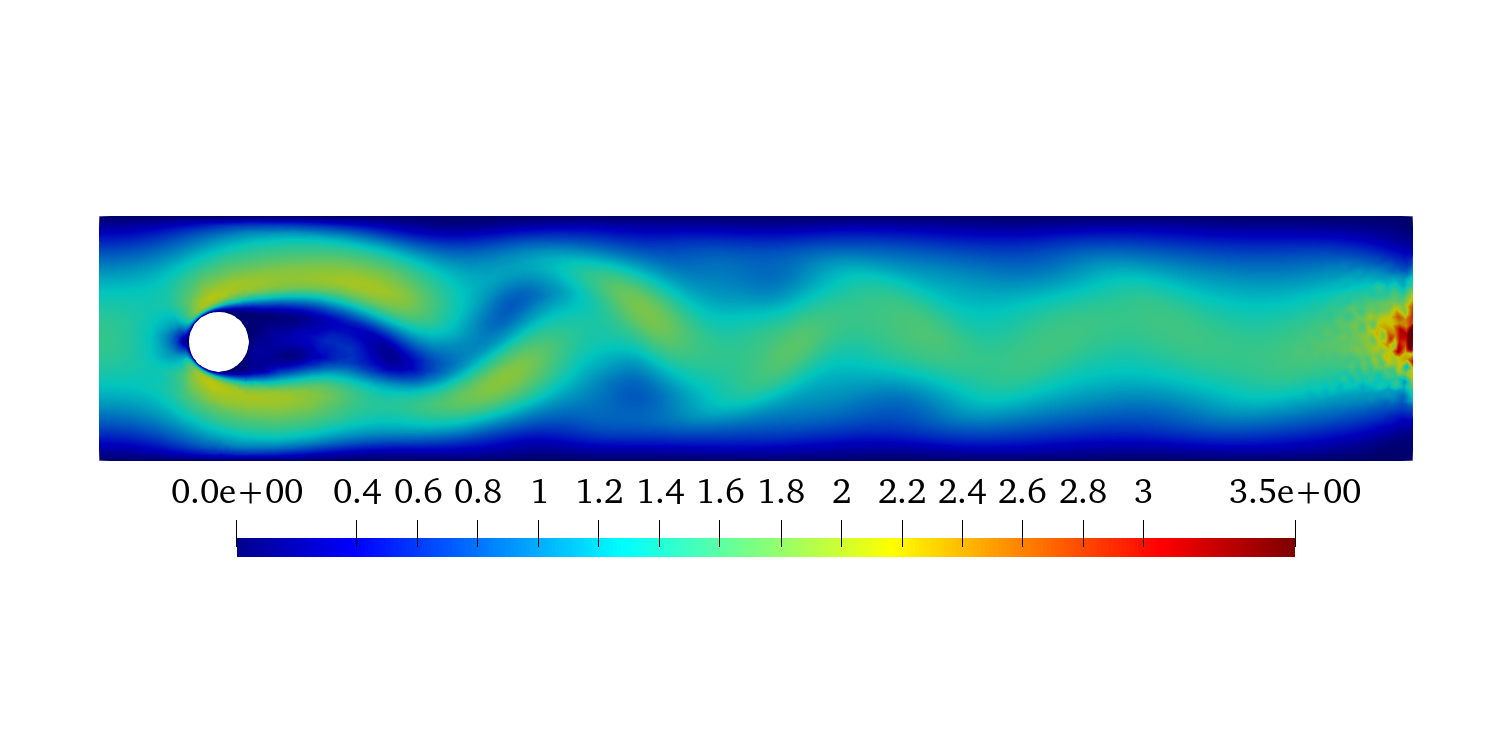}} 
    \subfloat[\optEFRlocGRAD{} - velocity ($t=4$)]{\includegraphics[width=0.5\textwidth, trim={3cm 7cm 3cm 7cm}, clip]{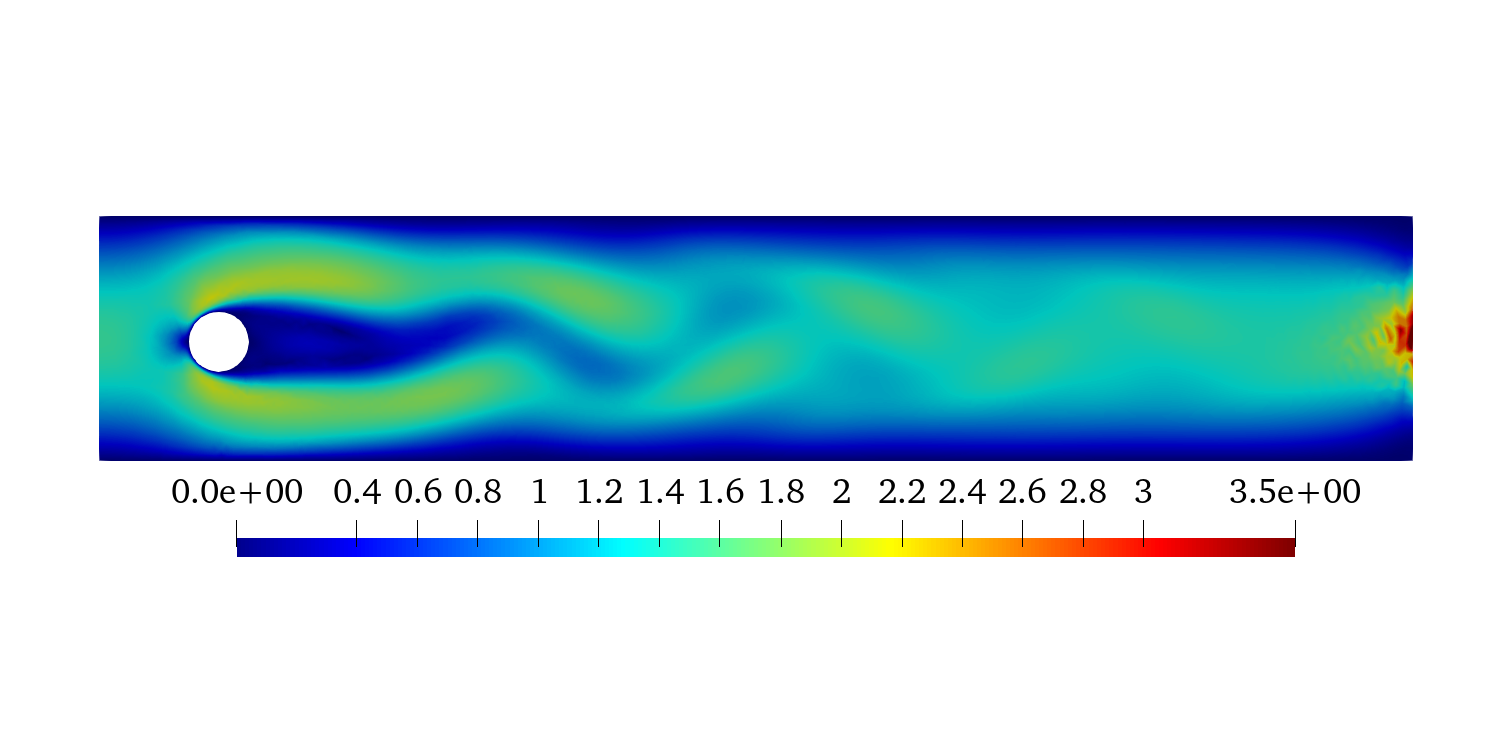}}\\
    \subfloat[\optEFRglob{} - velocity ($t=1.6$)]{\includegraphics[width=0.5\textwidth, trim={3cm 7cm 3cm 7cm}, clip]{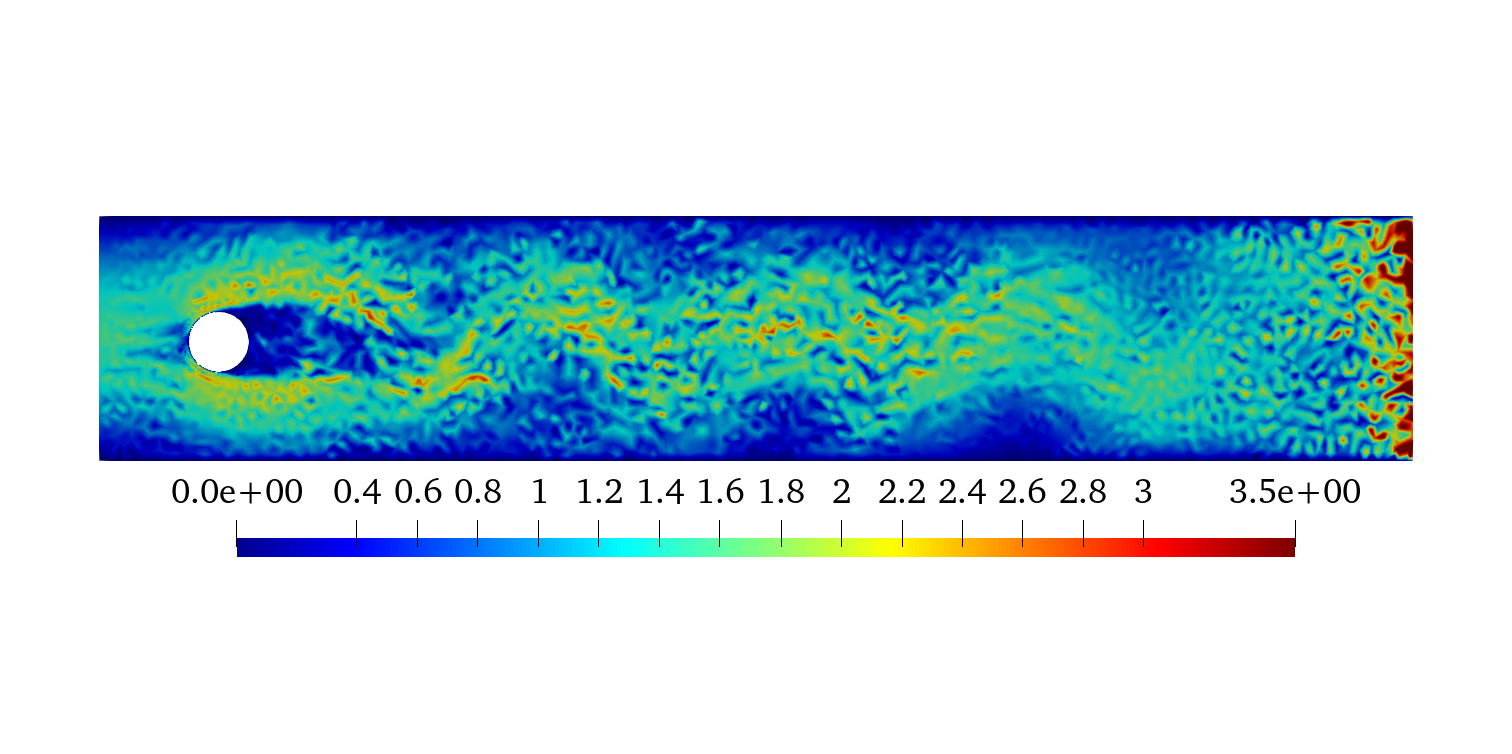}} 
    \subfloat[\optEFRglobGRAD{} - velocity ($t=4$)]{\includegraphics[width=0.5\textwidth, trim={3cm 7cm 3cm 7cm}, clip]{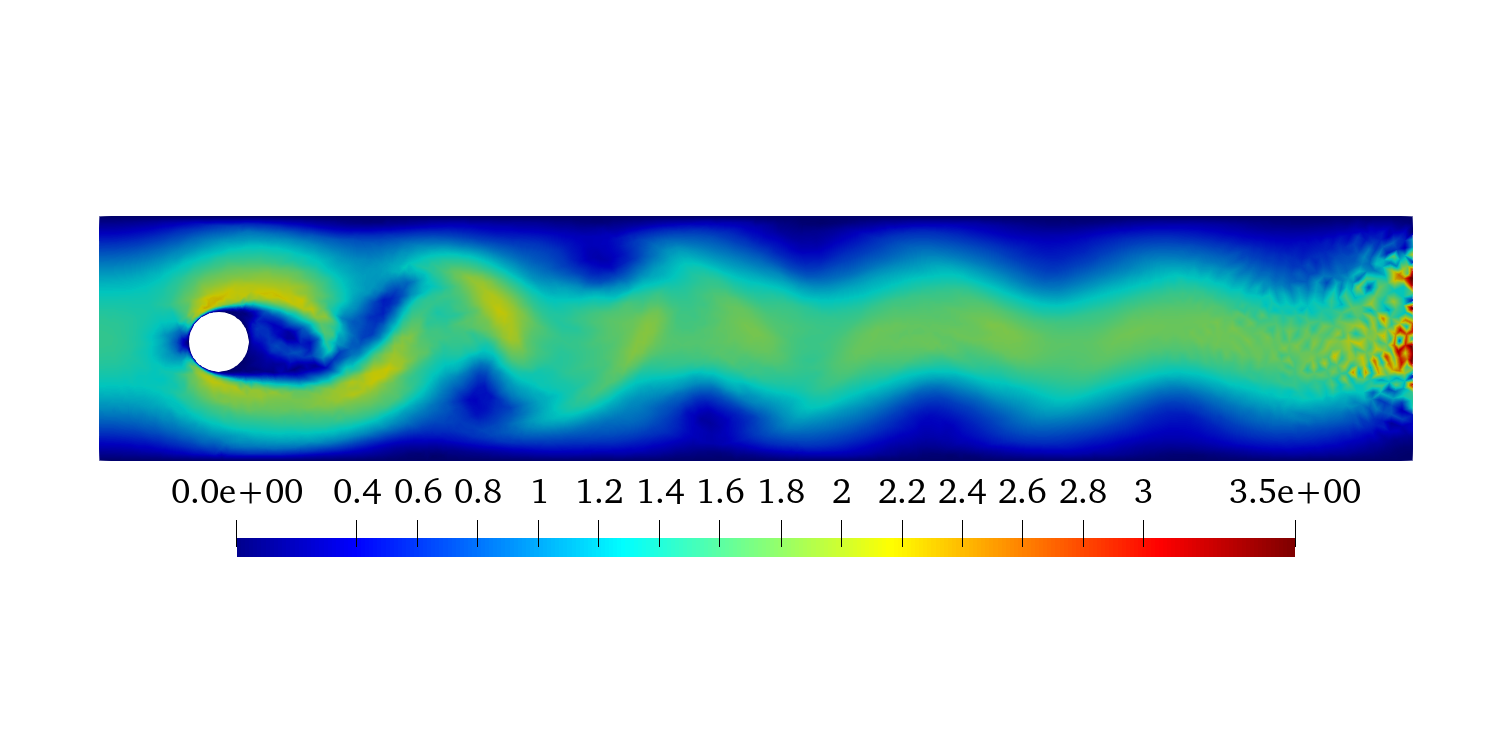}}
    \\
    \subfloat[\optEFRglobPRESS{} - velocity ($t=4$)]{\includegraphics[width=0.5\textwidth, trim={3cm 7cm 3cm 7cm}, clip]{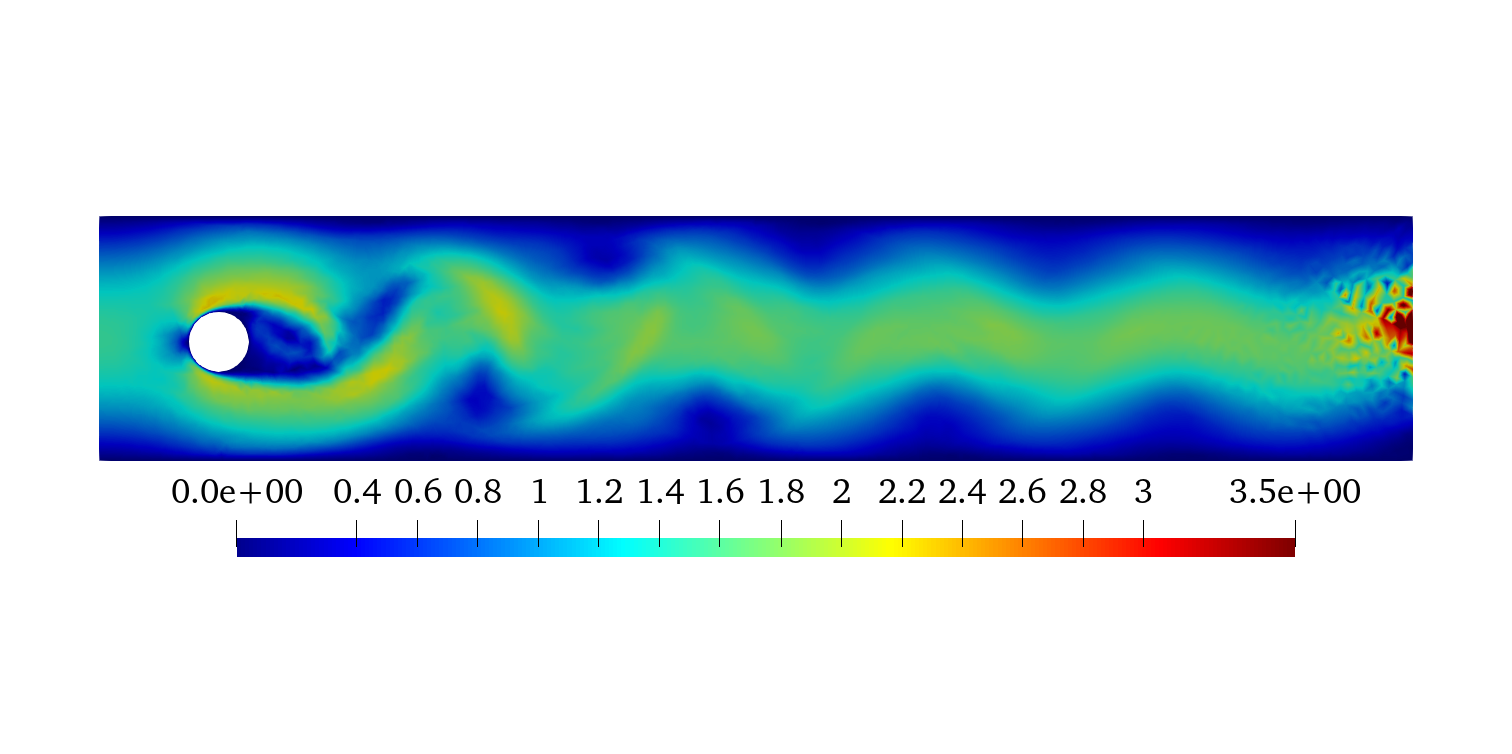}}
    \subfloat[Standard EFR ($\delta=\eta$, $\chi=5\Delta t$) - velocity ($t=4$)]{\includegraphics[width=0.5\textwidth, trim={3cm 6.5cm 3cm 7cm}, clip]{images/updated_Re1000/u_EFR_fixed.png}}
    \caption{Velocity fields at final time instance $t=4$ for standard EFR, $\chi$-Opt-EFR results, and 
    reference DNS simulation.
    }
    \label{fig:opt-efr-u}
\end{figure*}

\begin{figure*}[htpb!]
    \centering
    \subfloat[DNS - pressure ($t=4$)]{\includegraphics[width=0.5\textwidth, trim={3cm 7cm 3cm 7cm}, clip]{images/updated_Re1000/p_reference.png}}\\
    \subfloat[\optEFRloc{} - pressure ($t=4$)]{\includegraphics[width=0.5\textwidth, trim={3cm 7cm 3cm 7cm}, clip]{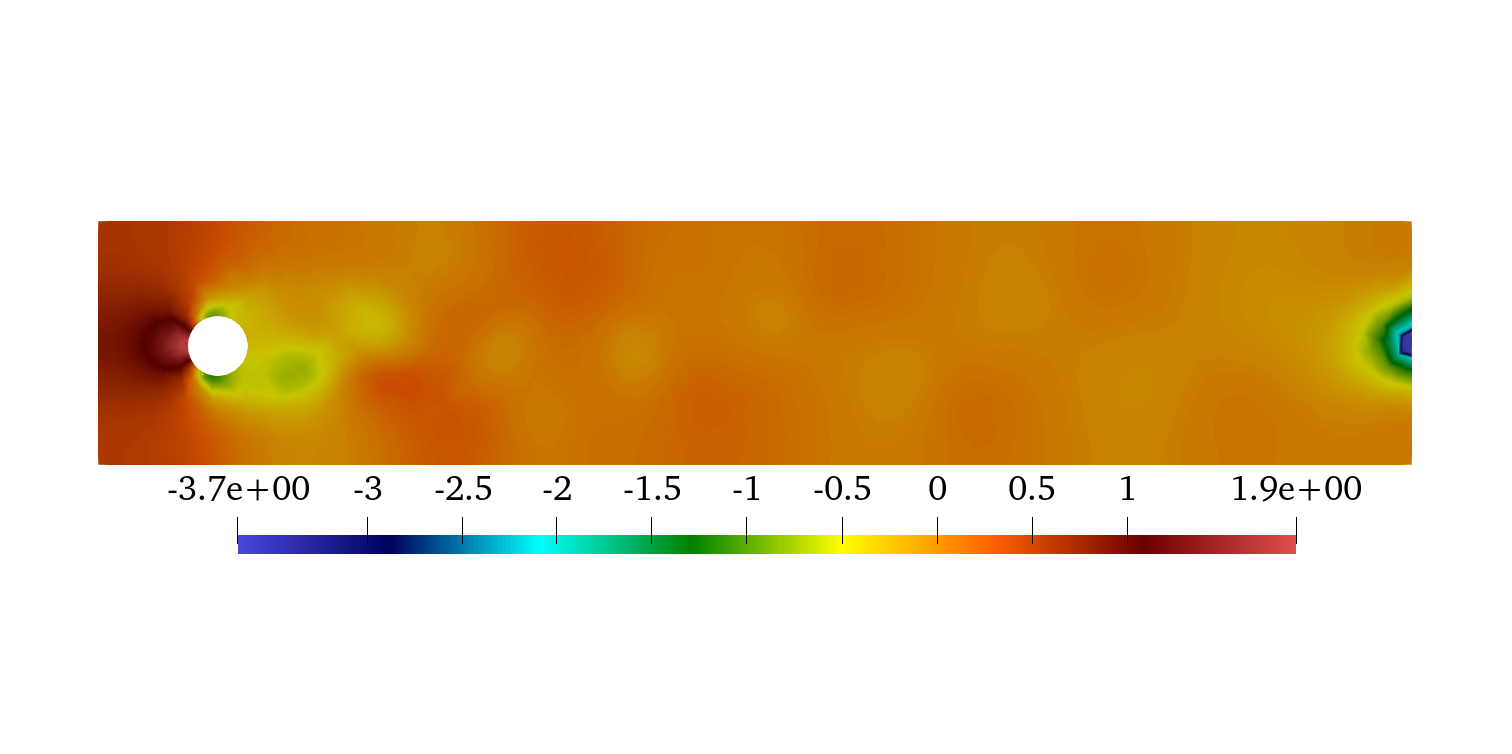}} 
    \subfloat[\optEFRlocGRAD{} - pressure ($t=4$)]{\includegraphics[width=0.5\textwidth, trim={3cm 6.5cm 3cm 7cm}, clip]{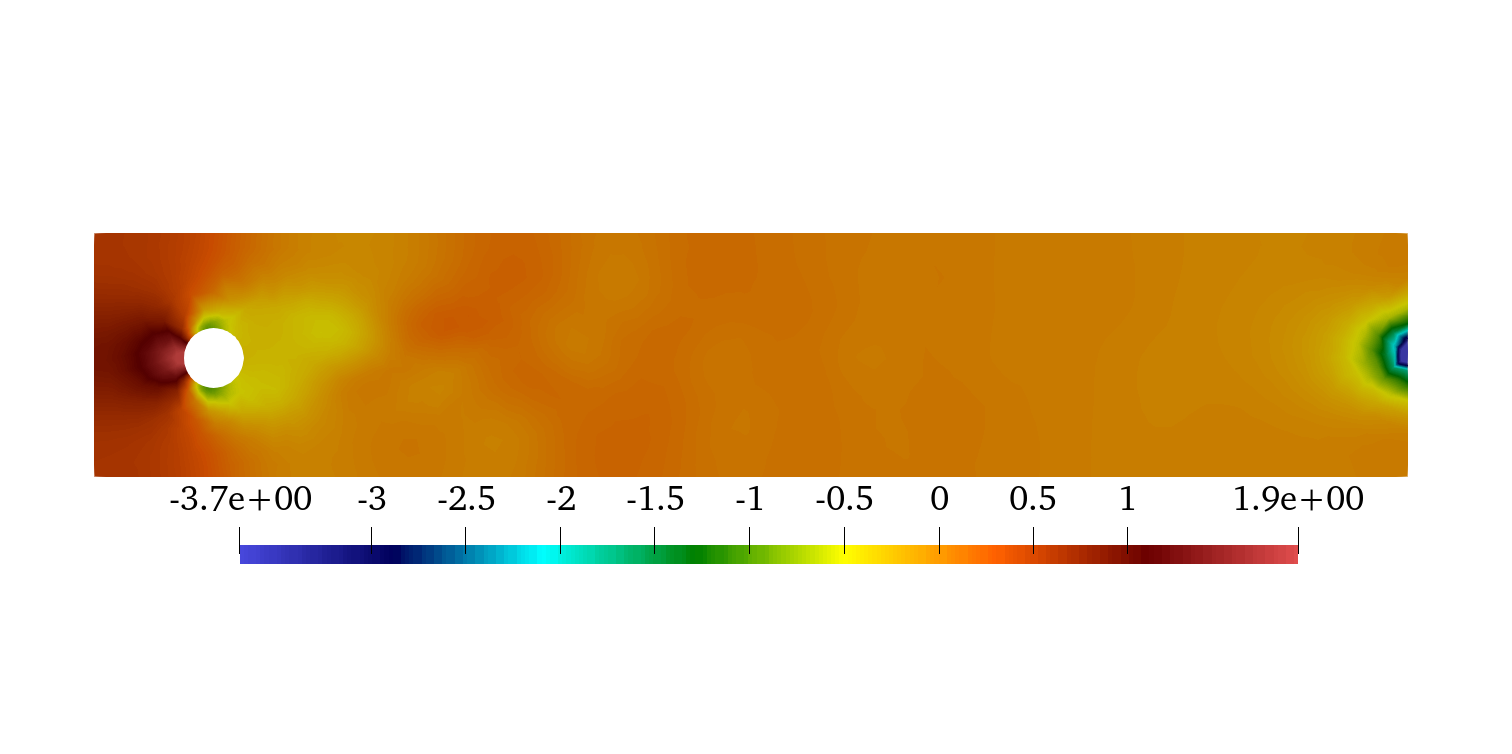}}\\
    \subfloat[\optEFRglob{} - pressure ($t=1.6$)]{\includegraphics[width=0.5\textwidth, trim={3cm 6.5cm 3cm 7cm}, clip]{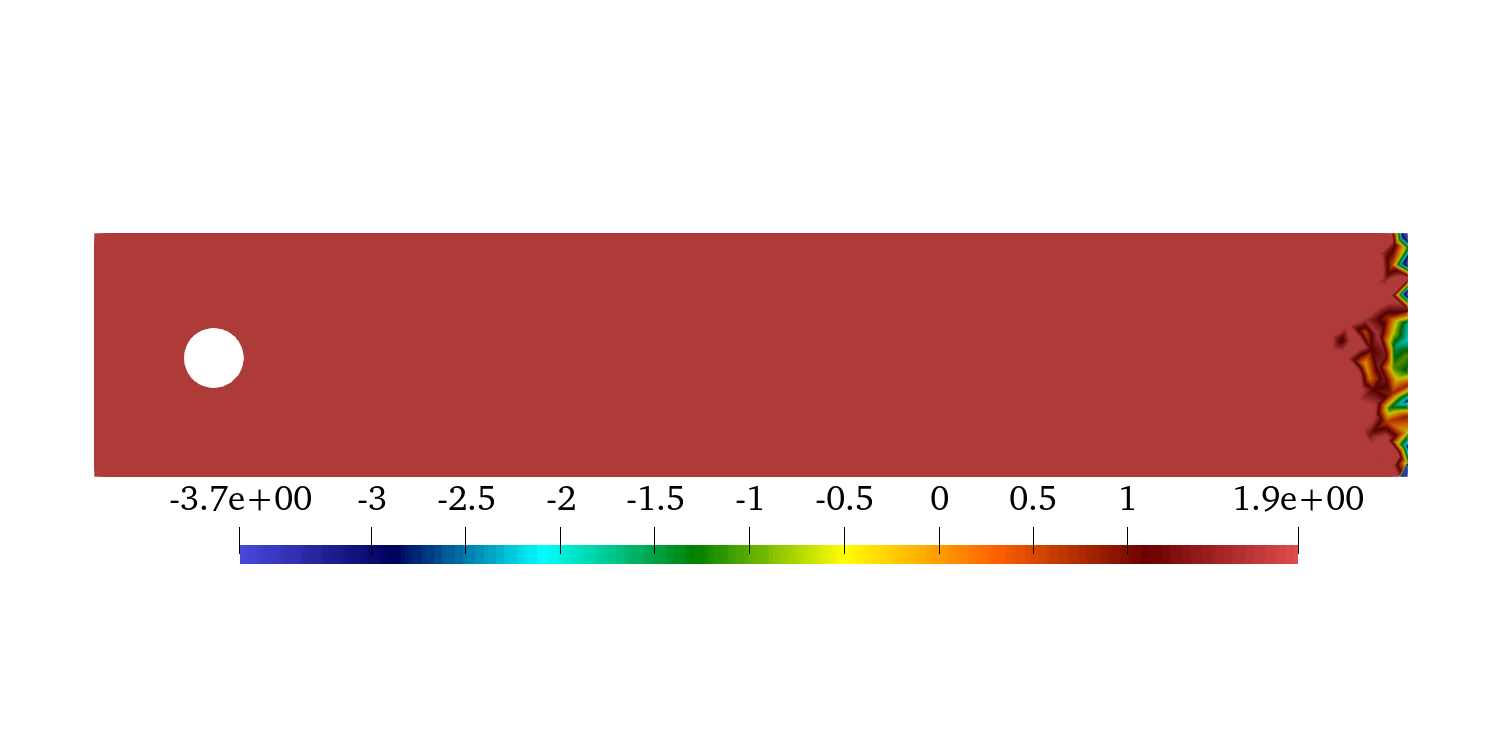}} 
    \subfloat[\optEFRglobGRAD{} - pressure ($t=4$)]{\includegraphics[width=0.5\textwidth, trim={3cm 6.5cm 3cm 7cm}, clip]{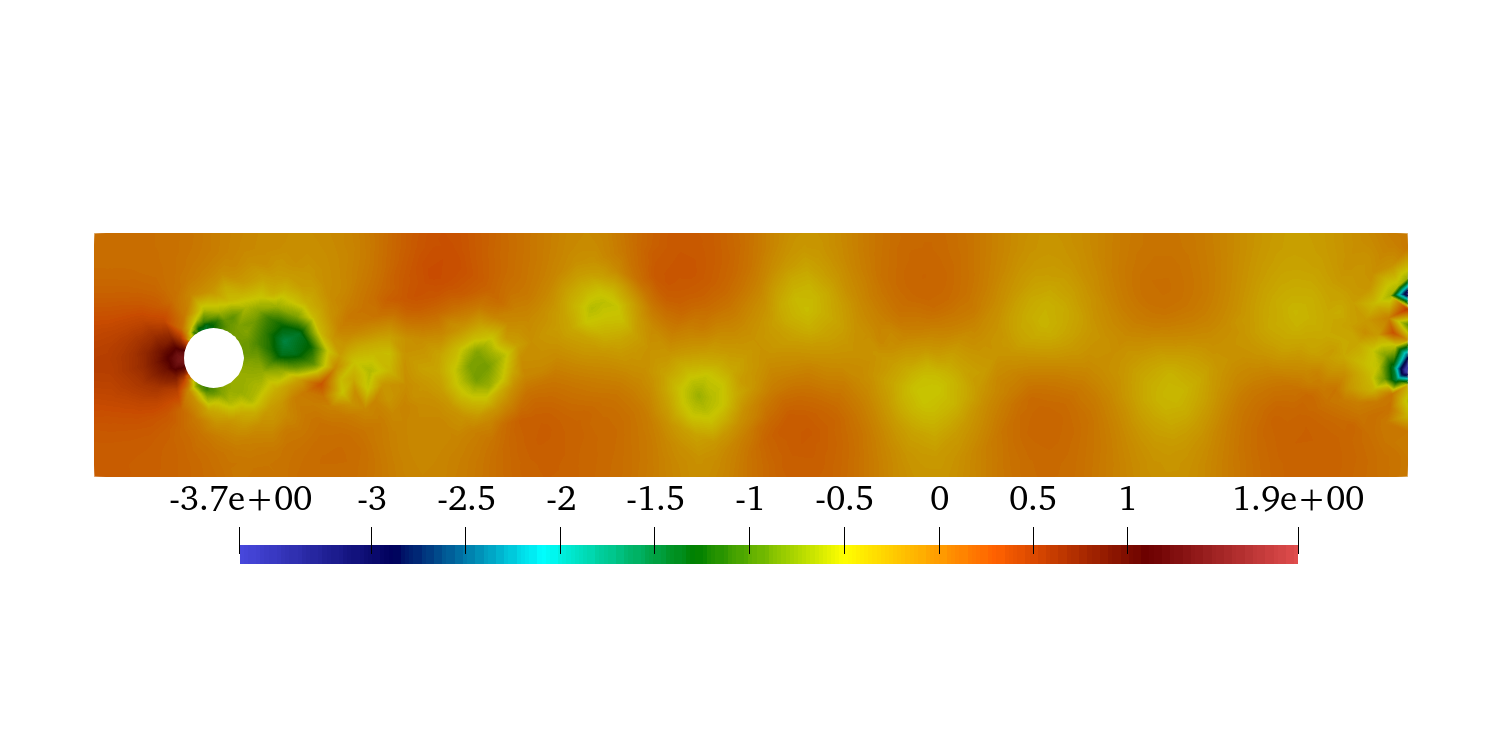}}
    \\
    \subfloat[\optEFRglobPRESS{} - pressure ($t=4$)]{\includegraphics[width=0.5\textwidth, trim={3cm 7cm 3cm 7cm}, clip]{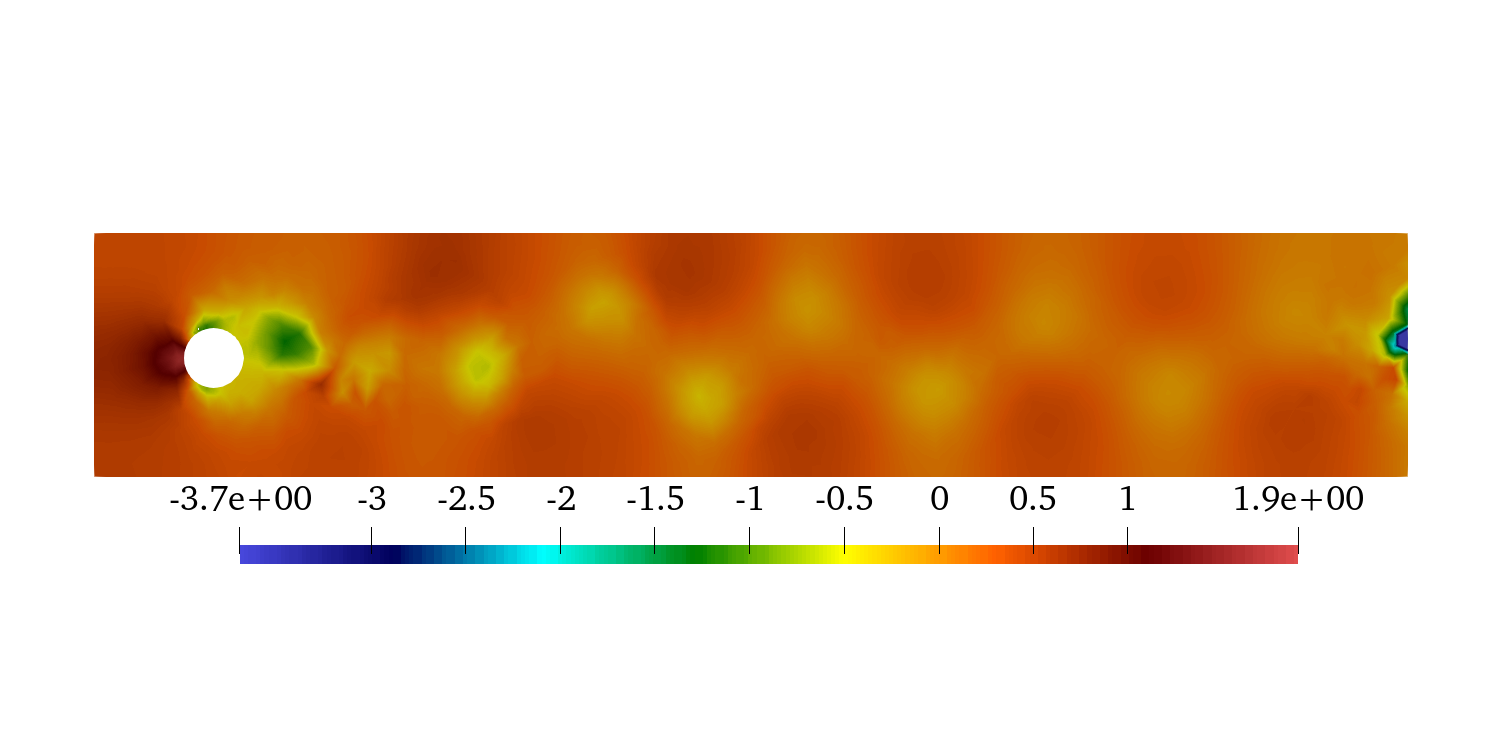}}
    \subfloat[Standard EFR ($\delta=\eta$, $\chi=5\Delta t$) - pressure ($t=4$)]{\includegraphics[width=0.5\textwidth, trim={3cm 6.5cm 3cm 7cm}, clip]{images/updated_Re1000/p_EFR_fixed.png}}
    \caption{Pressure fields at final time instance $t=4$ for Standard EFR, $\chi$-Opt-EFR results, and for the reference DNS simulation.}
    \label{fig:opt-efr-p}
\end{figure*}

On the other hand, when we only consider the \textbf{global} velocity norm discrepancy, the simulation blows up at $t=1.6$, as can be seen in Figure \ref{fig:errs-chi}, and in Figures \ref{fig:opt-efr-u} and \ref{fig:opt-efr-p} (D).

The simulation is more stable when 
the velocity gradients and the pressure contributions are also included into the objective function (approaches \optEFRglobGRAD{} and \optEFRglobPRESS{}).
However, the last two plots on the right in Figure \ref{fig:chi} (A) show an oscillating and chaotic behavior of $\chi(t)$.
To better understand the distribution of the parameters' values, we present in Figure \ref{fig:chi} (B) the relative 
number of optimal $\chi^n$ values attained by the
Opt-EFR algorithms, with $n=1, \dots, N_T$.
In particular, for a given value, $n_{points}$ is the number of times $\chi^n$ attains the given value,
while $n_{total}$ is the total number of samples, namely 1000. Figure \ref{fig:chi} shows us that the number of $\chi(t)$ elements close to $1$ is around 
$25\%$ in both algorithms \optEFRglobGRAD{} and \optEFRglobPRESS{}. This means the optimization algorithm converges statistically more often to values close to $1$, deactivating the relaxation step.

\RA{A comparison among the methods with respect to the $L^2$ relative errors, i.e., $E_{L^2}^u(t)$ and $E_{L^2}^p(t)$, is not included here for the sake of brevity. The reader can find this analysis in the supplementary material (section \ref{subsec:supp-results-chi}).}

\begin{figure*}[htpb!]
    \centering\includegraphics[width=0.9\textwidth]{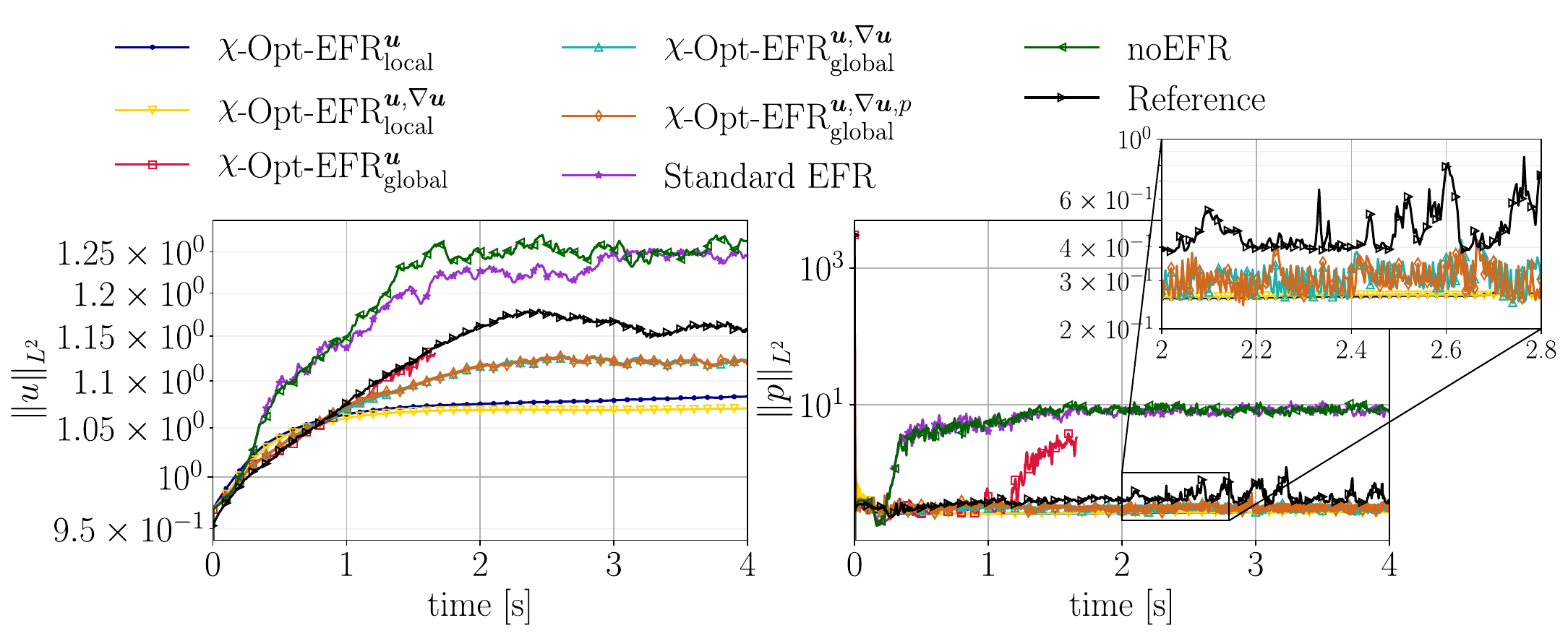}
    \caption{Velocity and pressure norms of $\chi$-Opt-EFR simulations, of the reference solutions projected on the coarse mesh, and of standard EFR. \RB{The pressure plot (right) also includes a box with a zoomed-in area in the time interval $[2, 2.8]$.}}
    \label{fig:norms-chi}
\end{figure*}

\begin{figure}[htpb!]
    \centering
    \includegraphics[width=0.5\textwidth]{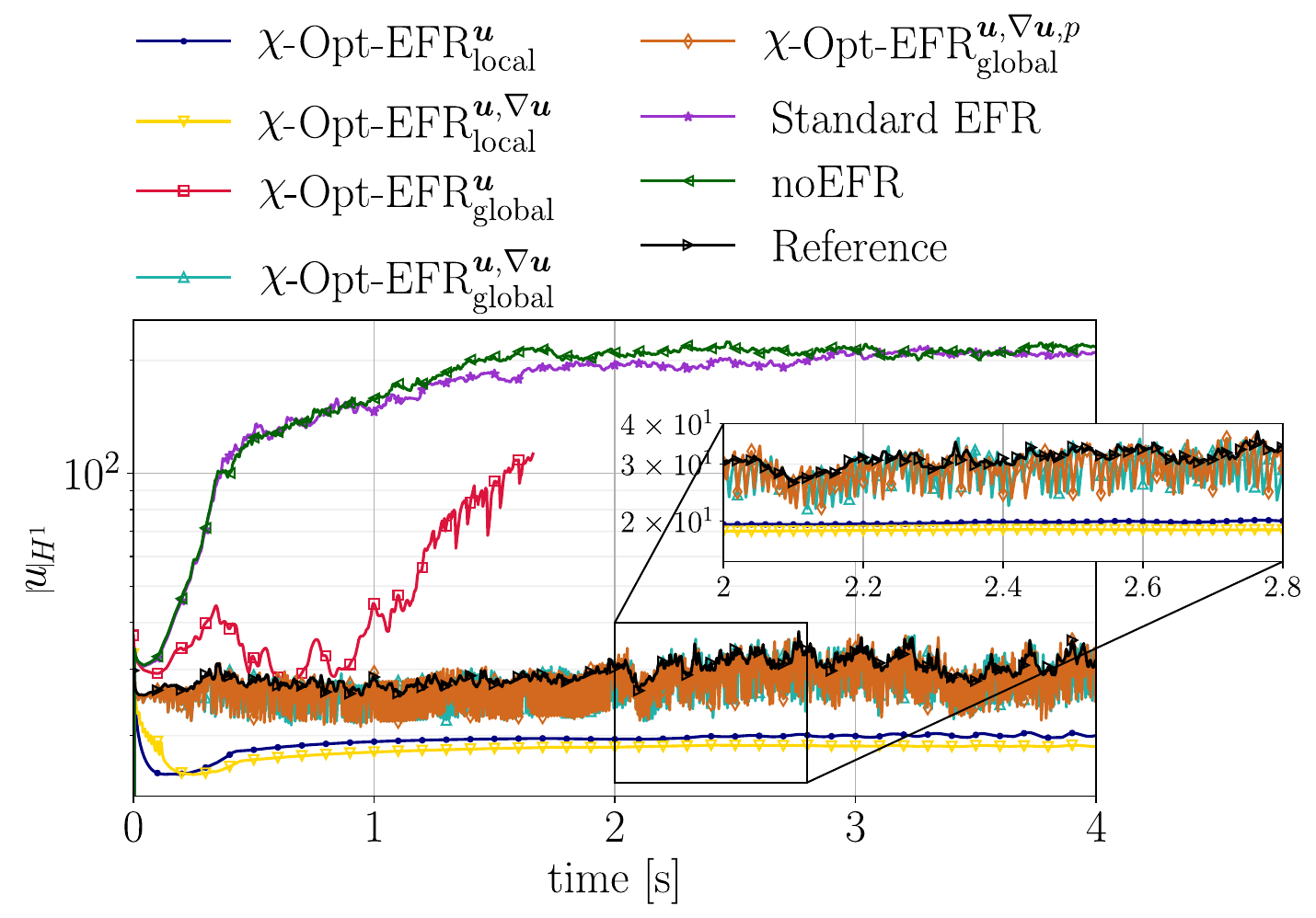}
    \caption{Velocity $H^1$ norms of $\chi$-Opt-EFR simulations, of the reference solutions, projected on the coarse mesh, and of standard EFR. \RB{The plot also includes a box with a zoomed-in area in the time range $[2, 2.8]$.}}
    \label{fig:norm-chi-h1}
\end{figure}

The $L^2$ velocity and pressure norms are displayed in Figure \ref{fig:norms-chi}, whereas Figure \ref{fig:norm-chi-h1} displays the velocity $H^1$ seminorms. We recall that the objective functions in algorithms \optEFRglob{}, \optEFRglobGRAD{}, and \optEFRglobPRESS{} include the information on this specific metric. Hence, we expect better results in these cases with respect to local optimization approaches.
As expected, Figures \ref{fig:norms-chi} and \ref{fig:norm-chi-h1} show that:

\begin{itemize}
    \item \optEFRloc{} and \optEFRlocGRAD{} simulations produce less accurate approximations than in the cases of global optimizations;
    \item \optEFRglob{} shows a good agreement in the $L^2$ velocity norm (until $t=1.6$)
with respect to the reference counterpart because it is exactly the metric optimized in the optimization algorithm, while the pressure norm and the velocity gradients are not accurately predicted. Indeed, this approach yields a poor approximation of the velocity and pressure fields, as can be seen in Figures \ref{fig:opt-efr-u} and \ref{fig:opt-efr-p} (D);
 \item Approaches \optEFRglobGRAD{} and \optEFRglobPRESS{} provide the best results, and show a good agreement in the $H^1$ velocity seminorm and in the pressure norm. \RB{This is confirmed by the zoomed-in area in Figure \ref{fig:norm-chi-h1}, which highlights the oscillatory behavior of the Opt-EFR norms. This behavior seems to reproduce the DNS reference norm, while the standard EFR is overdiffusive, and, hence, has a constant norm in time};
 \item We expect that the addition of the pressure contribution in \optEFRglobPRESS{} would be beneficial and improve the pressure accuracy with respect to \optEFRglobGRAD{}. However, in both cases, the pressure norm has a chaotic trend in time. The graphical comparison between Figure \ref{fig:opt-efr-p} (E) and (F) shows that the addition of the pressure discrepancy in the objective function only slightly increases the pressure accuracy.
\end{itemize}

\remark{ 
To highlight the reason why the procedure \optEFRglob{} blows up at $t=1.6$, in Figure \ref{fig:norms-chi} we display the pressure norm, which starts increasing after $t=1$, until it completely explodes, as can be seen in Figure \ref{fig:opt-efr-p} in subplot (D). This happens because the optimization algorithm for parameter $\chi$ does not converge to a global optimum.
}

%% file: sections/results-Opt-EF.tex
\begin{figure*}[htpb!]
    \centering
    \subfloat[]{\includegraphics[width=0.75\textwidth]{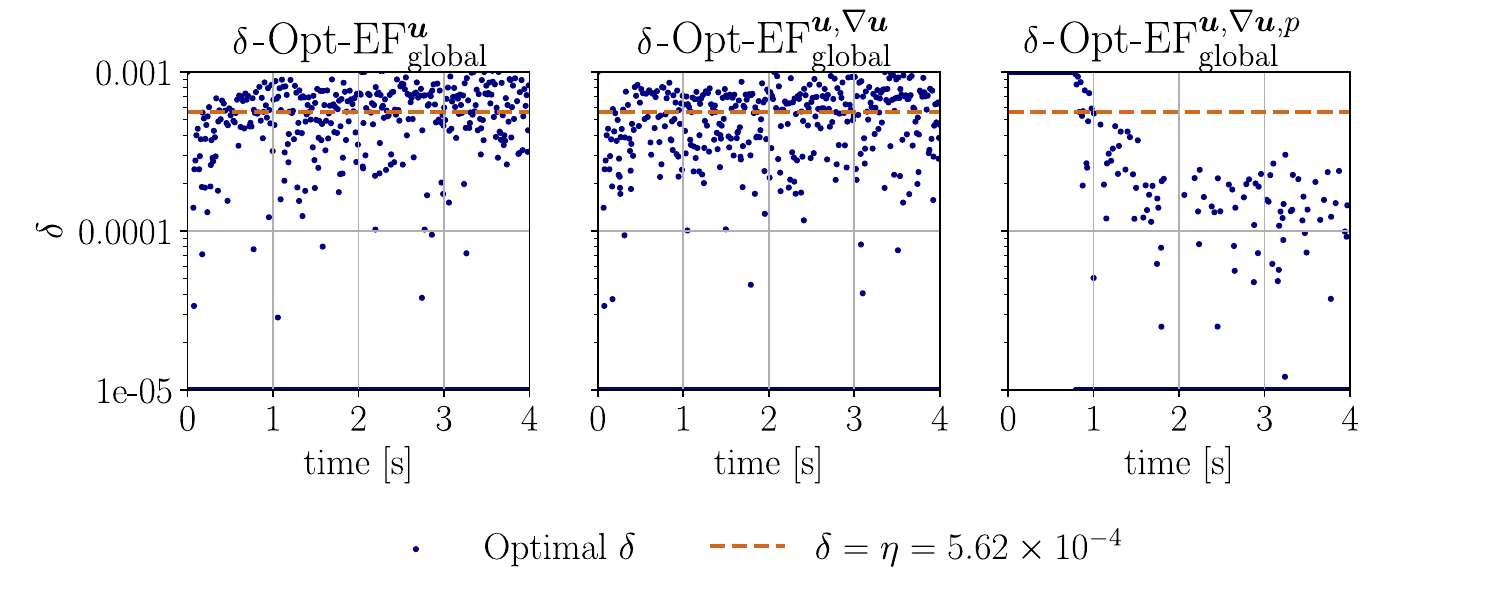}}\\
    \subfloat[]{\includegraphics[width=0.75\textwidth]{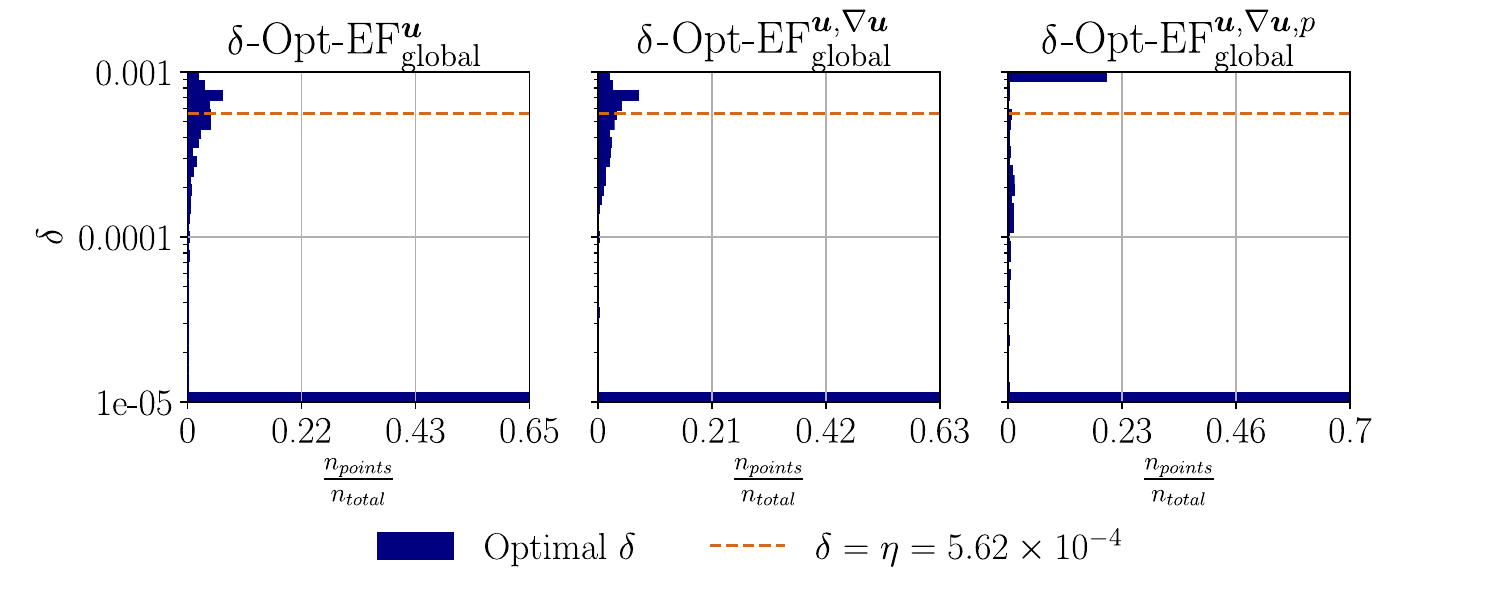}}
    \caption{Optimal value of $\delta(t)$ in the $\delta$-Opt-EF algorithms (A), and count of observations of the parameter values (B).}
    \label{fig:delta}
\end{figure*}

In this section, we present
results of the $\delta$-Opt-EF algorithms, namely the ones optimizing parameter $\delta(t)$, for a fixed 
relaxation parameter, 
$\chi=1$.
Figure \ref{fig:delta} (A) shows the trend of the optimized filter parameter in time.

The results here are similar to the results obtained in Subsection \ref{sec:results-chi-opt} with the global optimization of parameter $\chi(t)$, namely 
$\delta(t)$ 
is oscillating, but in this case in the prescribed range $[\num{1e-5}, \num{1e-3}]$. In particular, in the \optEFglobPRESS{} algorithm, the optimal parameter converges to $\num{1e-3}$ in the first time steps, and then it oscillates but most of the time converges to $\num{1e-5}$.
Figure \ref{fig:delta} (B) presents the distribution of the parameter 
values, and it shows that the most frequently occurring value is $\num{1e-5}$ in all cases. We note that the oscillating behavior between this minimum value and higher values is necessary to alleviate the spurious oscillations appearing in the solution.

\begin{figure*}[htpb!]
    \centering
    \includegraphics[width=0.9\textwidth]{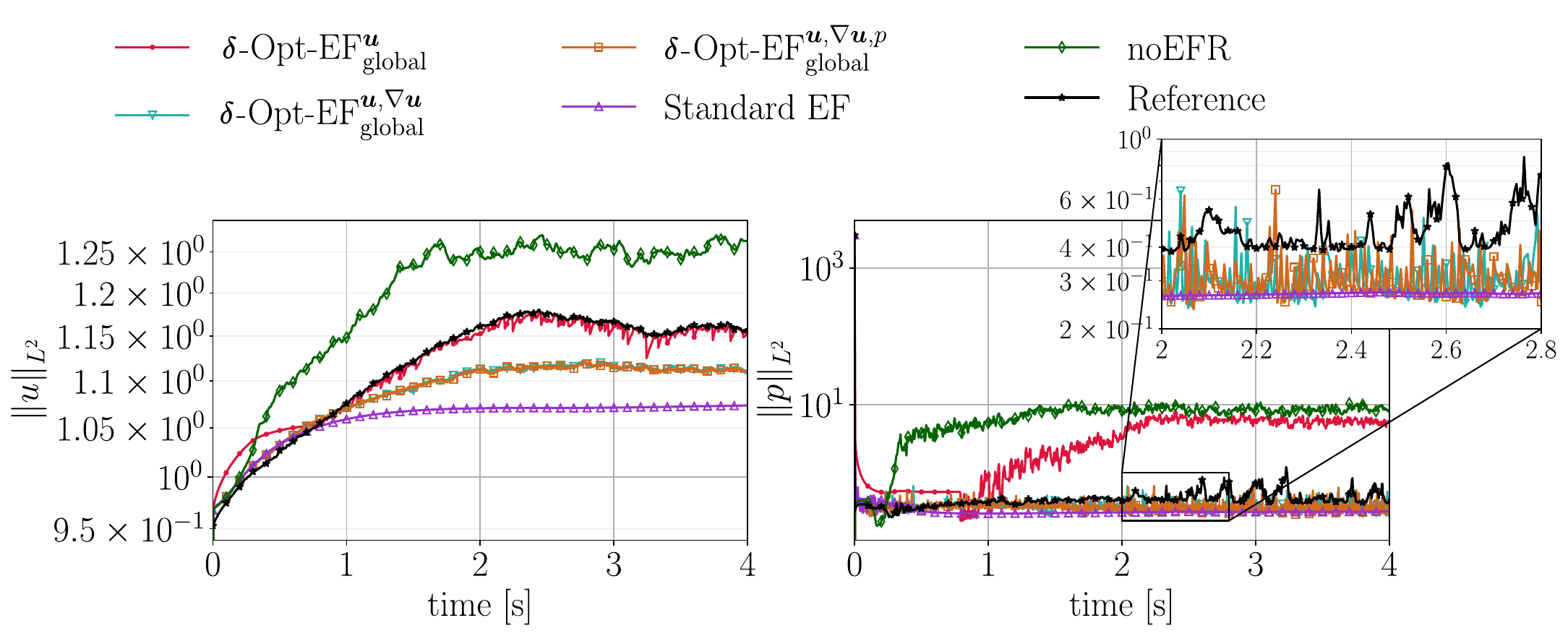}
    \caption{Velocity and pressure norms of $\delta$-Opt-EF simulations, of the EF simulation with fixed $\delta$, and of the reference solutions, projected on the coarse mesh. \RB{The pressure plot (right) also includes a box with a zoomed-in area in the time interval $[2, 2.8]$.}}
    \label{fig:norms-delta}
\end{figure*}

\begin{figure}[htpb!]
    \centering\includegraphics[width=0.5\textwidth]{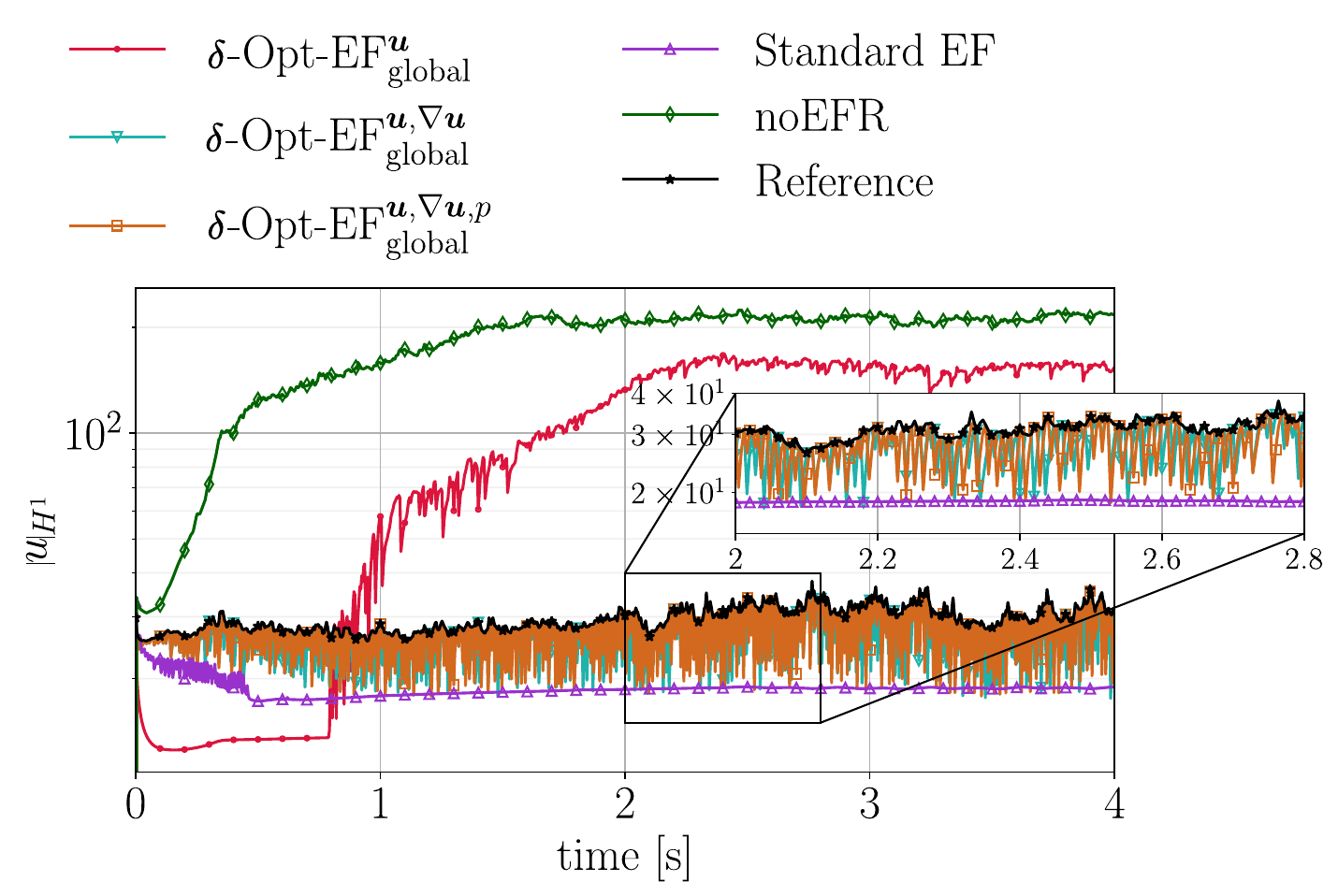}
    \caption{Velocity $H^1$ norms of $\delta$-Opt-EF simulations, of the EF simulation with fixed $\delta$, and of the reference solutions, projected on the coarse mesh. \RB{The plot also includes a box with a zoomed-in area in the time interval $[2, 2.8]$.}}
    \label{fig:norm-delta-h1}
\end{figure}

\begin{figure*}[htpb!]
    \centering
    \subfloat[DNS - velocity ($t=4$)]
    {\includegraphics[width=0.5\textwidth, trim={3cm 7cm 3cm 7cm}, clip]{images/updated_Re1000/u_reference.png}}\\
    \subfloat[\optEFglob{} - velocity ($t=4$)]{\includegraphics[width=0.5\textwidth, trim={3cm 7cm 3cm 7cm}, clip]{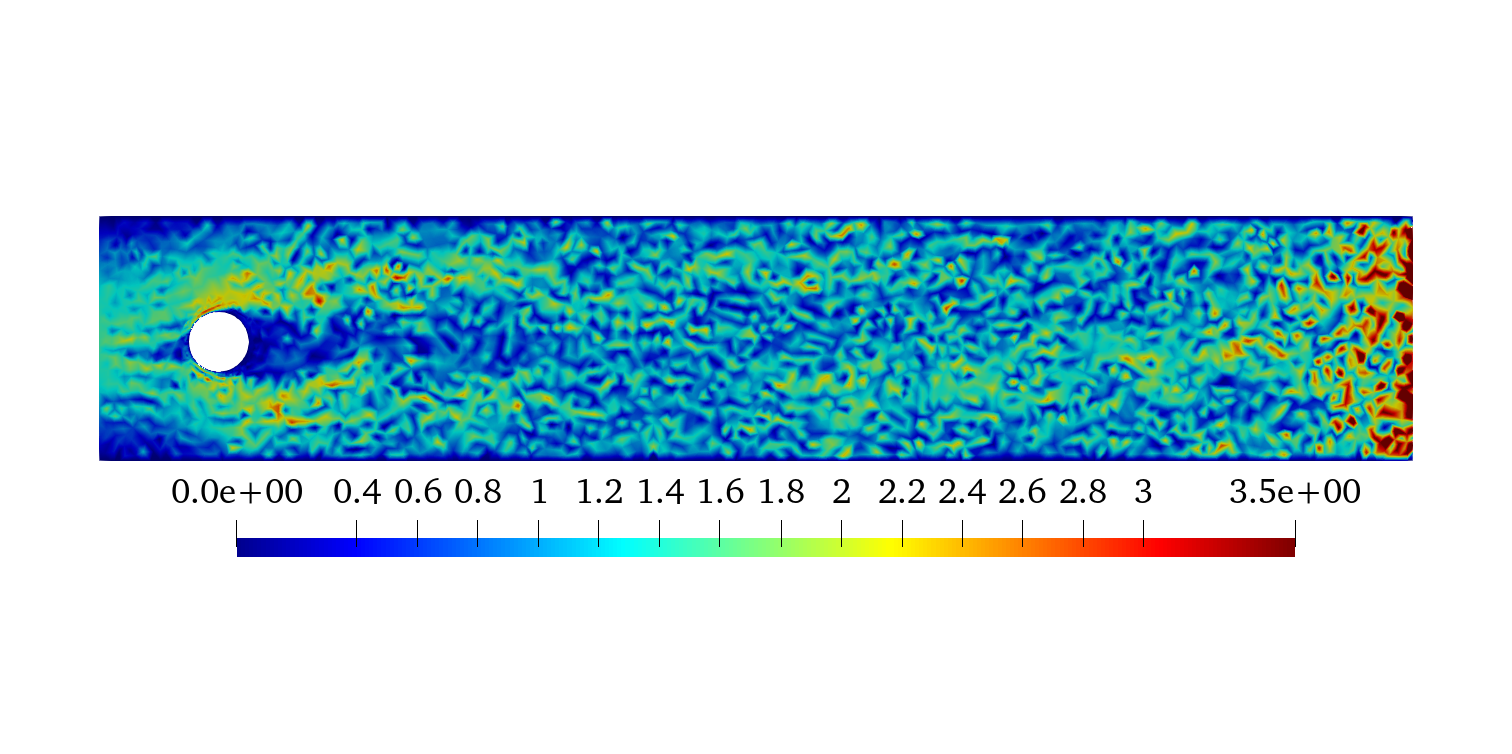}} 
    \subfloat[\optEFglobGRAD{} - velocity ($t=4$)]{\includegraphics[width=0.5\textwidth, trim={3cm 7cm 3cm 7cm}, clip]{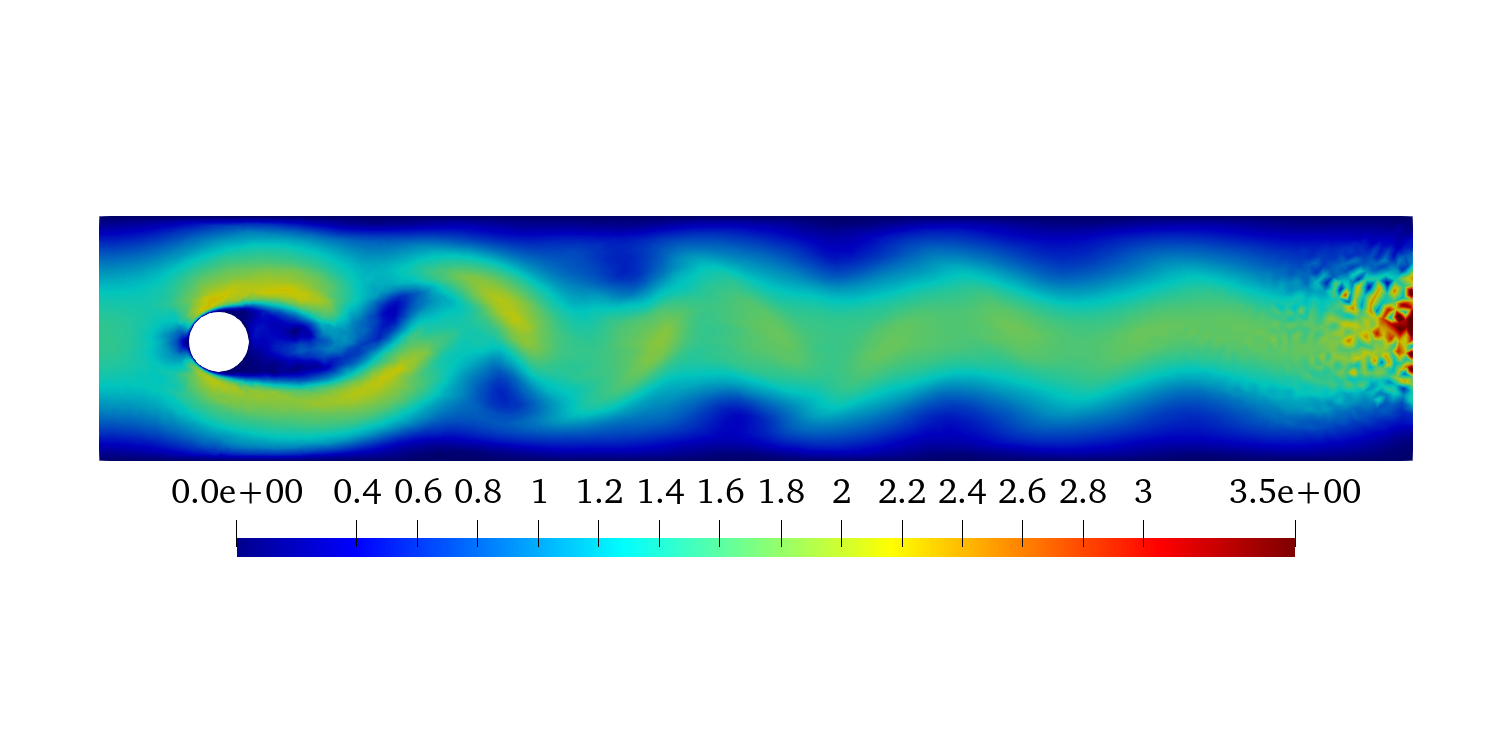}}\\
    \subfloat[\optEFglobPRESS{} - velocity ($t=4$)]{\includegraphics[width=0.5\textwidth, trim={3cm 7cm 3cm 7cm}, clip]{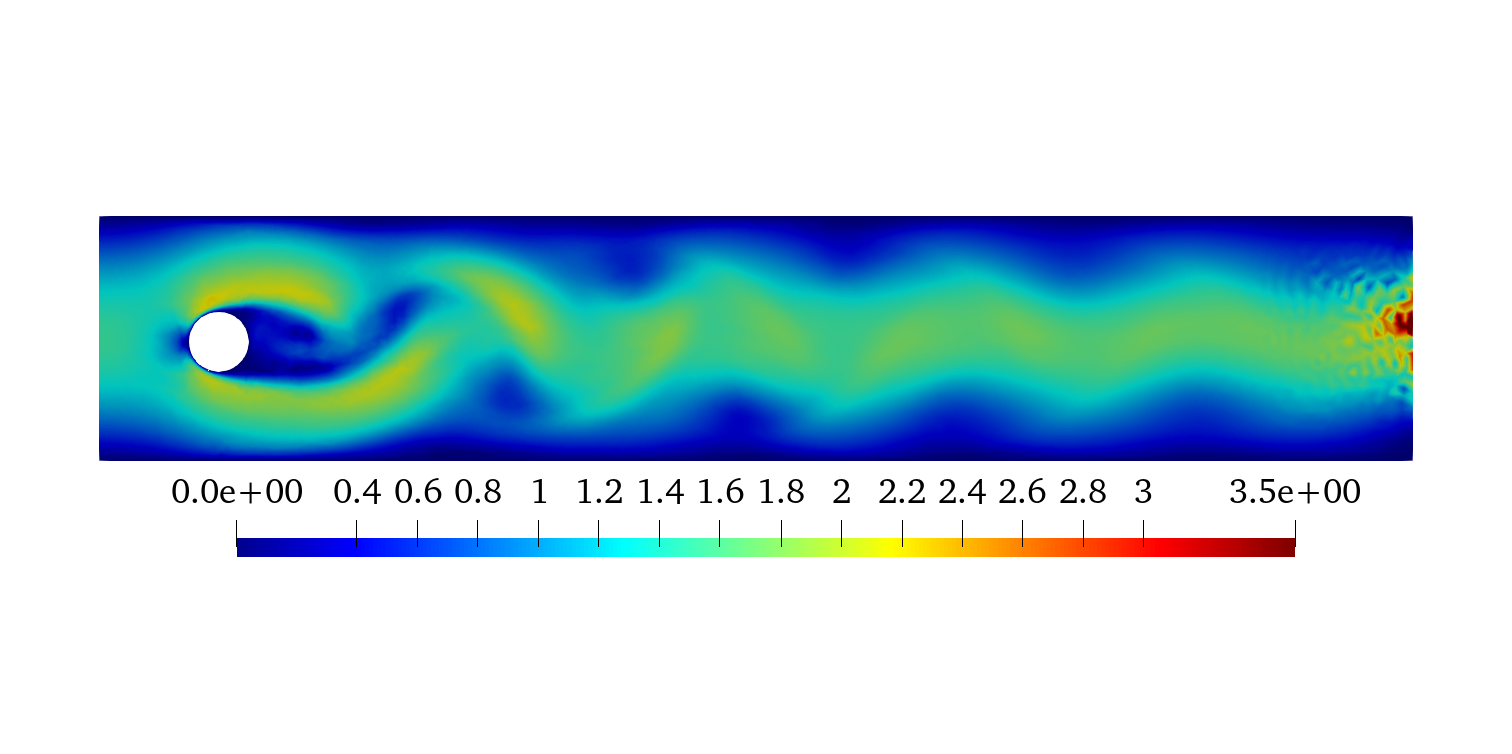}}
    \subfloat[Standard EF ($\delta=\eta$) - velocity ($t=4$)]{\includegraphics[width=0.5\textwidth, trim={3cm 6.5cm 3cm 7cm}, clip]{images/updated_Re1000/u_deltakolm_EF.png}}
    \caption{Velocity fields at final time $t=4$ for standard EF, $\delta$-Opt-EF, and 
    reference DNS. 
    }
    \label{fig:opt-ef-u}
\end{figure*}

\begin{figure*}[htpb!]
    \centering
    
    \subfloat[DNS - pressure ($t=4$)]{\includegraphics[width=0.5\textwidth, trim={3cm 6.5cm 3cm 7cm}, clip]{images/updated_Re1000/p_reference.png}}
   \\
    \subfloat[\optEFglob{} - pressure ($t=4$)]{\includegraphics[width=0.5\textwidth, trim={3cm 6.5cm 3cm 7cm}, clip]{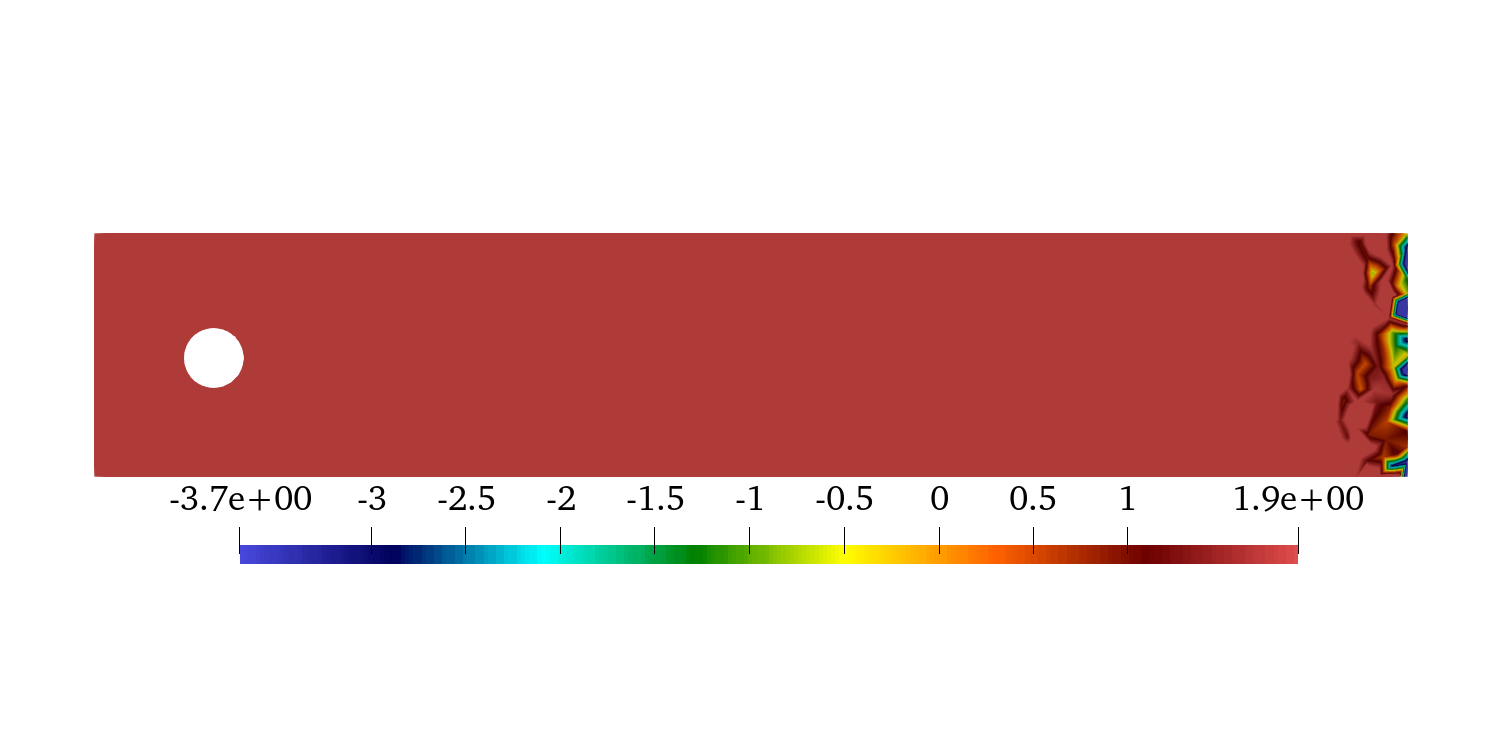}} 
    \subfloat[\optEFglobGRAD{} - pressure ($t=4$)]{\includegraphics[width=0.5\textwidth, trim={3cm 6.5cm 3cm 7cm}, clip]{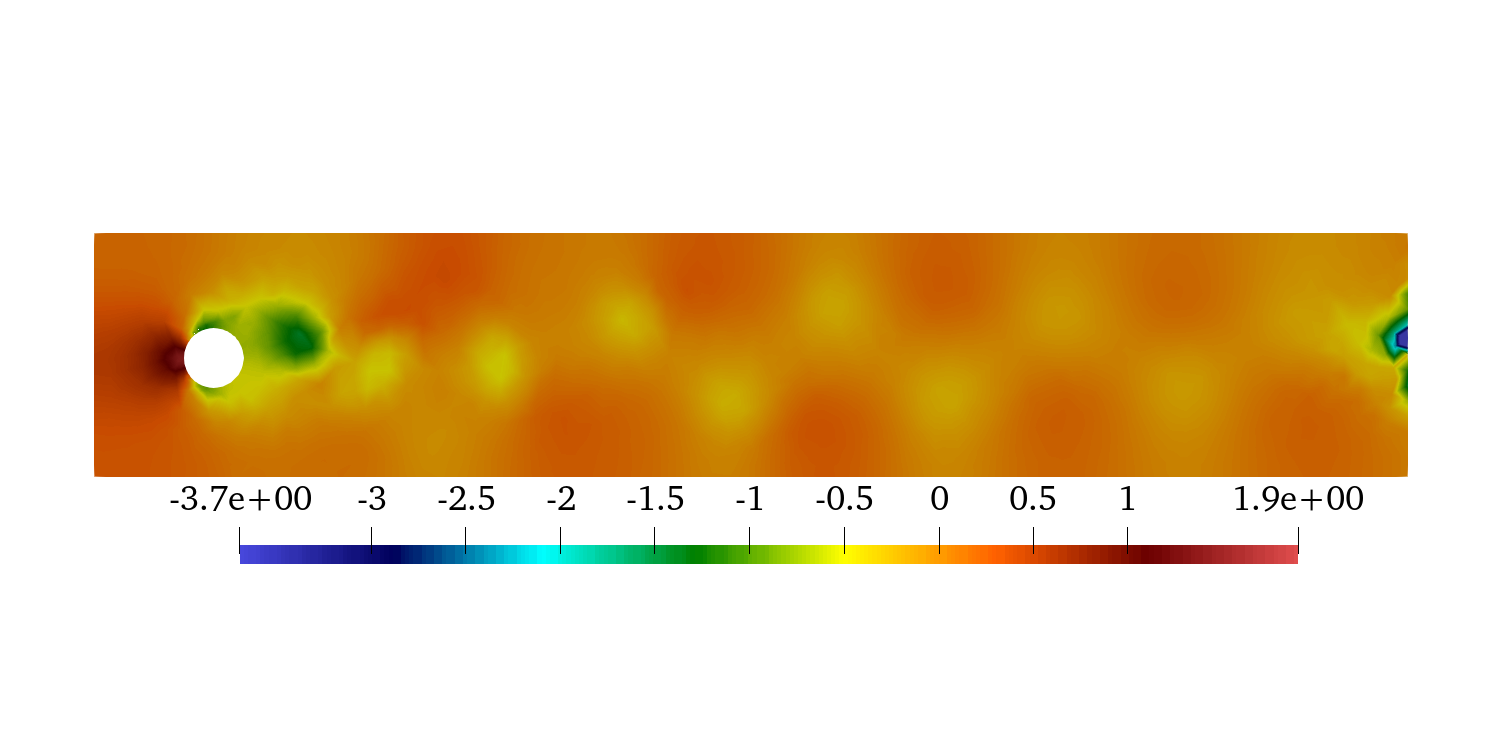}}\\
    \subfloat[\optEFglobPRESS{} - pressure ($t=4$)]{\includegraphics[width=0.5\textwidth, trim={3cm 6.5cm 3cm 7cm}, clip]{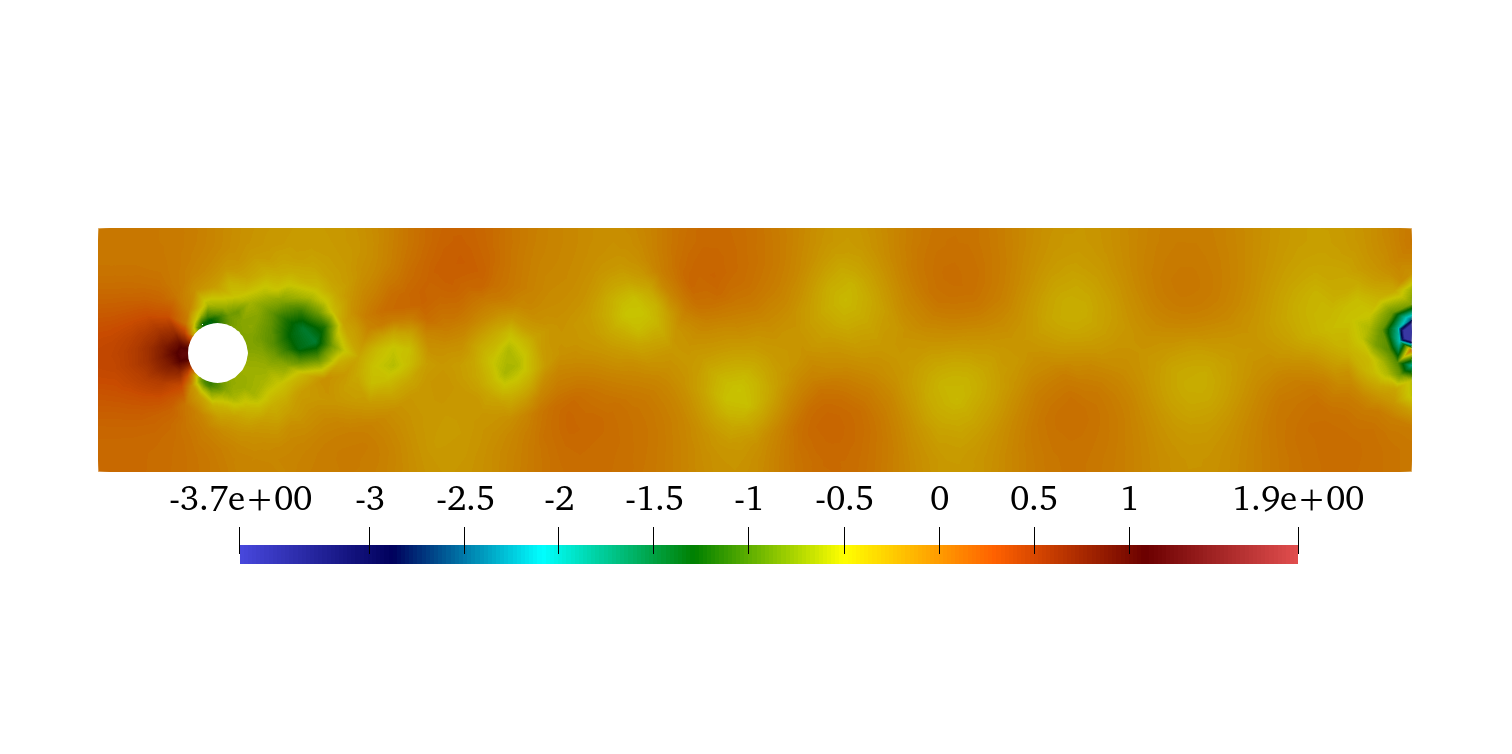}}
     \subfloat[Standard EF ($\delta=\eta$) - pressure ($t=4$)]{\includegraphics[width=0.5\textwidth, trim={3cm 6.5cm 3cm 7cm}, clip]{images/updated_Re1000/p_deltakolm_EF.png}}
    \caption{Pressure fields at final time $t=4$ for standard EF, $\delta$-Opt-EF, and 
    reference DNS. 
    }
    \label{fig:opt-ef-p}
\end{figure*}

The trend of the norms in Figures \ref{fig:norms-delta} and \ref{fig:norm-delta-h1}, and the graphical results in Figures \ref{fig:opt-ef-u} and \ref{fig:opt-ef-p} lead to the following conclusions: 
\begin{itemize}
    \item The optimized algorithms give a better representation of the reference fields, outperforming the over-diffusive approximation of standard EF with fixed $\delta$;
    \item Comparing 
    the $L^2$ norms (Figures \ref{fig:norms-delta}), we can see that the most accurate result 
    is obtained using the \optEFglob{} algorithm, which directly includes the measured metric inside the objective function. However, this algorithm provides inaccurate results, as seen both from the $L^2$ pressure norm plots and the $H^1$ velocity norm plots. This is confirmed by the graphical results (Figures \ref{fig:opt-ef-u} and \ref{fig:opt-ef-p} (B)), which show noisy and poor approximations;
    \item The \optEFglobGRAD{} and the \optEFglobPRESS{} approaches are more accurate than the \optEFglob{} one, as can be seen from the velocity $H^1$ and the pressure $L^2$ norm plots. \RB{These approaches are also closer to the reference DNS data than the standard EF, as highlighted in the zoomed-in areas.} Also the graphical results (Figures \ref{fig:opt-ef-u} and \ref{fig:opt-ef-p} (C) and (D)) confirm a good agreement with the reference DNS data (displayed in Figures \ref{fig:opt-ef-u} and \ref{fig:opt-ef-p} (A));
    \item As noticed for the $\chi$-Opt-EFR algorithms, the inclusion of the pressure contribution \optEFglobPRESS{} does not improve the pressure accuracy and negatively impacts the computational cost of the method.
\end{itemize}

\newpage

%% file: sections/results-D-Opt-EFR.tex
In this final part, we discuss the results of the double optimization algorithms $\delta \chi$-Opt-EFR.
Since the algorithms \optEFRglob{} and \optEFglob{} proved to be inaccurate in the previous subsections \ref{sec:results-chi-opt} and \ref{sec:results-delta-opt}, we here only consider algorithms \DoptEFRglobGRAD{} and \DoptEFRglobPRESS{}.

\begin{figure}[htpb!]
    \centering
    \subfloat[]{\includegraphics[width=\columnwidth]{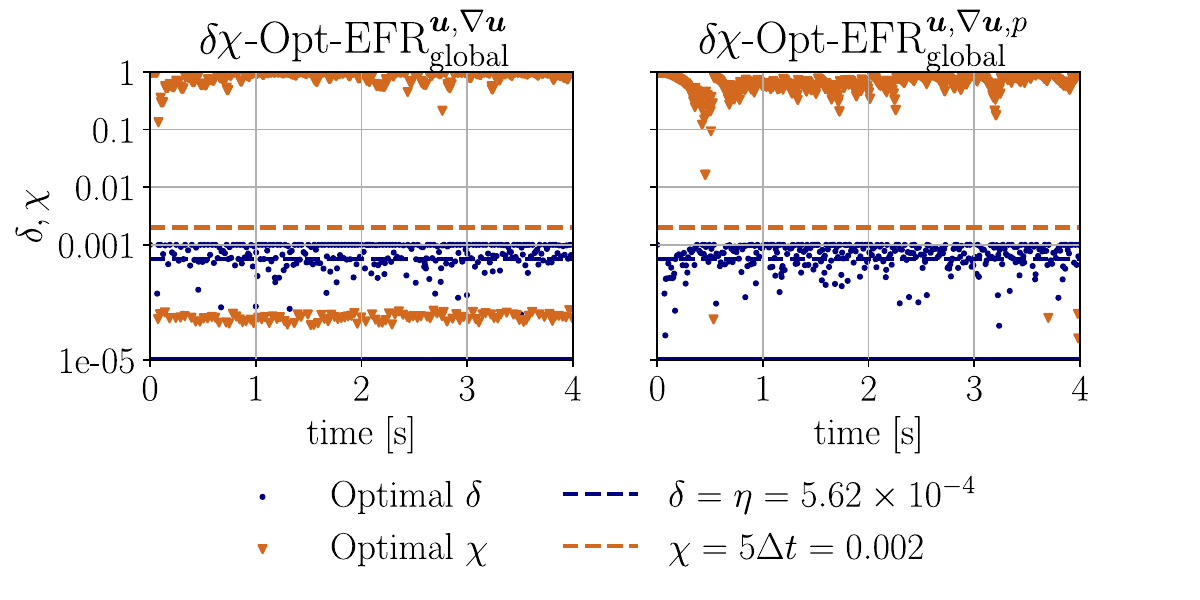}}\\
    \subfloat[]{\includegraphics[width=\columnwidth]{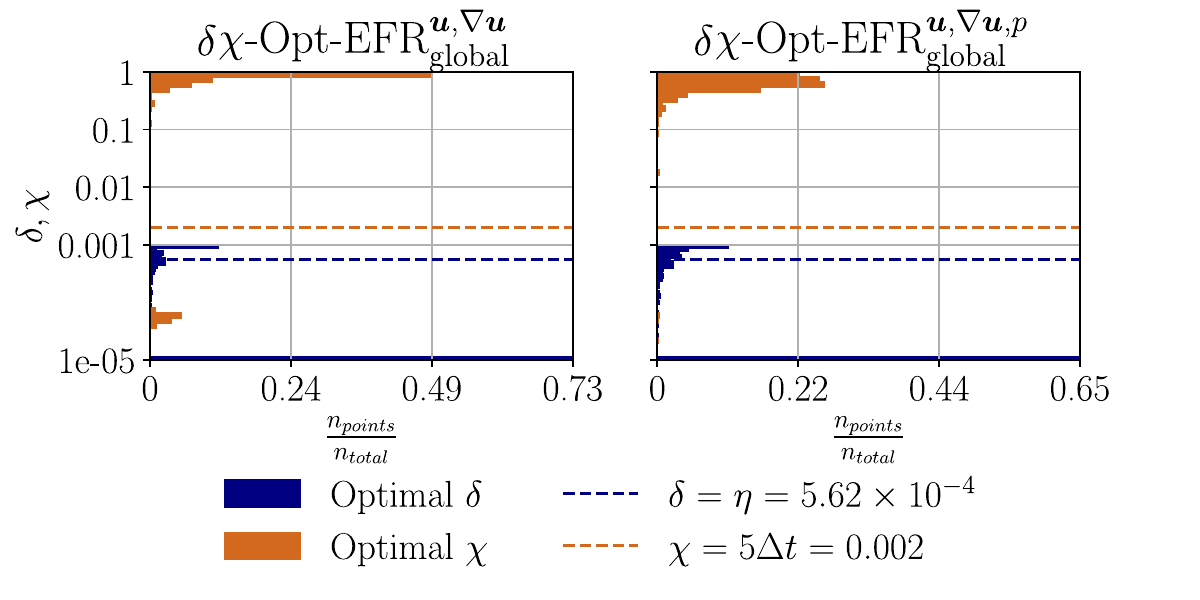}}
    \caption{Optimal value of $\delta(t)$ and $\chi(t)$ in the $\delta \chi$-Opt-EFR algorithms (A), and 
    relative parameter distribution (B).}
    \label{fig:double}
\end{figure}

Figure \ref{fig:double} (A) shows the optimal parameters' trend in time for the two proposed algorithms. 
The results for parameter $\delta(t)$ are similar in the two cases. However, the relaxation parameter $\chi(t)$ has statistically higher values in algorithm \DoptEFRglobPRESS{}, as can be seen from the histogram in Figure \ref{fig:double} (B). Moreover, in the double optimizations the optimal $\chi(t)$ is 
almost always in the range $[\num{1e-1}, 1]$, while in optimizations $\chi$-Opt-EFR$_{\text{global}}$ (Figure \ref{fig:chi}) 
it has a more oscillating behavior and it also converges to smaller values.

\begin{figure*}[htpb!]
    \centering
    \includegraphics[width=0.9\textwidth]{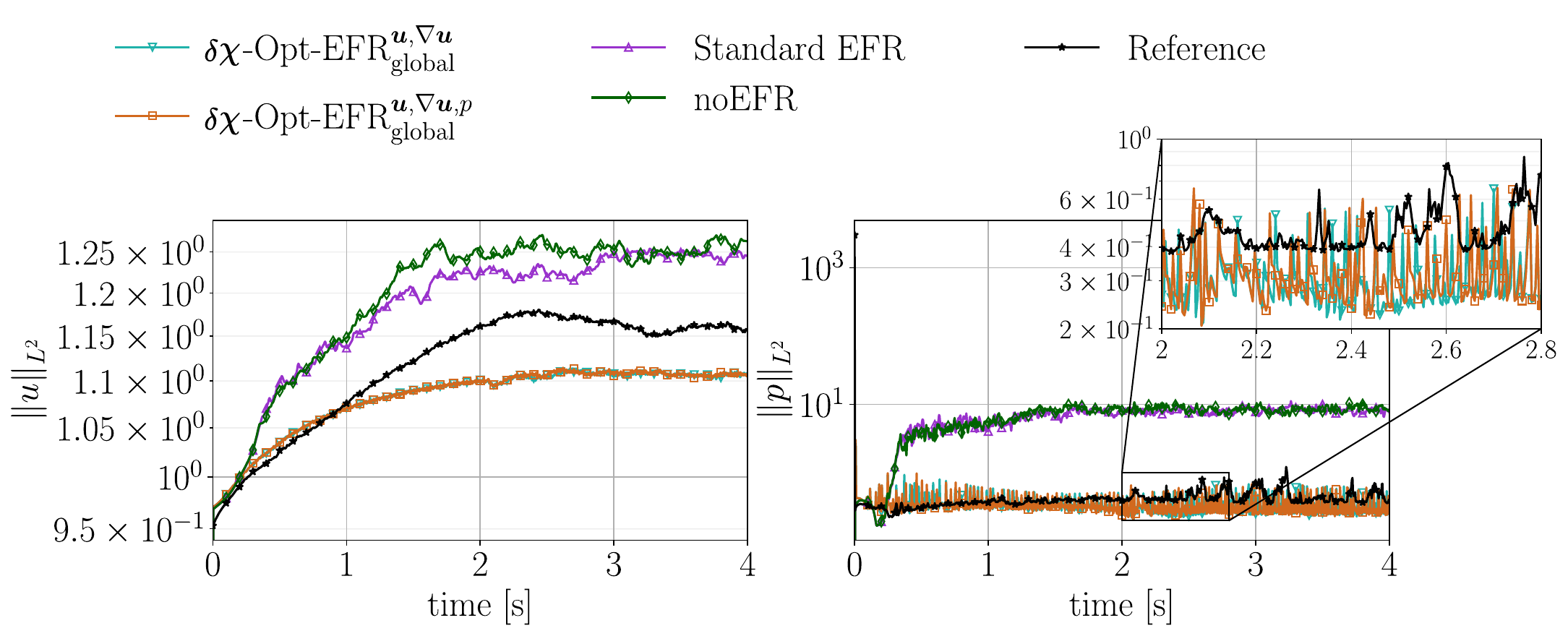}
    \caption{Velocity and pressure norms of $\delta \chi$-Opt-EFR simulations 
    and 
    reference solutions 
    projected on the coarse mesh. \RB{The pressure plot (right) also includes a box with a zoomed-in area in the time interval $[2, 2.8]$.}}
    \label{fig:norms-double}
\end{figure*}

\begin{figure}[htpb!]
    \centering
    \includegraphics[width=0.5\textwidth]{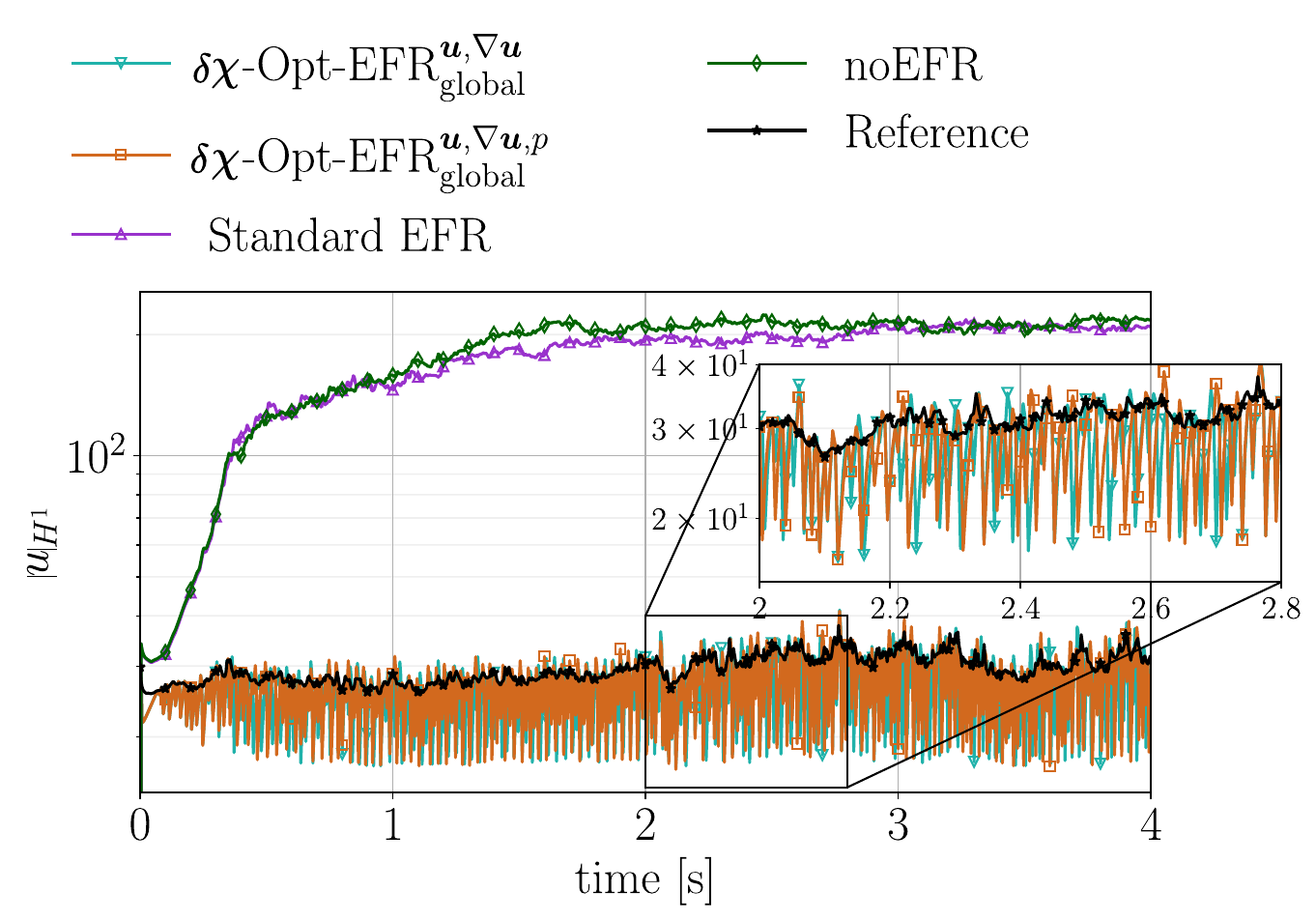}
    \caption{Velocity $H^1$ norms of $\delta \chi$-Opt-EFR simulations 
    and 
    reference solutions 
    projected on the coarse mesh. \RB{The plot also includes a box with a zoomed-in area in the time interval $[2, 2.8]$.}}
    \label{fig:norm-double-h1}
\end{figure}

\begin{figure*}[htpb!]
    \centering
    \subfloat[DNS - velocity ($t=4$)]{\includegraphics[width=0.5\textwidth, trim={3cm 7cm 3cm 7cm}, clip]{images/updated_Re1000/u_reference.png}} 
    \subfloat[Standard EFR - velocity ($t=4$)]{\includegraphics[width=0.5\textwidth, trim={3cm 7cm 3cm 7cm}, clip]{images/updated_Re1000/u_EFR_fixed.png}} 
    \\
    \subfloat[\DoptEFRglobGRAD{} - velocity ($t=4$)]{\includegraphics[width=0.5\textwidth, trim={3cm 7cm 3cm 7cm}, clip]{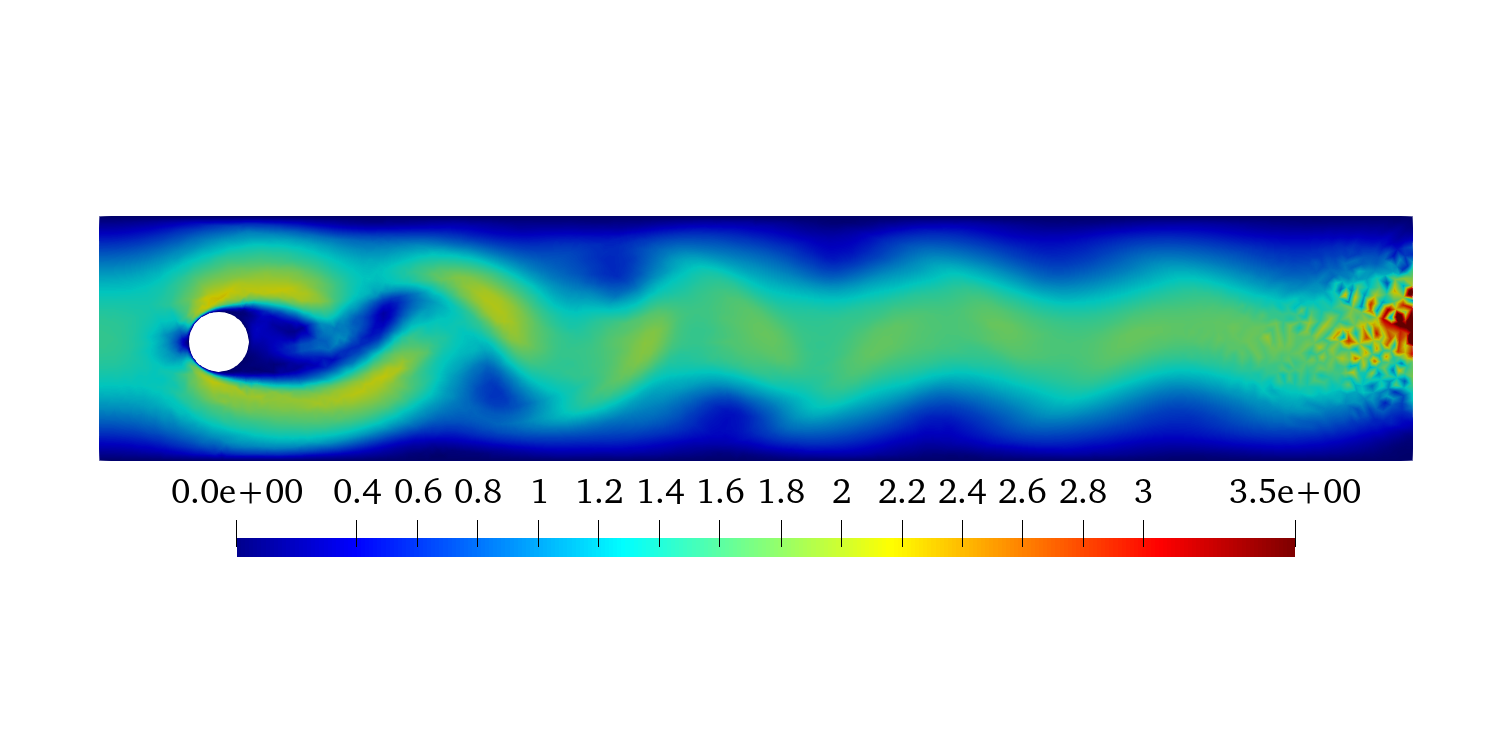}} 
    \subfloat[\DoptEFRglobPRESS{} - velocity ($t=4$)]{\includegraphics[width=0.5\textwidth, trim={3cm 7cm 3cm 7cm}, clip]{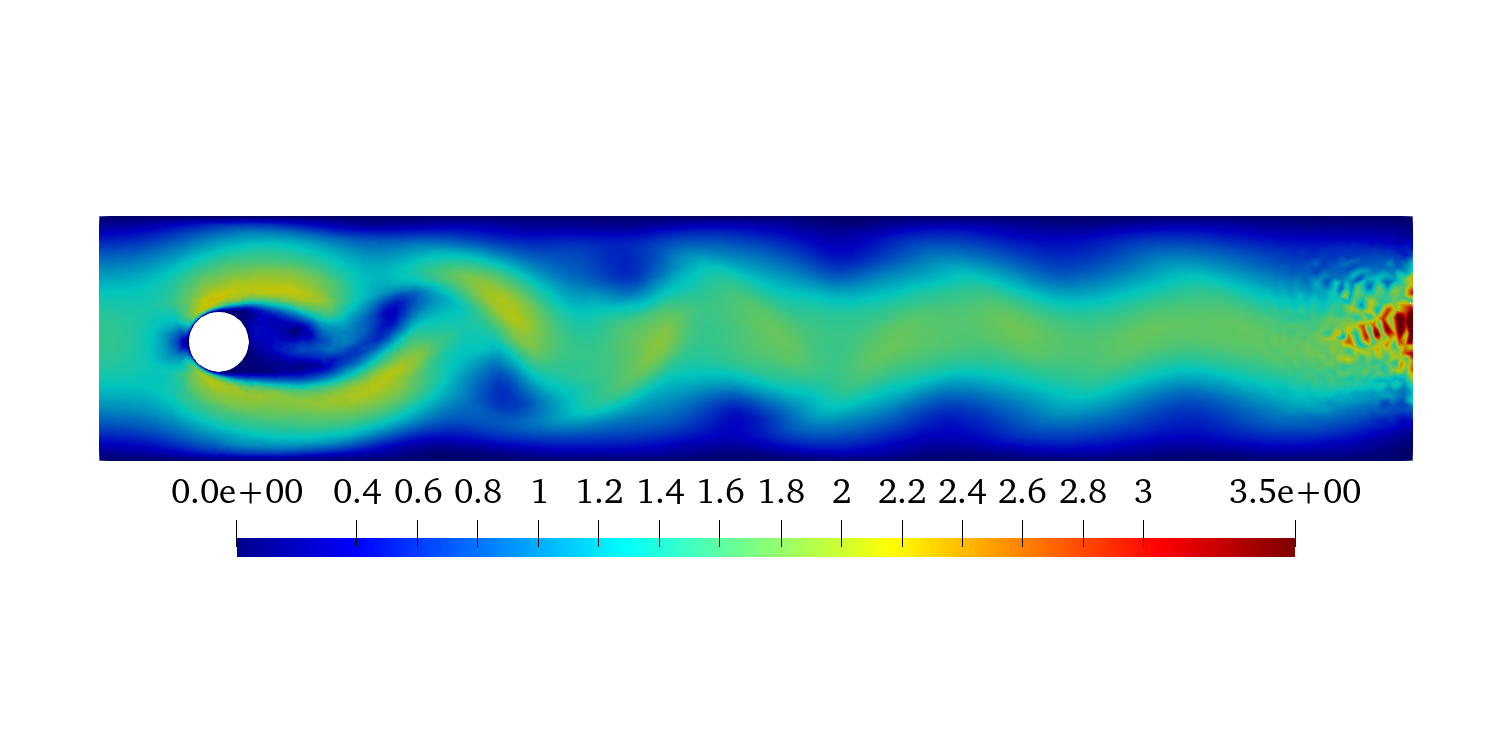}} 
    \caption{Velocity fields at final $t=4$ for $\delta \chi$-Opt-EFR simulations, 
    DNS reference, and 
    standard EFR.}
    \label{fig:d-opt-efr-u}
\end{figure*}

\begin{figure*}[htpb!]
    \centering
    \subfloat[DNS - pressure ($t=4$)]{\includegraphics[width=0.5\textwidth, trim={3cm 7cm 3cm 7cm}, clip]{images/updated_Re1000/p_reference.png}} 
    \subfloat[Standard EFR - pressure ($t=4$)]{\includegraphics[width=0.5\textwidth, trim={3cm 7cm 3cm 7cm}, clip]{images/updated_Re1000/p_EFR_fixed.png}} \\
    \subfloat[\DoptEFRglobGRAD{} - pressure ($t=4$)]{\includegraphics[width=0.5\textwidth, trim={3cm 6.5cm 3cm 7cm}, clip]{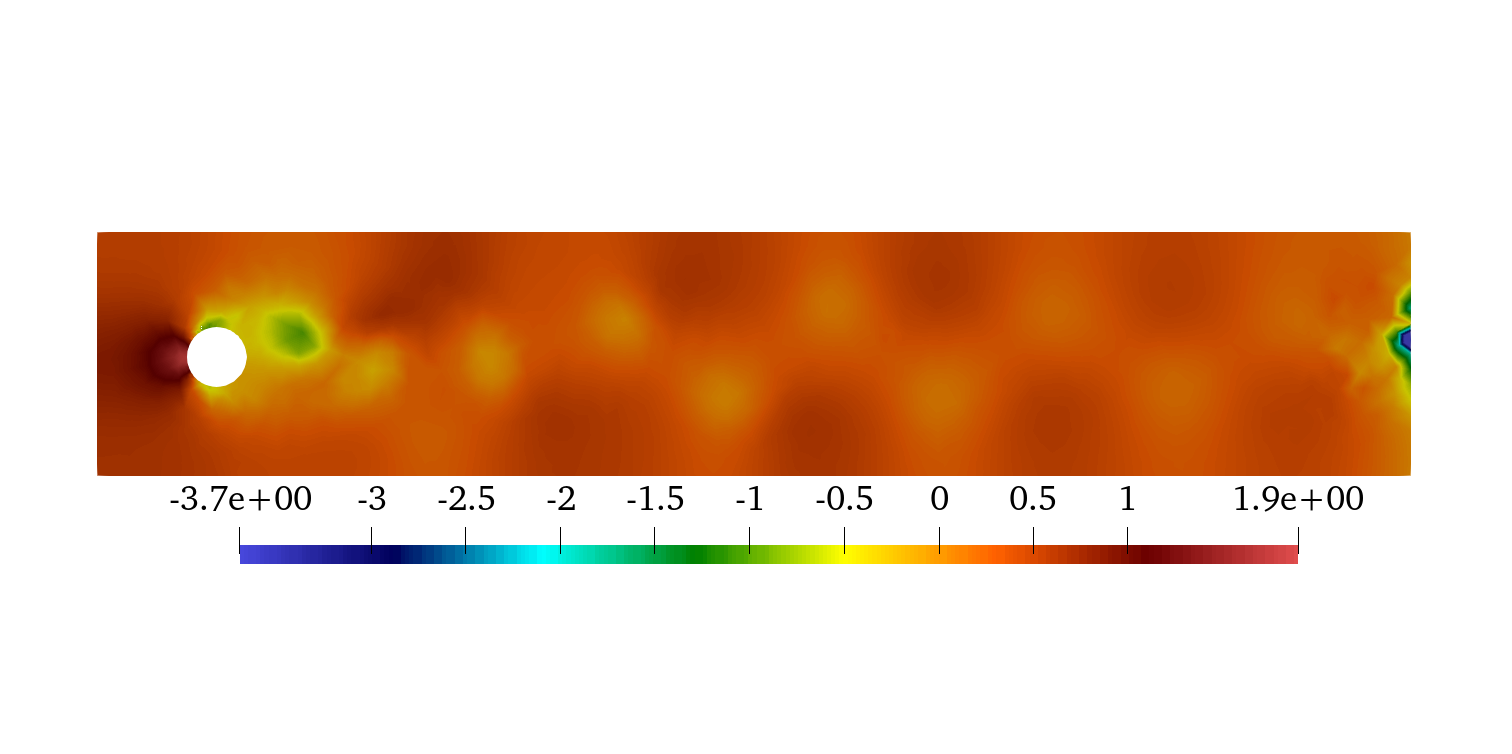}} 
    \subfloat[\DoptEFRglobPRESS{} - pressure ($t=4$)]{\includegraphics[width=0.5\textwidth, trim={3cm 6.5cm 3cm 7cm}, clip]{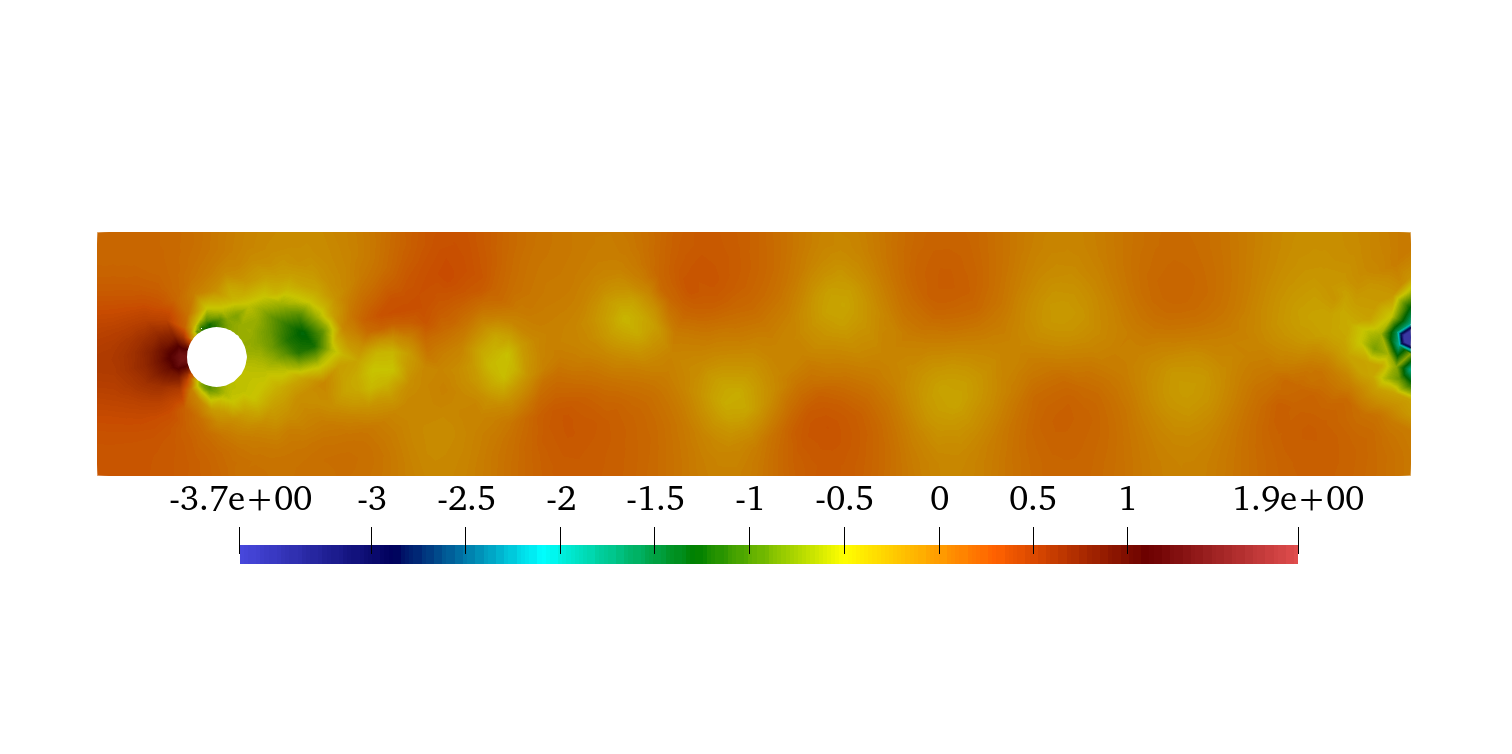}} 
    \caption{Pressure fields at final time $t=4$ for $\delta \chi$-Opt-EFR, 
    DNS reference, and 
    standard EFR.}
    \label{fig:d-opt-efr-p}
\end{figure*}

The results for the $L^2$ and $H^1$ norms (Figures \ref{fig:norms-double} and \ref{fig:norm-double-h1})
are similar \RB{for the two algorithms proposed (\DoptEFRglobGRAD{} and \DoptEFRglobPRESS{}), as also highlighted in the zoomed-in areas in the plots}.
If we also look at the graphical representations of the velocity and pressure fields at $t=4$ in Figures \ref{fig:d-opt-efr-u} and \ref{fig:d-opt-efr-p} (C) and (D), the behavior is quite similar. However, the final pressure field is slightly more accurate when the functional takes also into account the pressure contribution, namely in \DoptEFRglobPRESS{}.

\subsection{On the role of the optimization step}
\label{subsec:opt-step-final}

\RA{The previous analysis was performed with a fixed optimization step, $\Delta t_{opt}=10$.
However, the optimization step plays an important role in ensuring a good compromise between the accuracy and the efficiency of the described approach.
Hence, this section is dedicated to the comparison of the different types of optimization, in terms of \textbf{efficiency} and \textbf{accuracy}. For each 
algorithm, we also extend the analysis to different optimization steps, namely $\Delta t_{opt}=k\Delta t$, where the positive integer $k$ satisfies $k \in \{10,20, \cdots, 90, 100\}$.}

\begin{figure*}[htpb!]
    \centering
    \subfloat[Discrepancy on the velocity $L^2$ norm]{
    \includegraphics[width=1\textwidth]{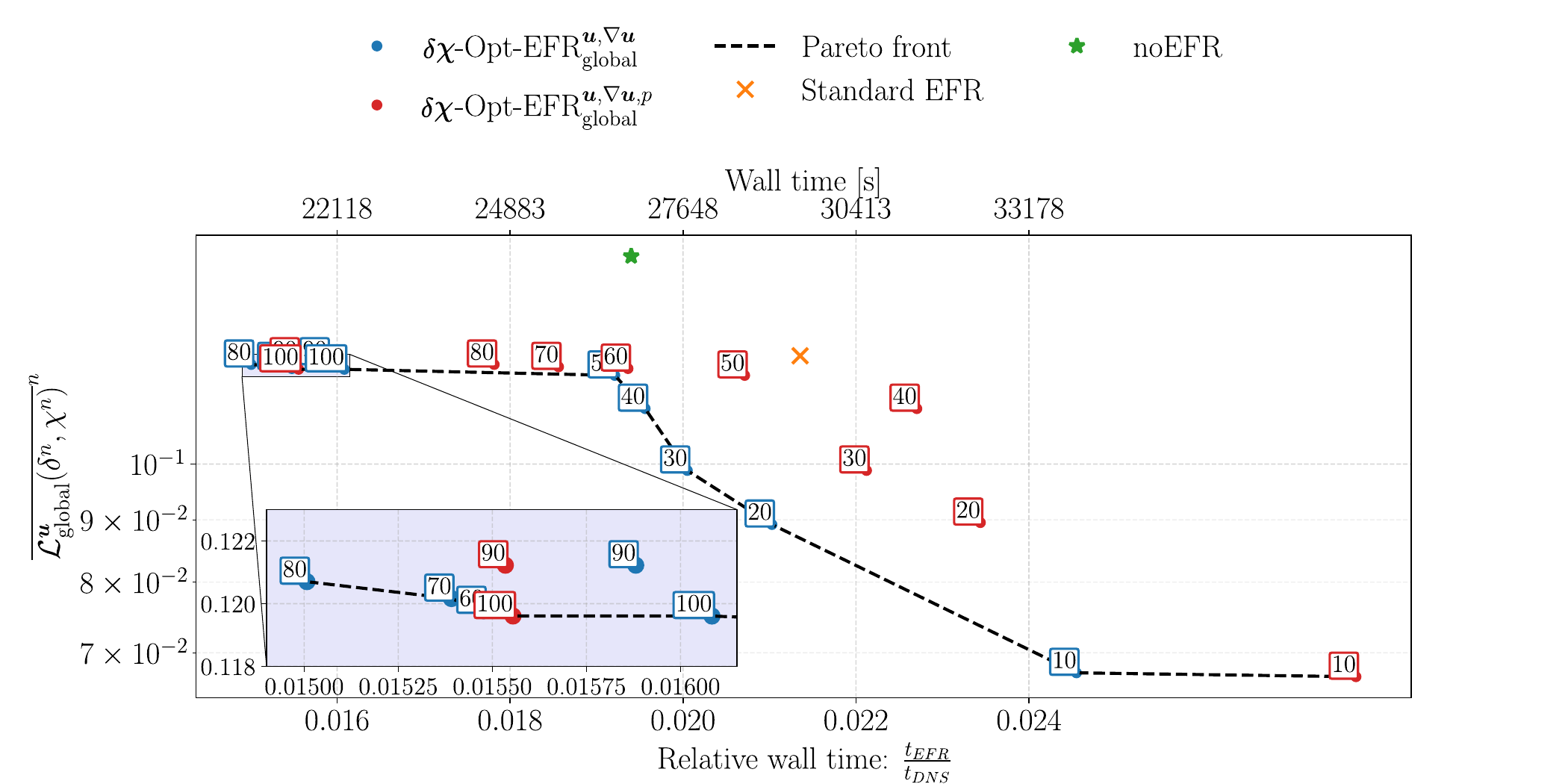}}\\
    \subfloat[Discrepancy on the pressure $L^2$ norm]{
    \includegraphics[width=1\textwidth]{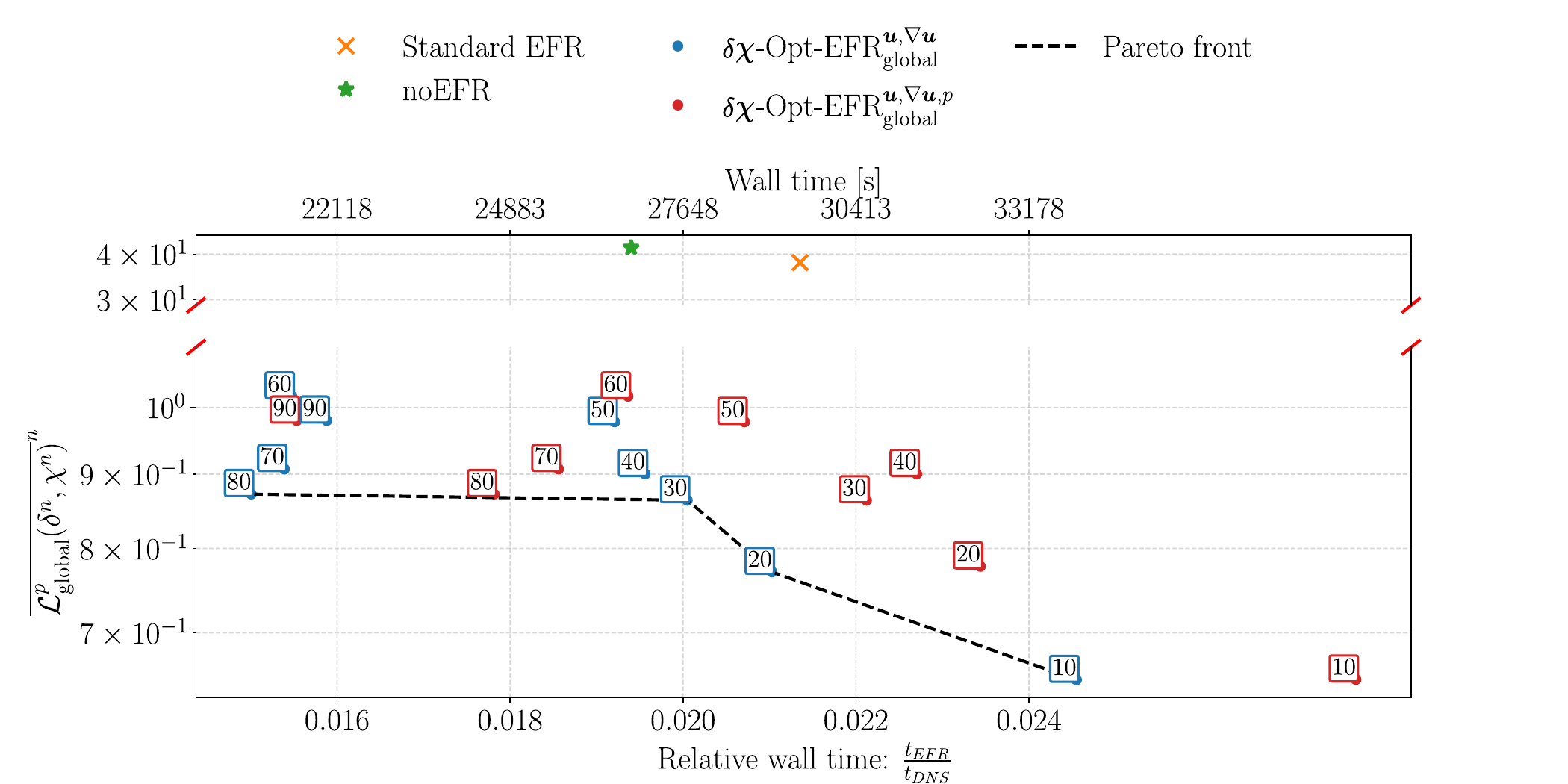}}
    \caption{Representation of how the optimization step affects the efficiency and the accuracy of the final velocity or pressure approximation in $\DoptEFRglobGRAD{}$ and $\DoptEFRglobPRESS{}$ simulations. The efficiency, on the $x$-axis, is measured through the relative wall-clock time with respect to the DNS time. The accuracy, on the $y$-axis, is measured with the $L^2$ velocity and pressure norms.}
    \label{fig:d-opt-efr-opt-step}
\end{figure*}

\begin{figure*}[htpb!]
    \centering
        \subfloat[Discrepancy on the velocity $H^1$ semi-norm]{
    \includegraphics[width=1\textwidth]{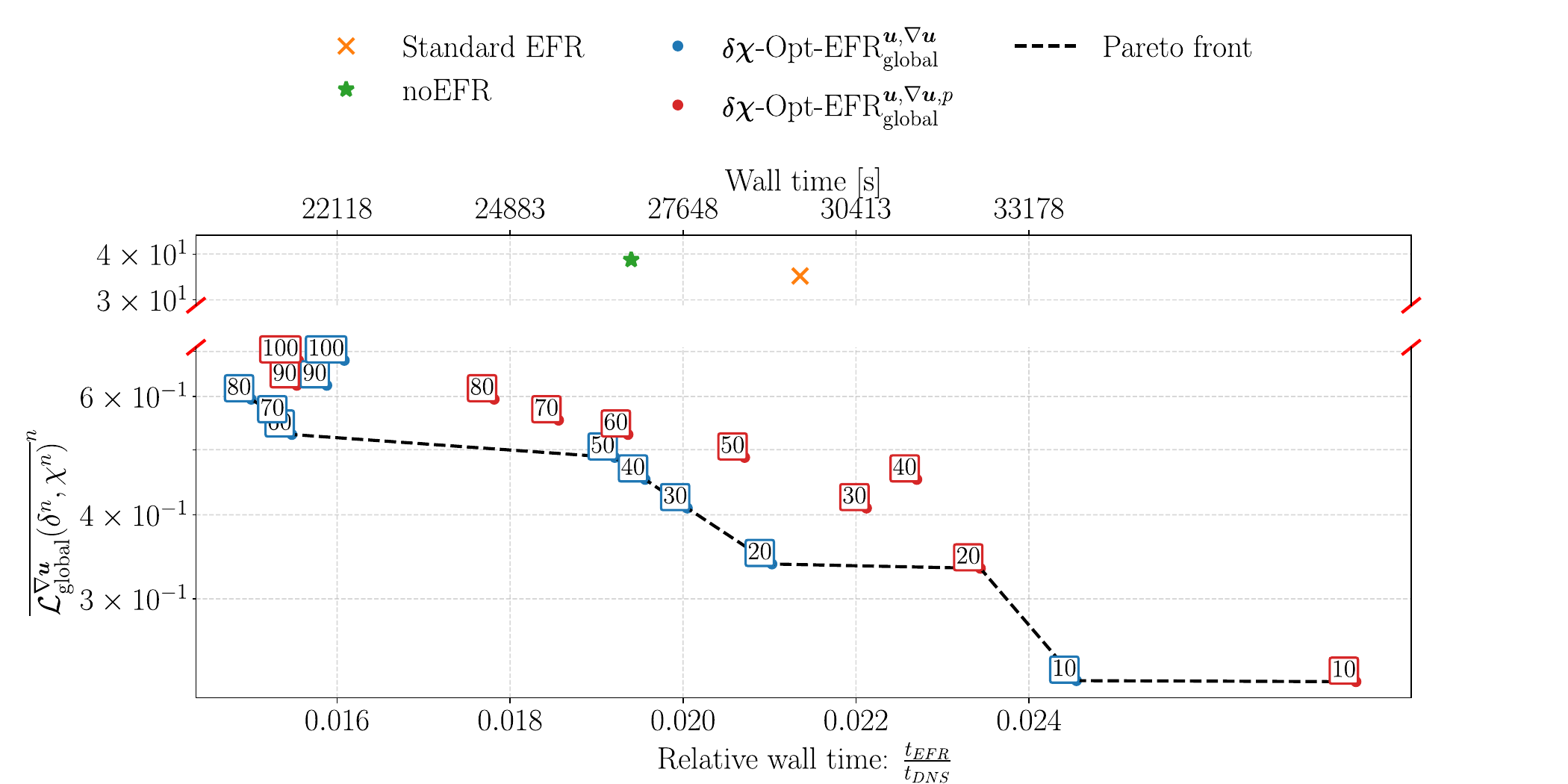}}
    \caption{Representation of how the optimization step affects the efficiency and the accuracy of the final velocity or pressure approximation in $\optEFRglobGRAD{}$ and $\optEFRglobPRESS{}$ simulations. The efficiency, on the $x$-axis, is measured through the relative wall-clock time with respect to the DNS time. The accuracy, on the $y$-axis, is measured with the $H^1$ velocity seminorm.}
    \label{fig:d-opt-efr-opt-step-h1}
\end{figure*}

\RA{We measure the \textbf{accuracy} using metrics $\overline{\mathcal{L}^{\bu}_{\mathrm{global}}(\boldsymbol{\mu}^n)}^n$, $\overline{\mathcal{L}^{\nabla \bu}_{\mathrm{global}}(\boldsymbol{\mu}^n)}^n$, and $\overline{\mathcal{L}^{p}_{\mathrm{global}}(\boldsymbol{\mu}^n)}^n$, where the operator $\overline{(\cdot)}^n$ is a time average, and the terms $\mathcal{L}^{(\cdot)}_{\mathrm{global}}(\boldsymbol{\mu}^n)$ are the global contributions in the optimization introduced in Equations \eqref{eq:global-contributions-general}.

We measure the \textbf{efficiency} 
in terms of the wall-clock time, expressed in relative percentage with respect to the DNS time (bottom $x$-axis) and in absolute time (top $x$-axis).

The results of the $\delta \chi$-Opt-EFR simulations with different $k$ values are presented in the Pareto plots in Figures \ref{fig:d-opt-efr-opt-step} and \ref{fig:d-opt-efr-opt-step-h1}, where the $k$ value is visualized in the boxes next to the corresponding markers.

Since in the plots the $x$-axis represents the computational time and the $y$ axis represents the errors, we expect the best-performing simulations to be located in the lower left part of the plots, which corresponds to low errors and 
low computational time.
Moreover, our goal is to discard simulations taking place in the upper left part, corresponding to high errors and low computational cost, and in the lower right part, corresponding to low errors and high computational cost.

Figures \ref{fig:d-opt-efr-opt-step} and \ref{fig:d-opt-efr-opt-step-h1} yield the following conclusions:

\begin{itemize}
    \item The standard EFR and the noEFR simulations lead to higher errors than the optimized-EFR approaches, as we expect. In particular, for the pressure and the velocity gradients, $\delta \chi$-Opt-EFR methods are able to improve the accuracy by \textbf{two} orders of magnitudes. Moreover, the standard EFR and the noEFR simulations also need more computational time than some of the $\delta \chi$-Opt-EFR simulations. This may be due to the fact that the noEFR and standard-EFR simulations are \textbf{unstable and noisy}, and, hence, the FEM Newton solver needs more iterations to converge to a solution.
    \item With respect to the relative wall-clock time, all the optimized-EFR simulations have \textbf{high gain in time} with respect to the DNS simulation. However, the absolute time suggests that the range of all the $\delta \chi$-Opt-EFR simulations is $[4, 9]$ hours. 
    This means that varying the optimization time window, the change in the computational time of the $\delta \chi$-Opt-EFR simulations may be on the order of hours.
    \item As expected, the simulations with a smaller optimization step, i.e., $10 \Delta t$ or $20 \Delta t$, need more time than the simulations with higher optimization steps, but are more accurate. 
    We also not\RB{e} that the accuracy among different optimization steps is still comparable in terms of the $L^2$ norms, but not in terms of the $H^1$ velocity seminorm. This may be due to the fact that the velocity gradients get too far from the global optimum before a novel optimization is performed. Hence, the following optimization step is not able anymore to converge and to ``fix" the spatial gradients.
    \item The \DoptEFRglobPRESS{} strategy slightly improves the pressure accuracy. However, it is more time consuming than the \DoptEFRglobGRAD{} without yielding a significant gain in the accuracy.
    
\end{itemize}
}

\RA{The same sensitivity study for $\chi$-Opt-EFR and $\delta$-Opt-EF can be found in the supplementary results (sections \ref{subsubsec:supp-pareto-chi} and \ref{subsubsec:supp-pareto-delta}). The above considerations for $\delta \chi$-Opt-EFR hold true also for the other two strategies.}
\RA{
\remark{The Pareto plots in Figures \ref{fig:d-opt-efr-opt-step} and \ref{fig:d-opt-efr-opt-step-h1} show how the optimization step influences the simulation wall-clock time. One would expect that simulations with larger time windows are faster than simulations with smaller time windows. However, it is important to underline that this general statement is not always satisfied. For example, simulations with $\Delta t_{opt} \in \{70, 80, 90, 100\}$ take more time than simulations with $\Delta t_{opt} \in \{40,50, 60\}$. This is because larger optimization steps lead to 
more unstable systems (i.e., spurious oscillations and presence of noise), and in this scenario the Newton method needs more iterations to converge to the solution.}
}
\subsection{Final Discussion}
\label{discussion}

We here consider a comparison among the algorithms having the best performance in the previous sections, namely $(\cdot)_{\text{global}}^{\bu, \nabla \bu, p}$.

\begin{figure*}[htpb!]
    \centering
    \includegraphics[width=0.9\textwidth]{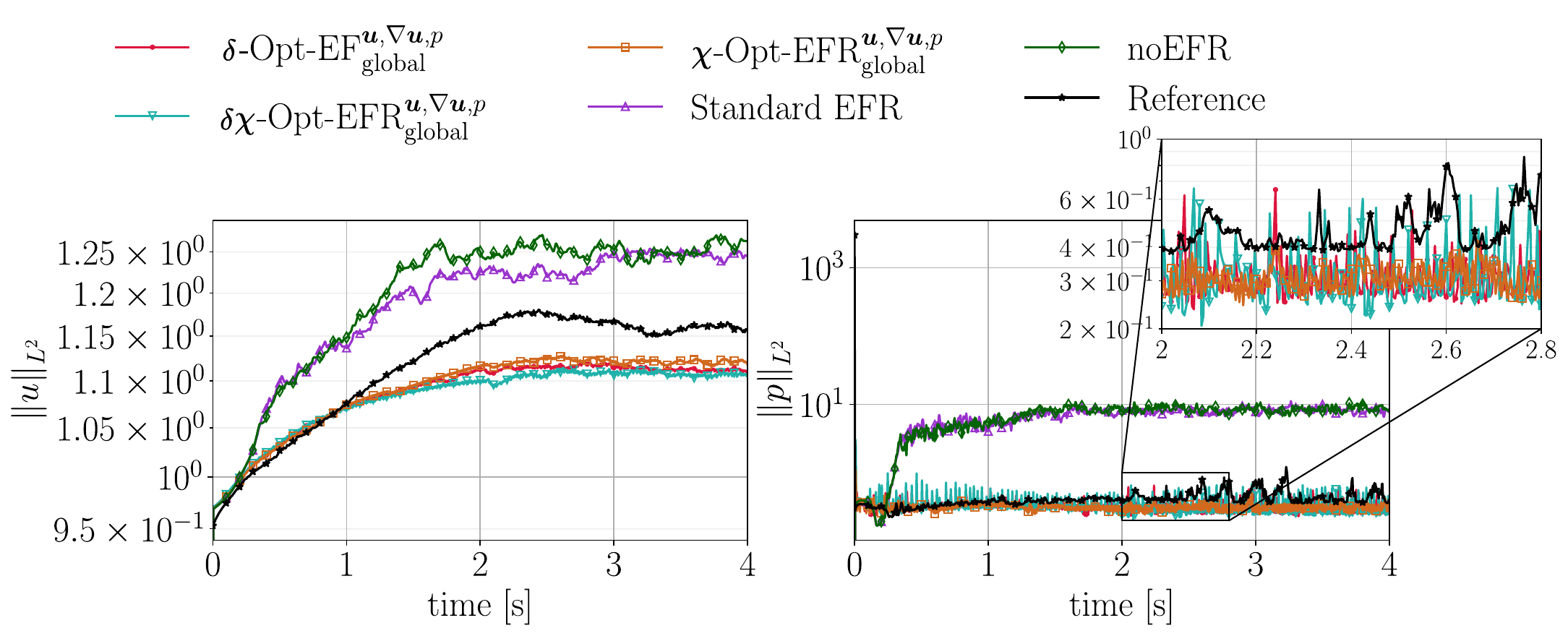}
    \caption{Velocity and pressure norms for algorithms $(\cdot)_{\text{global}}^{\bu, \nabla \bu, p}$ 
    and reference solutions projected on the coarse mesh. \RB{The pressure plot (right) also includes a box with a zoomed-in area in the time interval $[2, 2.8]$.}}
    \label{fig:norms-comparison}
\end{figure*}

\begin{figure}[htpb!]
    \centering
    \includegraphics[width=0.5\textwidth]{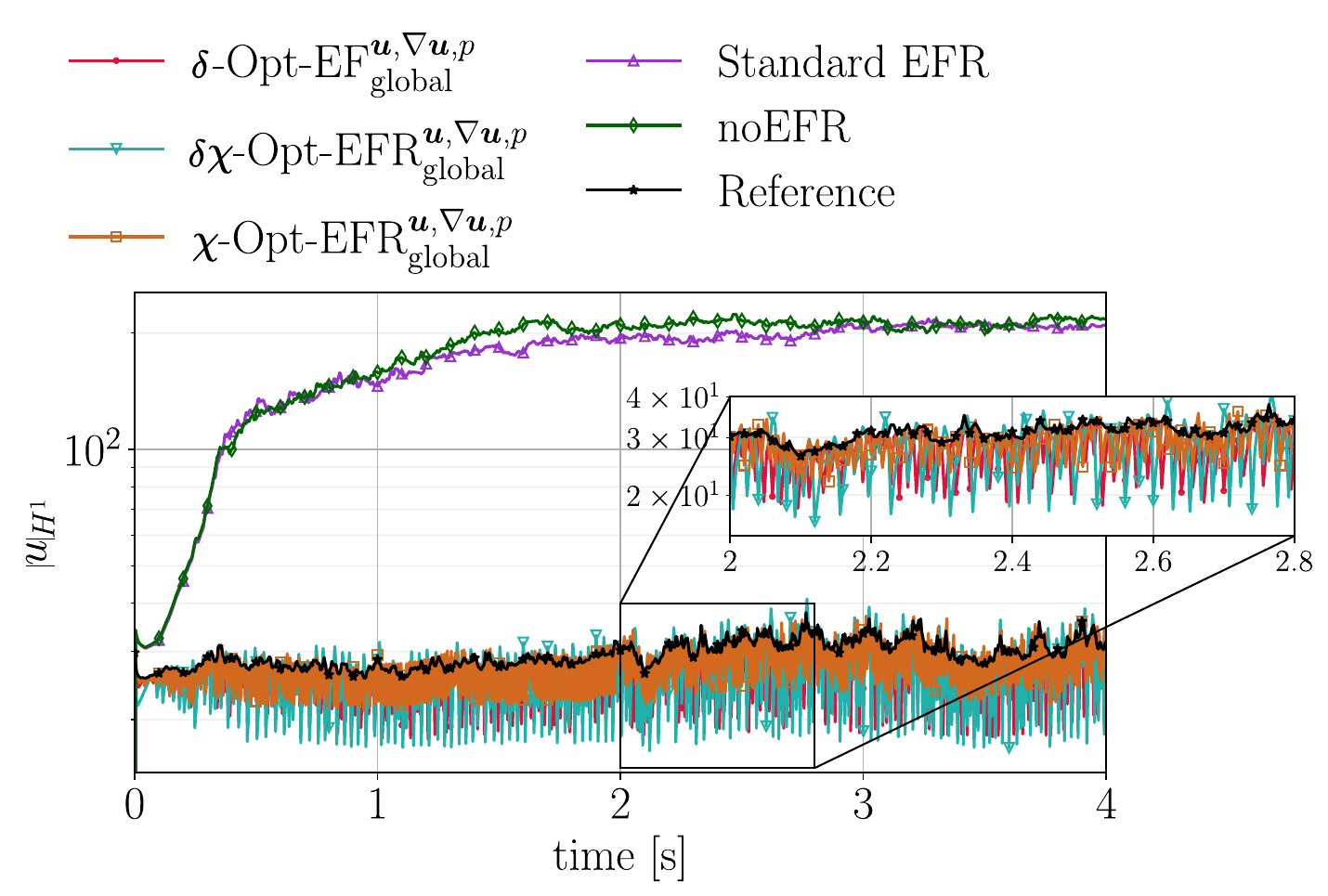}
    \caption{Velocity $H^1$ norms for algorithms $(\cdot)_{\text{global}}^{\bu, \nabla \bu, p}$ and reference solutions projected on the coarse mesh. \RB{The plot also includes a box with a zoomed-in area in the time interval $[2, 2.8]$.}}
    \label{fig:norm-comparison-h1}
\end{figure}

The comparison in the norms in Figures \ref{fig:norms-comparison} and \ref{fig:norm-comparison-h1} show that all the results are similar. However, in the $H^1$ velocity norm, algorithm \optEFRglobPRESS{} is less oscillating than the others \RB{(as one can notice from the zoomed-in areas)}, and is also slightly more accurate in the $L^2$ velocity norm.
\medskip

In general, there is a consistent enhancement in the accuracy between the standard and the Opt-EFR strategies.
This improvement is accompanied by slightly longer times in the case of Algorithm \DoptEFRglobGRAD{}, but considerably longer times for strategy $\DoptEFRglobPRESS{}$, since it includes an evolve step within each optimization pass. $\DoptEFRglobPRESS{}$ slightly improves the accuracy obtained with $\DoptEFRglobGRAD{}$ in terms of velocity gradients and pressure discrepancy.
\medskip

In addition, in what follows we report a series of considerations concerning the comparison among all the algorithms proposed in Sections \ref{sec:results-chi-opt}, \ref{sec:results-delta-opt}, and \ref{sec:results-double-opt}.

\begin{itemize}
    \item The $\chi$-Opt-EFR strategies optimizing \emph{local} quantities, namely algorithms \optEFRloc{} and \optEFRlocGRAD{}, provide over-diffusive results and are not able to fully capture the dynamics of the test case. For this reason, in Sections \ref{sec:results-delta-opt} and \ref{sec:results-double-opt}, we only consider the strategies with \emph{global} loss contributions. 
    
    \item The strategies of type $(\cdot)^{\bu}_{\text{global}}$, optimizing the discrepancy on the \emph{global} quantity $\|\bu\|_{L^2(\Omega)}$, lead to poor results. In particular, the optimized velocity field is characterized by spurious oscillations. That motivates the addition of the gradient contribution to the loss functional.
    
    \item The strategies of type $(\cdot)^{\bu, \nabla \bu}_{\text{global}}$, which include 
    the discrepancy in both $\|\bu\|_{L^2(\Omega)}$ and 
    $\| \nabla \bu\|_{L^2(\Omega)}$, provide the best results in terms of velocity accuracy. Therefore, the inclusion of the gradient contribution within the optimization process is indispensable for achieving favorable outcomes in velocity accuracy.
    
    \item The strategies considering the pressure contribution in the loss functional, namely $(\cdot)^{\bu, \nabla \bu, p}_{\text{global}}$, 
    yield similar results to the algorithms $(\cdot)^{\bu, \nabla \bu}_{\text{global}}$ in terms of velocity accuracy. However, the inclusion of the pressure discrepancy in the loss functional leads to better results in terms of pressure accuracy. 

    \item As discussed previously, the computational times of Opt-EFR algorithms are comparable to the standard EFR in strategy $(\cdot)^{\bu, \nabla \bu}_{\text{global}}$, but with considerably better results. However, the addition of the pressure contribution in the loss functional leads to longer computational times, without a consistent improvement in the accuracy. Moreover, in all the Opt-EFR algorithms, there is a consistent gain in time with respect to the DNS simulation.

    \item We can notice that all the algorithms of type $(\cdot)^{\bu, \nabla u}_{\text{global}}$ and $(\cdot)^{\bu, \nabla u, p}_{\text{global}}$ outperforms the results of standard EFR ($\delta=\eta$, $\chi=5\Delta t$). Moreover, \optEFglobGRAD{} and \optEFglobPRESS{} with $\Delta t_{opt}=10\Delta t$ or $\Delta t_{opt}=20\Delta t$ outperform the results of standard EF ($\delta=\eta$).
    Indeed, standard and fixed values for the filter and the relaxation parameter lead either to spurious oscillations or to over-diffusive fields. Moreover, the computational time of standard EFR is higher than that of the optimized methods 
    because the simulation is not stable and the Newton method needs more iterations to converge to a solution. This suggests that the optimized strategies work also as numerical stabilizations. 

    \item As a general observation, the velocity $H^1$ seminorm has 
    an oscillatory behavior in all the examined cases. This is linked to a similar chaotic behavior of the optimal parameter 
    values. Indeed, the optimization algorithm minimizes the discrepancy in the velocity gradients' norm, so:
    \begin{itemize}
        \item If the norm is far from the reference one, i.e., the velocity starts to exhibit a noisy and nonphysical behavior, which means that the parameters are too small and the algorithm converges to large optimal values;
        \item Consequently, after $\Delta t_{opt}$, the norm is closer to the reference DNS, and the simulation is more stable. Hence, the simulation does not need anymore large filter radii and/or relaxation parameters, and the algorithm converges to smaller optimal values.
    \end{itemize}
    This mechanism repeats as the simulation evolves, leading to oscillatory behaviors in both parameters and norms.
    \end{itemize}

\begin{figure}[htpb!]
    \subfloat[$H^1$ velocity relative errors statistics]{\includegraphics[width=0.8\textwidth]{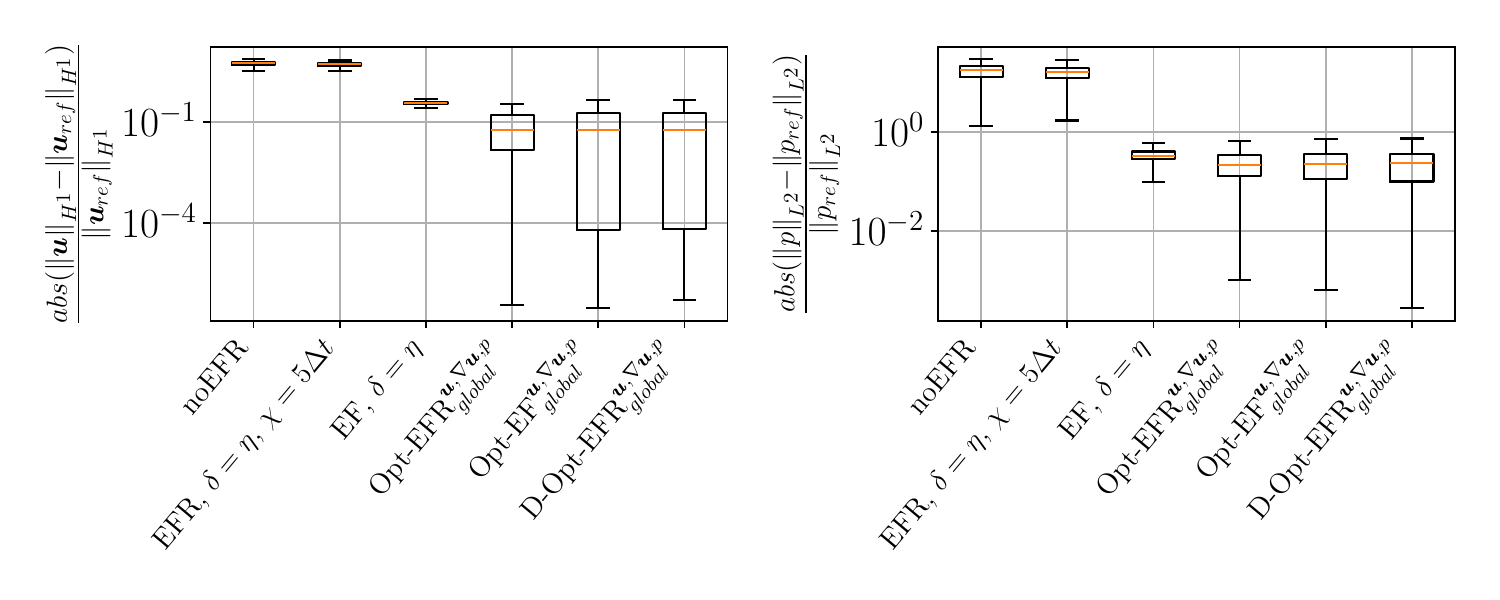}}\\
    \subfloat[$L^2$ pressure relative errors statistics]{\includegraphics[width=0.8\textwidth]{images/boxplot_double.pdf}}
    \caption{\RA{Statistics (median and confidence intervals) for noEFR, standard EFR/EF, and the methods of type $(\cdot)^{\bu, \nabla \bu, p}_{\text{global}}$. The results display the relative errors with respect to the DNS reference, and the plots show the time-averaged error (orange line) and the confidence interval (boxes). Specifically, Figure (A) displays the velocity $H^1$ norm discrepancy, while Figure (B) displays to the pressure $L^2$ norm results.
    The statistics are computed within the time window $[0, 4]$.}}
    \label{fig:boxplots}
\end{figure}

\begin{figure}[htpb!]
    \centering
    \includegraphics[width=0.9\textwidth]{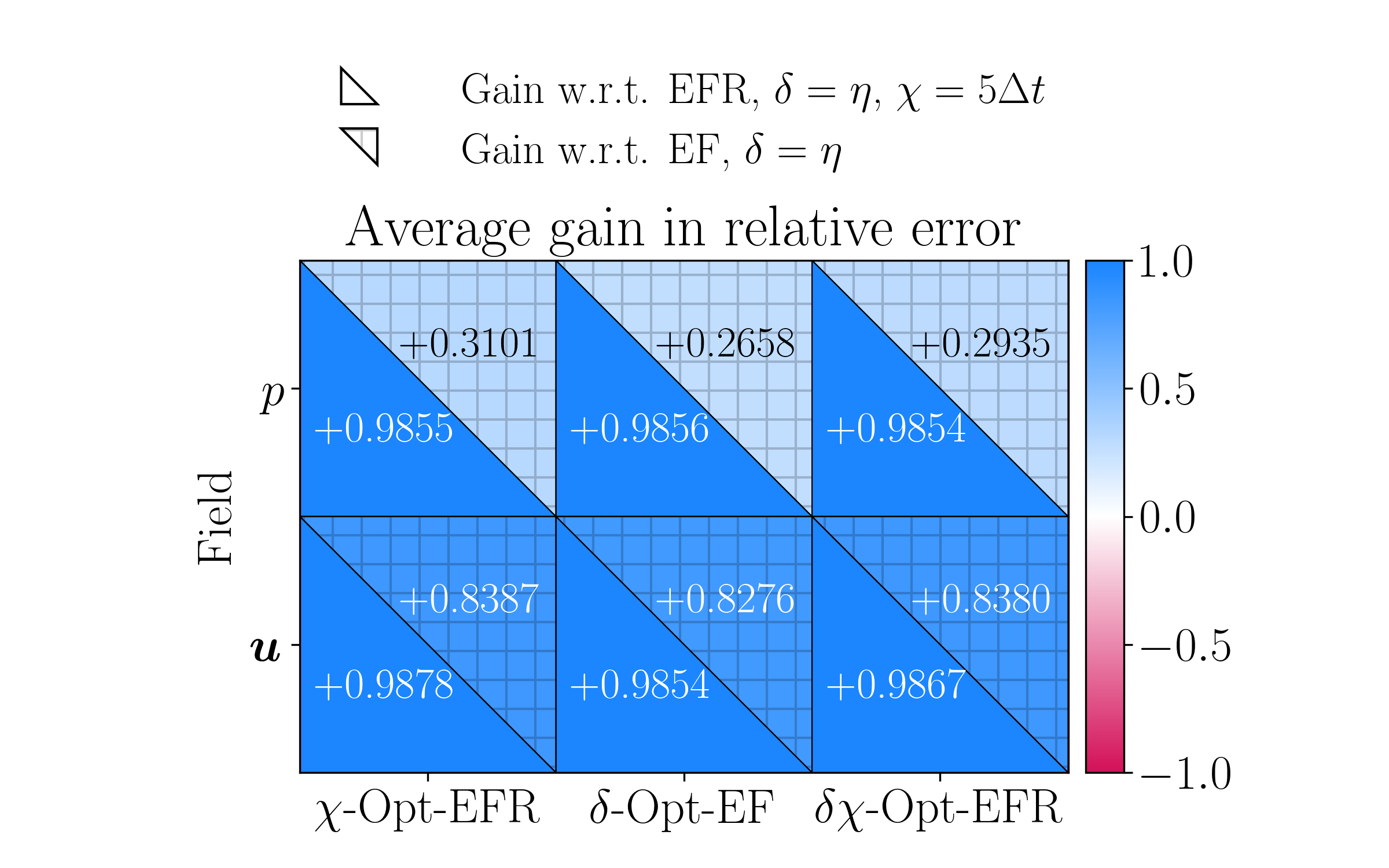}
    \caption{\RA{Gain in the relative errors in the $H^1$ velocity norm and $L^2$ pressure norm. The results show the time average gain of the best Opt-EFR methods with respect to the standard EFR and standard EF simulations.}}
    \label{fig:gains}
\end{figure}

As an overview of the results, we include two additional Figures:
\begin{itemize}
    \item Figure \ref{fig:boxplots} represents the time-average (orange lines) and confidence interval (boxes) 
    of the relative errors in the velocity $H^1$ norm and pressure $L^2$ norm.
    It confirms a significant improvement in the Opt-EFR strategies with respect to the baseline EFR/EF.
    \item Figure \ref{fig:gains} displays the gain increase in the accuracy of Opt-EFR results with respect to the baselines (standard EFR/EF).
    We measure the ``\emph{gain}" metric as the time-averaged relative error of the Opt-EFR solutions with respect to standard EFR and standard EF results.
    We conclude that Opt-EFR strategies significantly enhance the accuracy of standard methods. In particular, we reach an average gain 
    of \textbf{99} $\boldsymbol{\%}$ with respect to standard EFR, and of \textbf{83} $\boldsymbol{\%}$ (velocity) and \textbf{27-30} $\boldsymbol{\%}$ (pressure) with respect to standard EF.
\end{itemize}

\remark{Algorithm \ref{alg:opt-all} can be adapted to various test cases and different discretization types, such as finite element and finite volume methods. As expected, the related results will be problem-dependent. We remark, however, that this paper aims to provide a data-driven alternative to the standard EFR technique that enhances the simulation accuracy while maintaining a similar computational cost.}

%% file: sections/conclusions.tex
In this work, we make a significant contribution to the EFR regularization approach by proposing a \textbf{novel data-driven automated tuning of the filter and relaxation parameters}, $\delta$ and $\chi$, respectively.

We note that the literature \RA{provides} guidelines for the values of these parameters. \RA{However, existing approaches typically rely on: (i) constant parameter values; and (ii) heuristic choices or problem-dependent tuning, lacking a systematic framework for optimal selection.}
 \RA{ Our approach introduces an optimization-based methodology (the \textbf{Opt-EFR}) to dynamically adapt these parameters over time, improving the robustness and accuracy of EFR simulations. This novelty is especially relevant for under-resolved flow regimes, where fixed parameters often lead to suboptimal results.} 

\RA{In particular}, we propose and investigate three novel \RA{data-driven} optimization algorithms:
\begin{itemize}
    \item $\chi$-Opt-EFR, optimizing $\chi$ for fixed $\delta$; 
    \item $\delta$-Opt-EF, optimizing $\delta$ for fixed $\chi=1$;
    \item $\delta \chi$-Opt-EFR, optimizing both $\delta$ and $\chi$.
\end{itemize}

In each case, we investigate different objective functions, minimizing the discrepancy between the outcome of the EFR simulations (in an under-resolved regime) and a reference DNS (i.e., noEFR performed on a resolved grid).


We perform the numerical simulations on the test case of incompressible unsteady flow past a cylinder at a fixed Reynolds number $Re=1000$.

The \RA{key new findings of the proposed study are:}
\begin{enumerate}
    \item The optimized-EFR simulations always provide better results than noEFR and \RA{standard} EFR with fixed parameters.
    \item The best results are achieved by including in the objective function the discrepancies of the velocity \textbf{gradients} with respect to the reference simulation.
    \item The \textbf{global objective functions} always yield more accurate results than the corresponding local objective functions.
    \item \RA{The Opt-EFR simulations have a \textbf{similar or smaller computational time} than the standard EFR, even if they include an additional optimization step. This happens because the choice of optimal parameters improves the \textbf{numerical stability} of the system and allows for faster convergence.}
\end{enumerate}

\RA{
Despite these promising results, we note as a limitation the \textbf{reliance on data} for the optimization, which requires pre-computed DNS data to evaluate the discrepancy between the optimized EFR and the DNS reference. 
Additionally, the approach is \textbf{problem-dependent}, as the algorithm needs to be re-executed in different flow configurations.
}

As a possible future research direction, the authors would like to consider the application of machine learning approaches, such as \textbf{reinforcement learning}, in order to perform the training phase of the new optimization algorithms separately from the simulation. In this way, the training is performed only once and employs the data from a specific time window, allowing the prediction of the solution in a time extrapolation setting. \RA{This approach could significantly reduce the computational overhead and enable real-time parameter adaptation.} 

Another possible outlook for the project is to extend the new optimization strategies to the framework of reduced-order models by comparing the optimal full-order model parameters with the optimal reduced-order model parameters. 

\RA{We also highlight the \textbf{flexibility} of the Opt-EFR algorithms, that can be trivially extended to other flow scenarios, such as turbulent channel flows or aerodynamics simulations, where adaptive regularization plays a crucial role in improving prediction accuracy.}

%% file: sections/extra-final-sections.tex
\section*{Acknowledgments}
This work was partially funded by INdAM-GNCS: Istituto Nazionale di Alta Matematica –– Gruppo Nazionale di Calcolo Scientifico, and by European Union’s Horizon 2020 research and innovation programme under the Marie Skłodowska-Curie Actions, grant agreement 872442 (ARIA).
MS and GR thank the ``20227K44ME - Full and Reduced order modelling of coupled systems: focus on non-matching methods and automatic learning (FaReX)" project – funded by European Union – Next Generation EU  within the PRIN 2022 program (D.D. 104 - 02/02/2022 Ministero dell’Università e della Ricerca).
GR acknowledges the consortium iNEST (Interconnected North-East Innovation Ecosystem), Piano Nazionale di Ripresa e Resilienza (PNRR) – Missione 4 Componente 2, Investimento 1.5 – D.D. 1058 23/06/2022, ECS00000043, supported by the European Union's NextGenerationEU program.
MS thanks the ECCOMAS EYIC Grant ``CRAFT: Control and Reg-reduction in Applications
for Flow Turbulence". MG thanks the PON “Research and Innovation on Green related
issues” FSE REACT-EU 2021 project.
TI was funded by the National Science Foundation through grant DMS-2012253.
This manuscript reflects only the authors’ views and opinions and the Ministry cannot be considered responsible for them.

\subsection*{Author contributions}
Conceptualization: Anna Ivagnes, Maria Strazzullo, Traian Iliescu; Methodology: Anna Ivagnes, Maria Strazzullo, Traian Iliescu; Formal analysis and investigation: Anna Ivagnes; Writing - original draft preparation: Anna Ivagnes, Maria Strazzullo, Traian Iliescu; Writing - review and editing: Maria Strazzullo, Michele Girfoglio, Traian Iliescu; Funding acquisition: Gianluigi Rozza; Supervision: Traian Iliescu, Gianluigi Rozza.

\subsection*{Financial disclosure}

None reported.

\subsection*{Conflict of interest}
No potential competing interest was reported by the authors.

%% file: sections/supplementary.tex
\section{\RA{Methodological details on Opt-EFR algorithms}}
\RA{This part of the supplementary material specifies additional details on the $\chi$-Opt-EFR and $\delta$-Opt-EF methodologies.}
\label{sec:supp-methodology}

\subsection{\RA{The $\chi$-Opt-EFR algorithms: finding an optimal relaxation parameter}}
\label{sec:ml-efr}
\input{sections/methods-Opt-EFR}

\subsection{\RA{The $\delta$-Opt-EF algorithms: finding an optimal filter parameter}}
\label{sec:ml-ef}
\input{sections/methods-Opt-EF}

\subsection{\RA{The $\delta \chi$-Opt-EFR algorithms: finding the optimal filter and relaxation parameters}}
\label{sec:d-ml-efr}
\input{sections/methods-D-Opt-EFR}

\section{\RA{Additional numerical results}}
\label{sec:supp-results}

\RA{This part collects additional numerical results for the Opt-EFR strategies.
Specifically, we include the following comparative studies:
\begin{itemize}
    \item Algorithm $\chi$-Opt-EFR (section \ref{subsec:supp-results-chi}):
    \begin{itemize}
        \item Comparison of different $\chi$-Opt-EFR algorithms in the $L^2$ relative errors with respect to the DNS reference, highlighting the improvement with respect to standard EFR;
        \item A sensitivity study on how the optimization step influences the accuracy and computational time of simulations (section \ref{subsubsec:supp-pareto-chi}).
    \end{itemize}
    \item Algorithm $\delta$-Opt-EF (section \ref{subsec:supp-results-delta}):
    \begin{itemize}
        \item An analysis on how the results change by including the incompressibility constraint into the differential filter (section \ref{subsubsec:supp-incompressibility});
        \item A sensitivity study on how the optimization step influences the accuracy and computational time of simulations (section \ref{subsubsec:supp-pareto-delta}).
    \end{itemize}
\end{itemize}
}

\subsection{\RA{$\chi$-Opt-EFR: Additional results}}
\label{subsec:supp-results-chi}

\begin{figure*}[htpb!]
    \centering
    \includegraphics[width=0.75\textwidth]{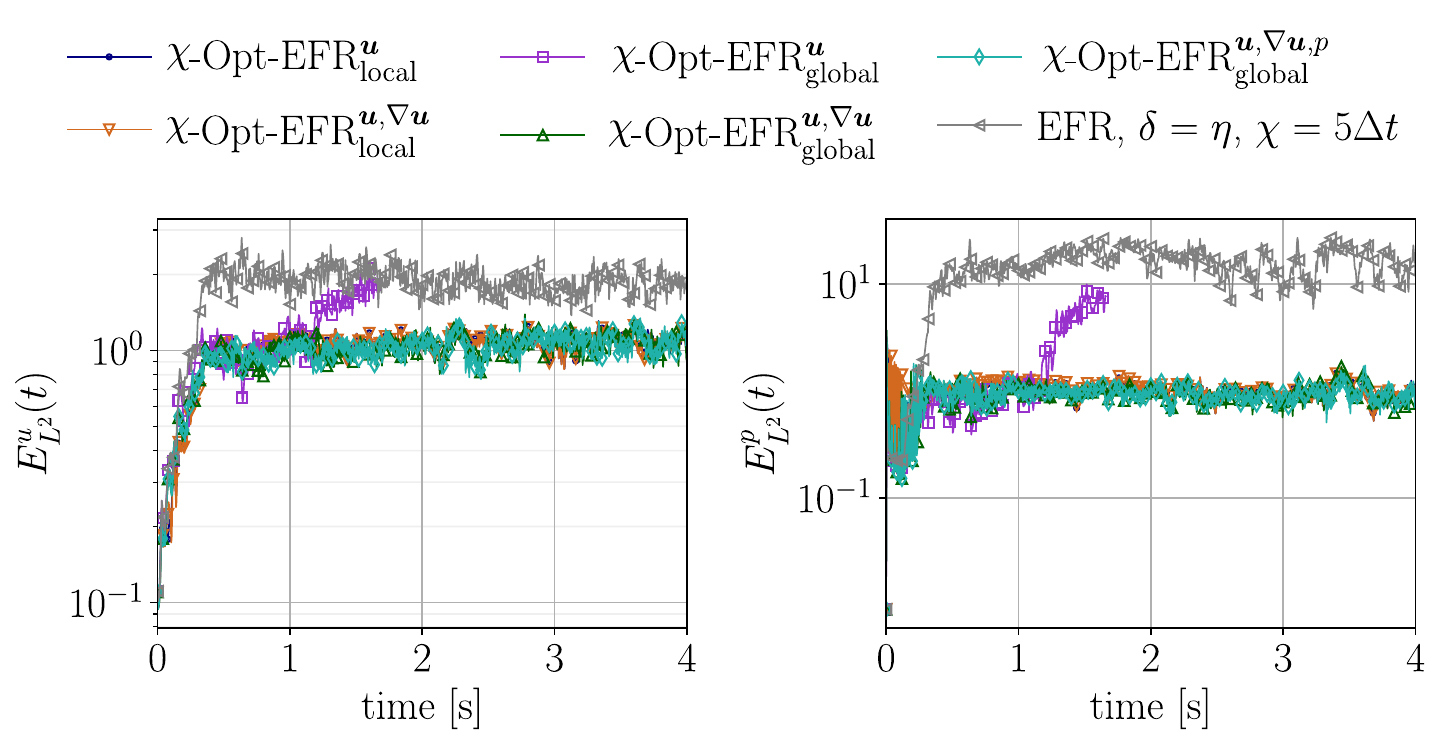}
    \caption{Relative velocity and pressure errors of $\chi$-Opt-EFR simulations with respect to the reference solutions 
    projected on the coarse mesh. The results for the standard EFR approach are also reported.}
    \label{fig:errs-chi}
\end{figure*}

We proceed with our analysis by investigating the accuracy with respect to the $L^2$ relative errors, i.e., $E_{L^2}^u(t)$ and $E_{L^2}^p(t)$, displayed in Figure \ref{fig:errs-chi}.
All the results look similar in the $L^2$ errors, except for the \optEFRglob{}, whose solution blows up at second $1.6$. Except for procedure \optEFRglob{}, all the optimized techniques yield more accurate results with respect to standard EFR (gray line in Figure \ref{fig:errs-chi}).

However, as we can see in the graphical results of Figures \ref{fig:opt-efr-u} and \ref{fig:opt-efr-p}, the $L^2$ relative error is not a significant metric to compare the different approaches in our test case. Indeed, even if the results for that metric are all similar, the flow dynamics are well captured in simulations of \optEFRglobGRAD{} and \optEFRglobPRESS{}, but not in the other approaches. In particular:
\begin{itemize}
    \item In \optEFRglob{}, the field is characterized by spurious oscillations in the whole domain;
    \item The local methods \optEFRloc{} and \optEFRlocGRAD{} are over-diffusive.
\end{itemize}

Hence, we will not consider the $L^2$ errors as metrics in the remaining part of the manuscript, but we will consider the $L^2$ velocity and pressure norms, and the $H^1$ velocity seminorm to evaluate the results.
We also not that the reason for analyzing \emph{global} features, instead of \emph{local} errors, is that, considering the optimization of local quantities, there is a phase shift in capturing the dynamics of the vortex shedding. Moreover, the local optimization is not able to fully recover all the vortices appearing in the reference DNS simulation.

\newpage

\subsubsection{On the role of the optimization step}
\label{subsubsec:supp-pareto-chi}




\RA{This section intestigates the role of the optimization step on the efficiency and accuracy of the $\chi$-Opt-EFR solution.}

\begin{figure*}[htpb!]
    \centering
    \subfloat[Discrepancy in the velocity $L^2$ norm]{
    \includegraphics[width=1\textwidth]{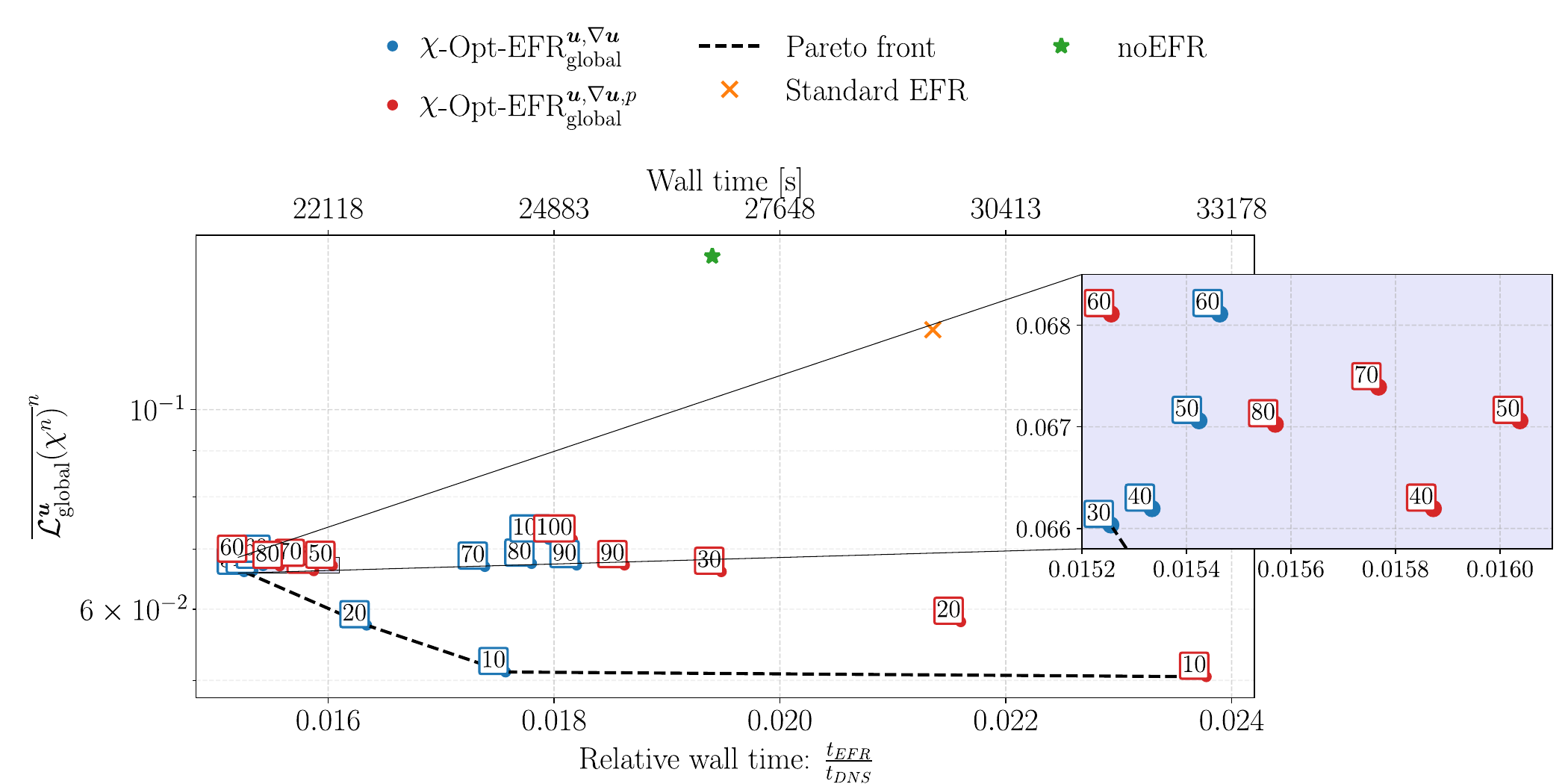}}\\
    \subfloat[Discrepancy in the pressure $L^2$ norm]{
    \includegraphics[width=1\textwidth]{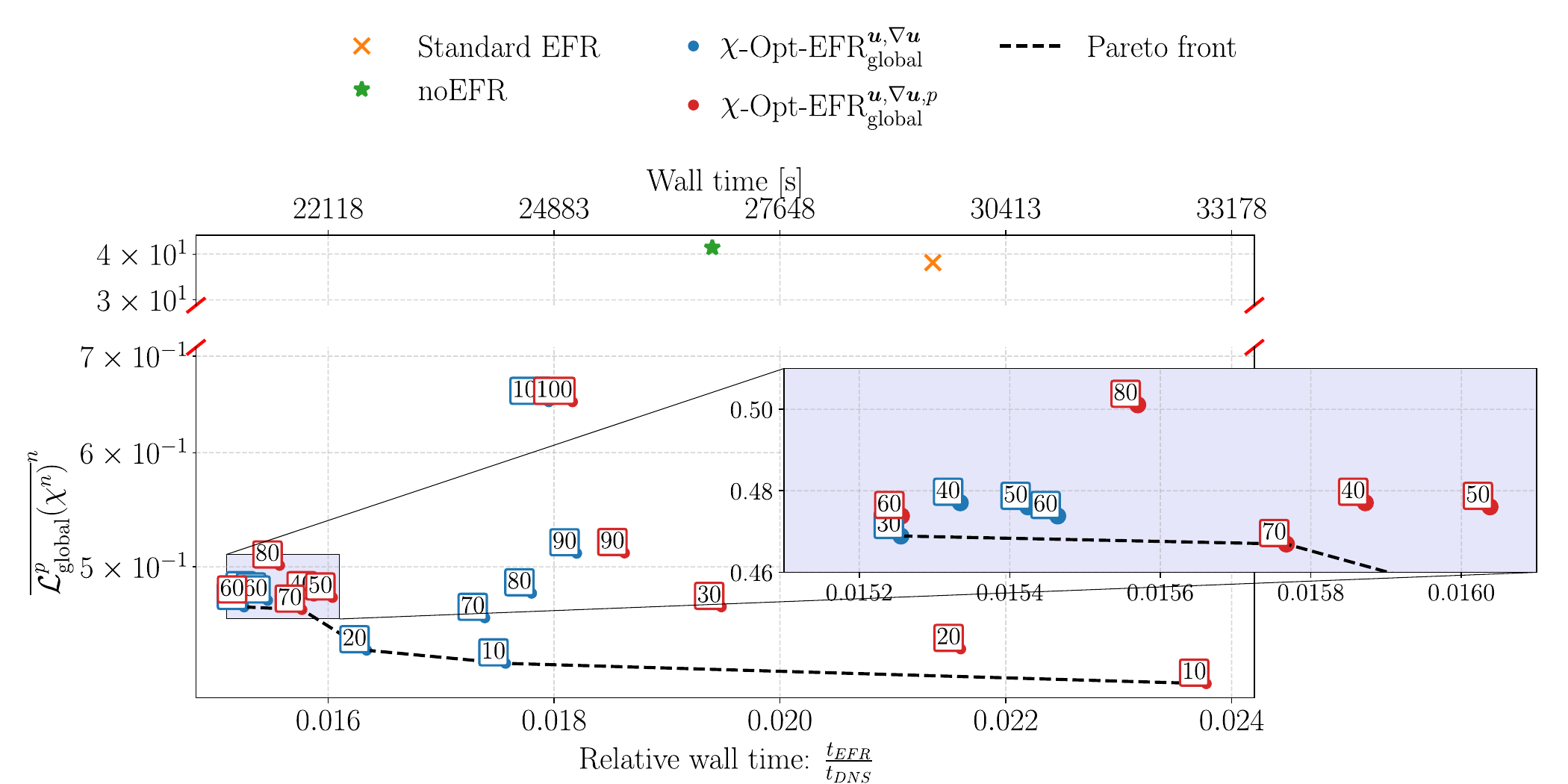}}
    \caption{Representation of how the optimization step affects the efficiency and the accuracy of the final velocity or pressure approximation in $\optEFRglobGRAD{}$ and $\optEFRglobPRESS{}$ simulations. The efficiency, on the $x$-axis, is measured through the relative wall clock time with respect to the DNS time. The accuracy, on the $y$-axis, is measured with respect to the $L^2$ velocity and pressure norms..}
    \label{fig:opt-efr-opt-step}
\end{figure*}

The results of the $\chi$-Opt-EFR simulations with different $k$ are presented in the Pareto plots in Figure \ref{fig:opt-efr-opt-step}, where the $k$ value is visualized in the boxes next to the corresponding markers.

\begin{figure*}[htpb!]
    \centering
        \subfloat[Discrepancy on the velocity $H^1$ semi-norm]{
    \includegraphics[width=1\textwidth]{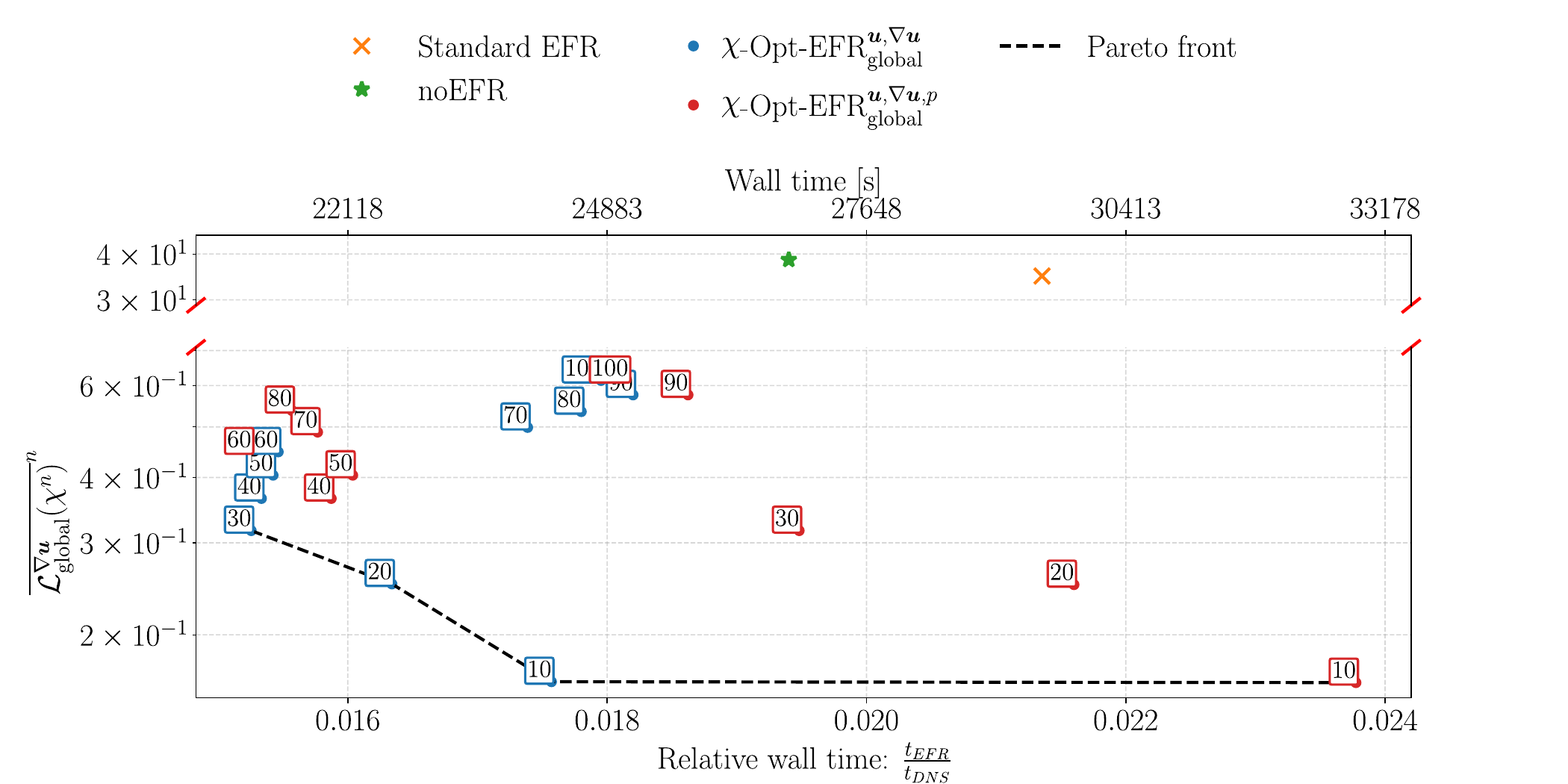}}
    \caption{Representation of how the optimization step affects the efficiency and the accuracy of the final velocity or pressure approximation in $\optEFRglobGRAD{}$ and $\optEFRglobPRESS{}$ simulations. The efficiency, on the $x$-axis, is measured through the relative wall clock time with respect to the DNS time. The accuracy, on the $y$-axis, is measured with respect to the $H^1$ velocity seminorm.}
    \label{fig:opt-efr-opt-step-h1}
\end{figure*}

Figures \ref{fig:opt-efr-opt-step} and \ref{fig:opt-efr-opt-step-h1} yield \RA{similar conclusions as the $\delta \chi$-Opt-EFR (section \ref{subsec:opt-step-final}). Specifically, all the optimized simulations show an increased accuracy with respect to standard EFR, with a similar computational time.}

\subsection{\RA{$\delta $-Opt-EF: Additional results}}
\label{subsec:supp-results-delta}

\subsubsection{On the incompressibility constraint}
\label{subsubsec:supp-incompressibility}

In EF simulations of incompressible flows, the incompressibility constraint is in general not satisfied when applying the differential filter in Step (II).
In particular, the filter preserves the above-mentioned constraint under periodic boundary conditions, but not under no-slip boundary conditions~\cite{bertagna2016deconvolution,layton2008numerical,van2006incompressibility}.
A technique that can be exploited to tackle this problem is the so-called \textbf{grad-div stabilization}. This approach consists in penalizing the violation of the incompressibility constraint by adding in the differential filter equations a term of the form:
\begin{equation}
\gamma \nabla (\nabla \cdot \overline{\boldsymbol{w}^{n+1}}),
    \label{graddiv}
\end{equation}
with $\gamma=100$, as in \cite{heavner2017locally}.
The reader can find more details on the grad-div stabilization in \cite{decaria2018determination, john2017divergence, layton2009accuracy}.

As stated in \cite{strazzullo2022consistency}, EFR simulations with $\chi \neq 1$ are typically characterized by low values of the velocity divergence, and the grad-div stabilization is not needed. 
However, such stabilization becomes necessary in EF simulations, namely when $\chi = 1$.

Although the investigation in Section \ref{sec:results-delta-opt} does not employ the grad-div stabilization, next we highlight how the results would change in case it is considered.

\begin{figure}[htpb!]
    \centering
    \includegraphics[width=0.5\textwidth]{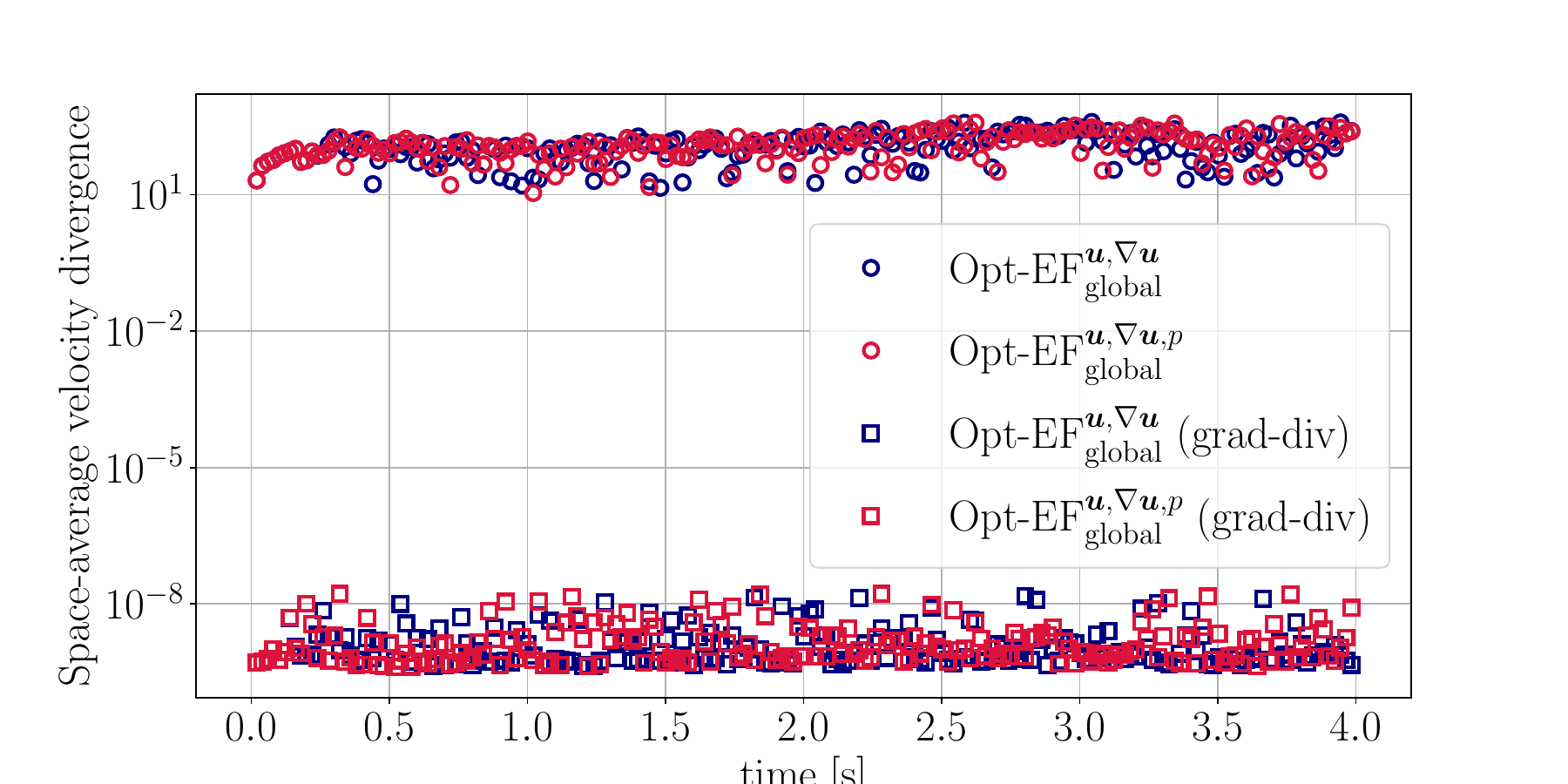}
    \caption{Variation in time of the space-averaged divergence of the velocity field, 
    for \optEFglobGRAD{} and \optEFglobPRESS{} with and without grad-div stabilization.}
    \label{fig:divergence-value}
\end{figure}

Figure \ref{fig:divergence-value} displays
the variation in time of the space-average of the velocity divergence 
for the \optEFglobGRAD{} and \optEFglobPRESS{}, with and without grad-div stabilization. It clearly shows that in the simulations without the grad-div stabilization the incompressibility constraint is not satisfied, while the grad-div stabilization fixes the issue.

\begin{figure*}[htpb!]
    \centering
    \includegraphics[width=0.9\linewidth]{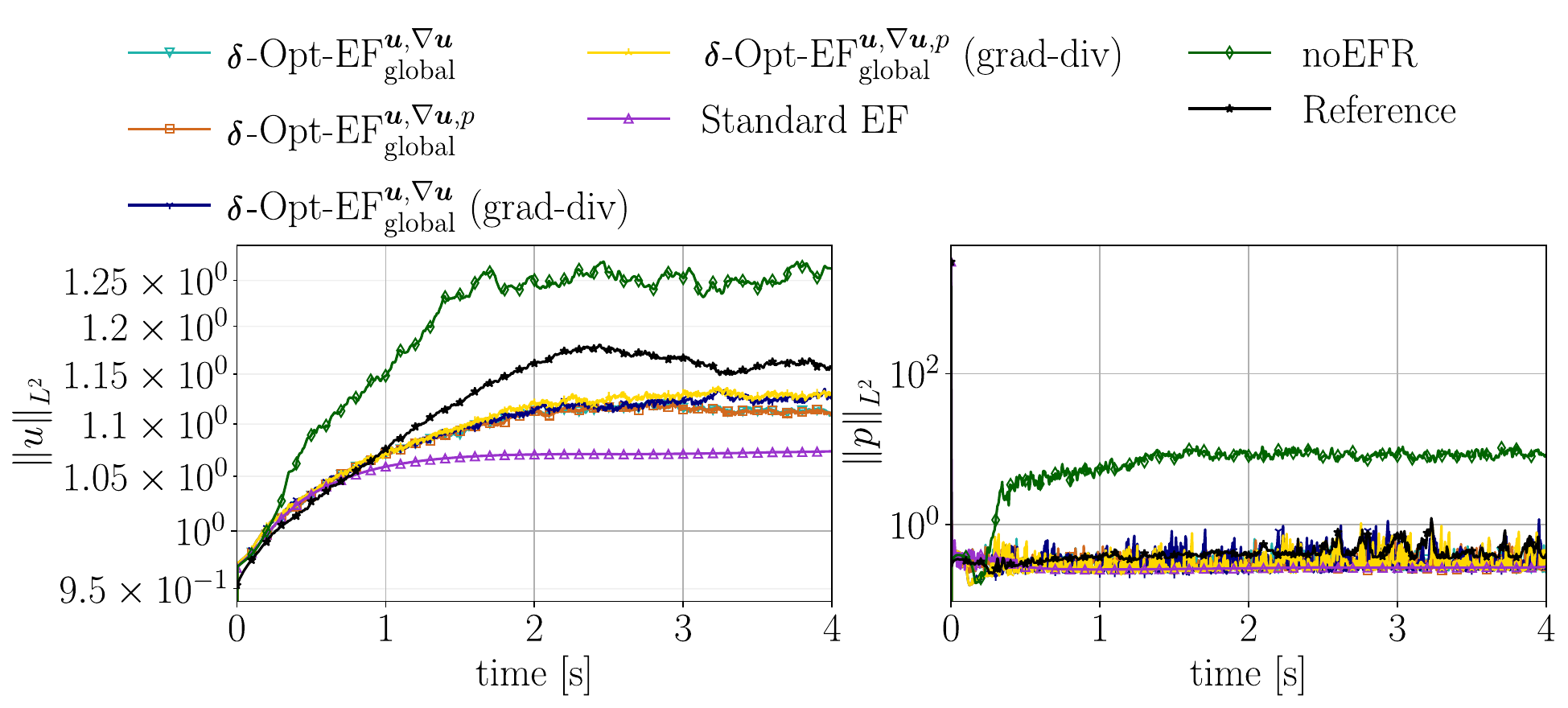}
    \caption{$L^2$ norms of the velocity and pressure fields in $\delta$-Opt-EF simulations, with and without 
    grad-div stabilization.}
    \label{fig:norm-l2-graddiv}
\end{figure*}

\begin{figure*}[htpb!]
    \centering
    \includegraphics[width=0.7\linewidth]{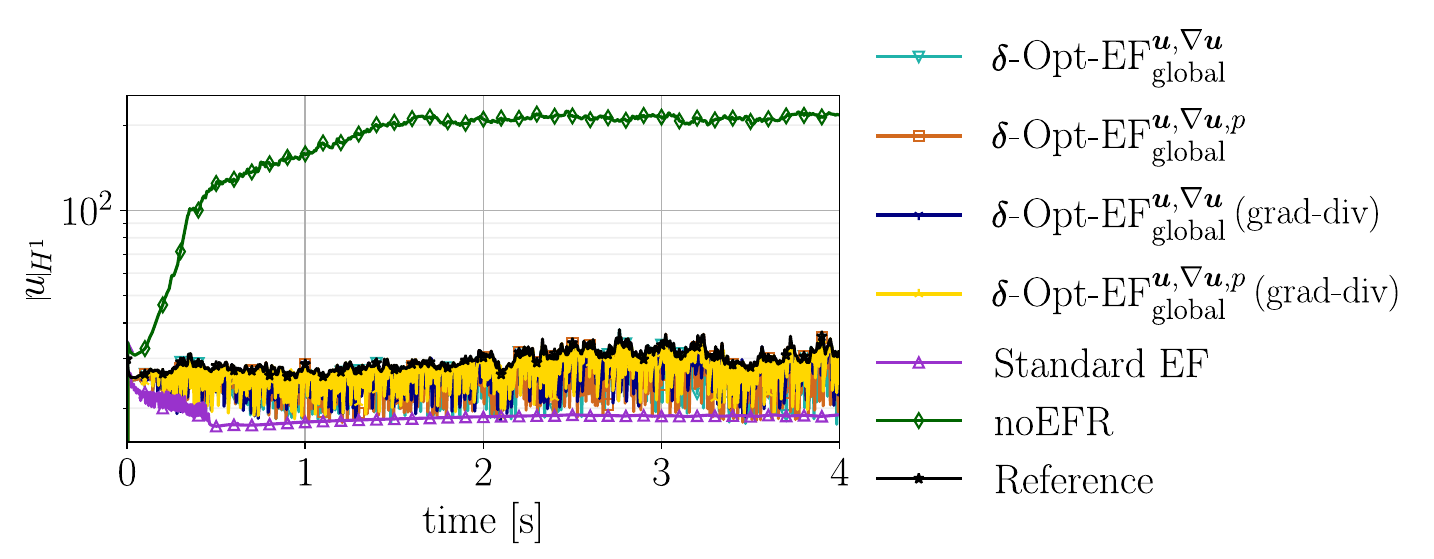}
    \caption{$H^1$ seminorm of the velocity field in $\delta$-Opt-EF simulations, with and without 
    grad-div stabilization.}
    \label{fig:norm-h1-graddiv}
\end{figure*}

Figures \ref{fig:norm-l2-graddiv} and \ref{fig:norm-h1-graddiv}, similarly to Figures \ref{fig:norms-delta} and \ref{fig:norm-delta-h1}, display the time evolution of the velocity and pressure $L^2$ norms, and of the velocity $H^1$ seminorm, respectively.
The plots suggest that the results with the grad-div stabilization are similar to the results obtained without stabilizing the differential filter. In particular, Figure \ref{fig:norm-l2-graddiv} shows that the velocity norm is slightly closer to the reference norm in the last time instances of the simulation, while the pressure maintains the same chaotic behavior observed for the non-stabilized simulations.


\begin{figure*}[htpb!]
    \centering
    \subfloat[DNS]{\includegraphics[width=0.5\textwidth, trim={3cm 7cm 3cm 7cm}, clip]{images/updated_Re1000/u_reference.png}}
    \\
    \subfloat[\optEFglobGRAD{} without grad-div]{\includegraphics[width=0.5\textwidth, trim={3cm 7cm 3cm 7cm}, clip]{images/updated_Re1000/u_delta_global_grad.png}}
    \subfloat[\optEFglobPRESS{} without grad-div]{\includegraphics[width=0.5\textwidth, trim={3cm 7cm 3cm 7cm}, clip]{images/updated_Re1000/u_delta_global_press.png}}\\
    
    \subfloat[\optEFglobGRAD{} with grad-div]{\includegraphics[width=0.5\textwidth, trim={6cm 10cm 6cm 15cm}, clip]{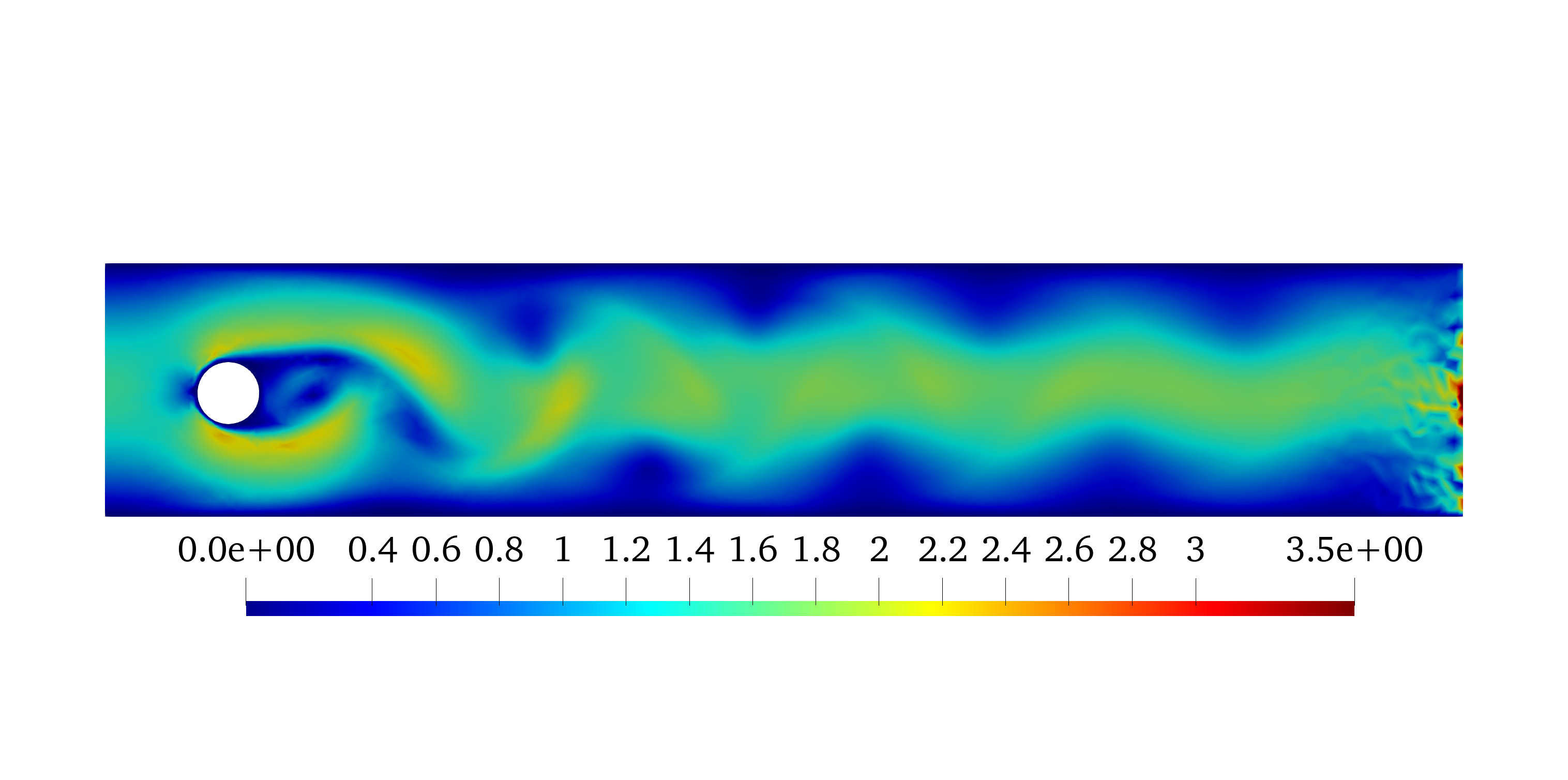}}
    \subfloat[\optEFglobPRESS{} with grad-div]{\includegraphics[width=0.5\textwidth, trim={6cm 10cm 6cm 15cm}, clip]{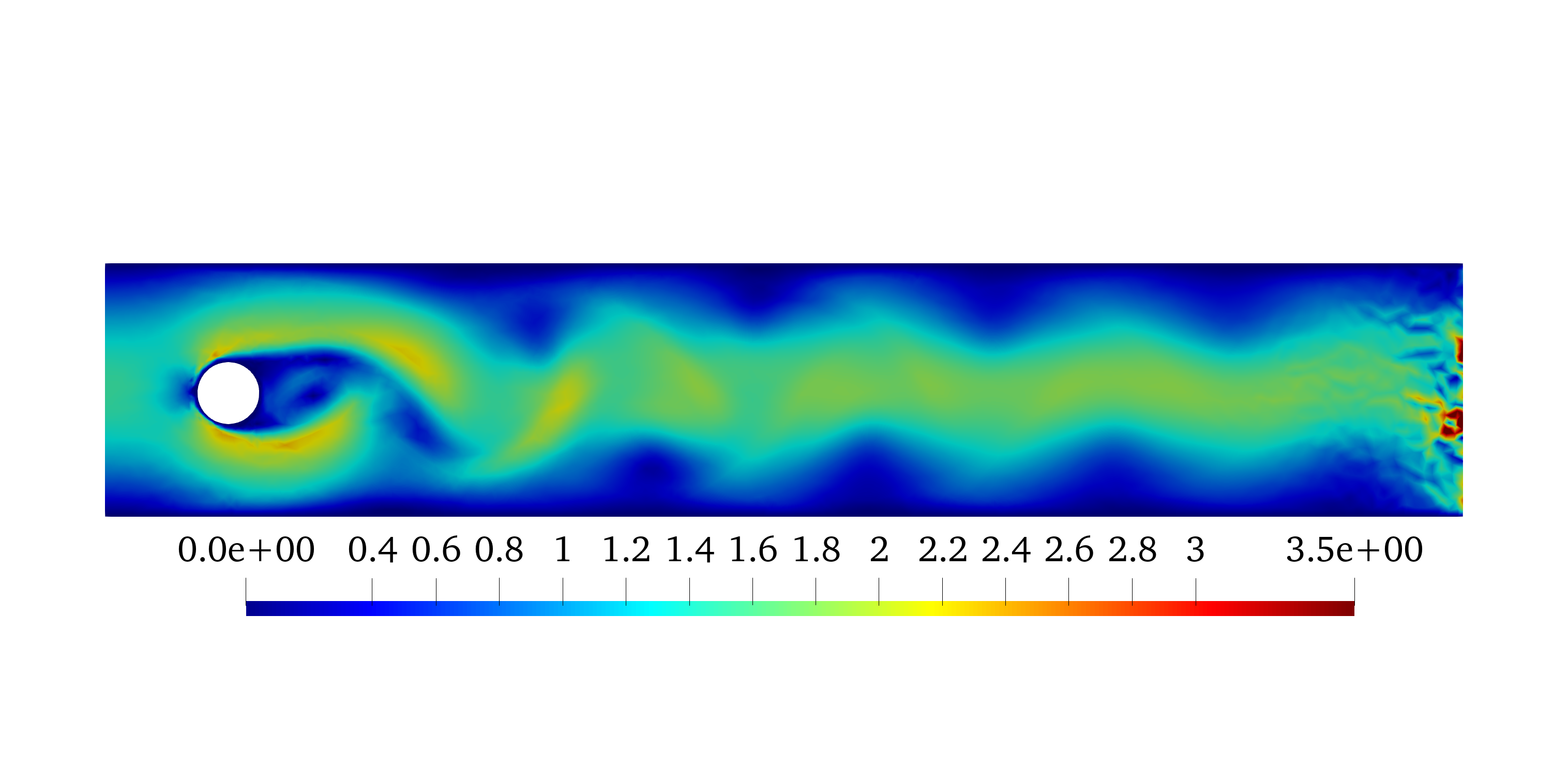}}
    
    \caption{Velocity fields for $\delta$-Opt-EF (with and without 
    grad-div stabilization) and reference DNS 
    at the final time instance ($t=4$).}
    \label{fig:graddiv-fields-vel}
\end{figure*}

\begin{figure*}[htpb!]
    \centering
    
    \subfloat[DNS]{\includegraphics[width=0.5\textwidth, trim={3cm 7cm 3cm 7cm}, clip]{images/updated_Re1000/p_reference.png}}\\
    \subfloat[\optEFglobGRAD{} without grad-div]{\includegraphics[width=0.5\textwidth, trim={3cm 7cm 3cm 7cm}, clip]{images/updated_Re1000/p_delta_global_grad.png}}
    \subfloat[\optEFglobPRESS{} without grad-div]{\includegraphics[width=0.5\textwidth, trim={3cm 7cm 3cm 7cm}, clip]{images/updated_Re1000/p_delta_global_press.png}}\\
    
    \subfloat[\optEFglobGRAD{} with grad-div]{\includegraphics[width=0.5\textwidth, trim={6cm 10cm 6cm 15cm}, clip]{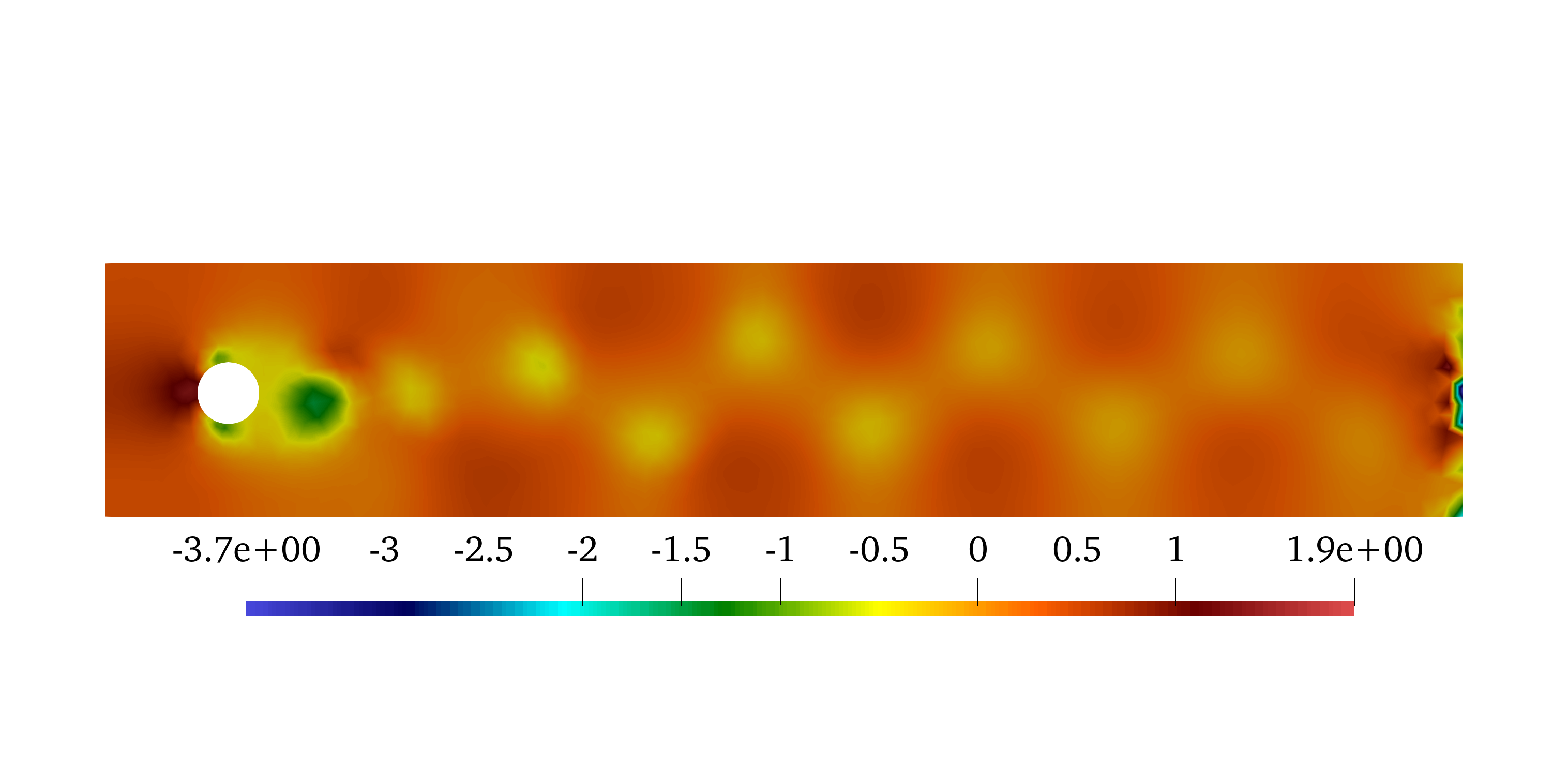}}
    \subfloat[\optEFglobPRESS{} with grad-div]{\includegraphics[width=0.5\textwidth, trim={6cm 10cm 6cm 15cm}, clip]{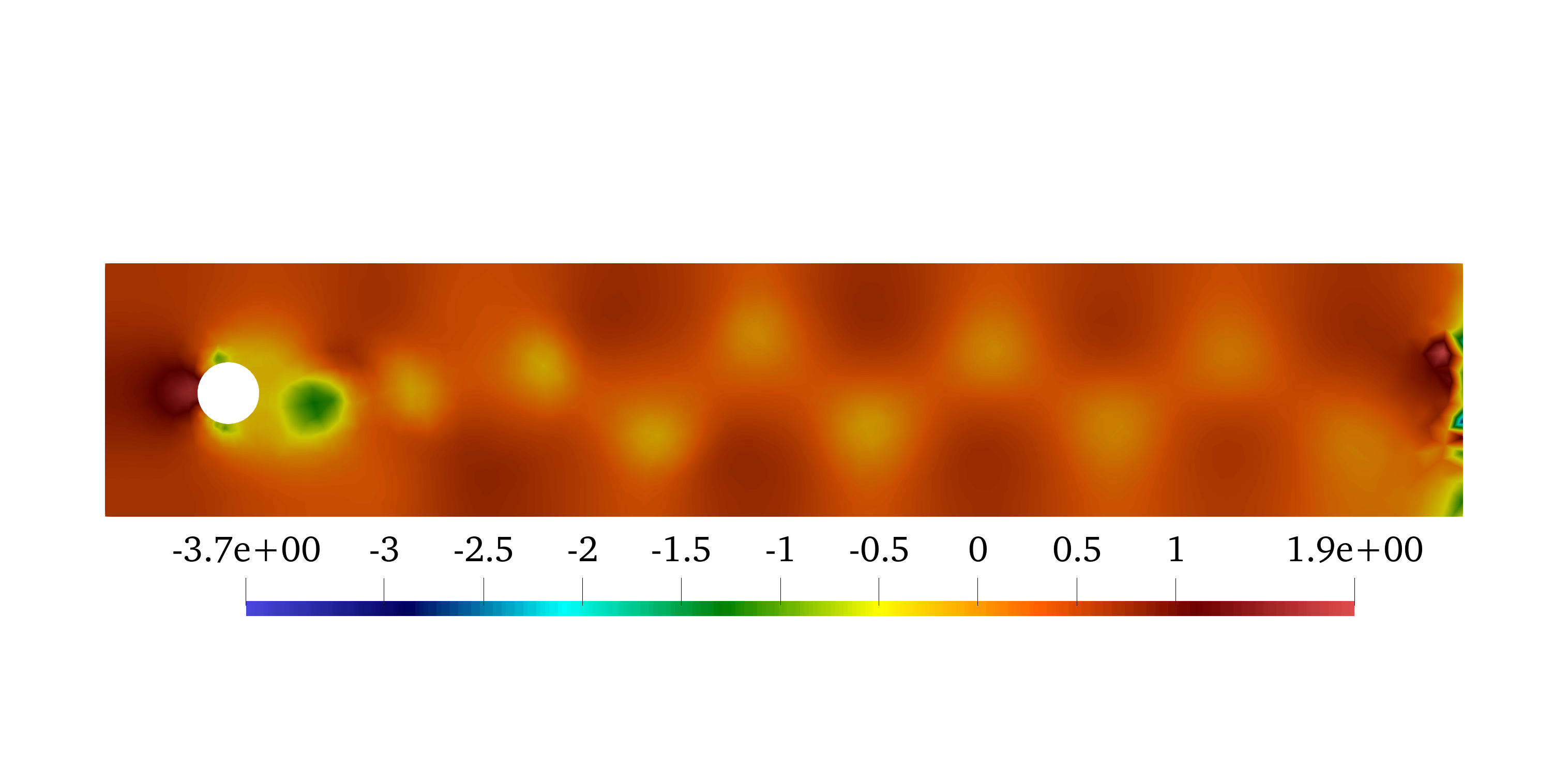}}
    
    \caption{Pressure fields for $\delta$-Opt-EF (with and without 
    grad-div stabilization) and reference DNS 
    at the final time instance ($t=4$).}
    \label{fig:graddiv-fields-press}
\end{figure*}

We also provide a comparison of the global performance of simulations with and without grad-div stabilization in Table \ref{tab:graddiv}. The table reports the time average of the global contributions defined in Equations \eqref{eq:global-contributions-delta}.
We can notice that the grad-div simulations provide a slightly improved velocity, as already illustrated in the left plot in Figure \ref{fig:norm-l2-graddiv}. However, the pressure accuracy is considerably lower.
Moreover, if we look at the graphical representations in Figures \ref{fig:graddiv-fields-vel} and \ref{fig:graddiv-fields-press}, we can notice that,
at the final time instance, the grad-div simulation and the DNS lift have opposite signs, while the no-grad-div simulations better capture the flow behavior.



%
\begin{table*}[h]
\centering
\caption{Comparison between simulations with and without grad-div stabilization in terms of accuracy.}
\label{tab:graddiv}
{
\begin{tabular}{ccccc}
\\ \toprule
&&$\overline{\mathcal{L}_{\mathrm{global}}^{\bu}(\delta^n)}^n$&$\overline{\mathcal{L}_{\mathrm{global}}^{\nabla \bu}(\delta^n)}^n$&$\overline{\mathcal{L}_{\mathrm{global}}^{p}(\delta^n)}^n$
\\
\midrule
\multirow{2}{*}{\textbf{Without} grad-div} & \optEFglobGRAD{}&$0.0578$&$0.188$&$0.420$
\\
\cline{2-5}
&\optEFglobPRESS{}&$0.0585$&$0.189$&$0.425$
\\
\midrule
\multirow{2}{*}{\textbf{With} grad-div} & \optEFglobGRAD{}&$0.0520$&$0.174$&$0.575$
\\
\cline{2-5}
&\optEFglobPRESS{}&$0.0469$&$0.167$&$0.560$
\\ \bottomrule
\end{tabular}}
\end{table*}

\newpage

\subsubsection{On the role of the optimization step}
\label{subsubsec:supp-pareto-delta}
As in Section \ref{subsubsec:supp-pareto-chi}, we include here an analysis of the accuracy depending on the optimization step.

\begin{figure*}[htpb!]
    \centering
    \subfloat[Discrepancy on the velocity $L^2$ norm]{
    \includegraphics[width=0.9\textwidth]{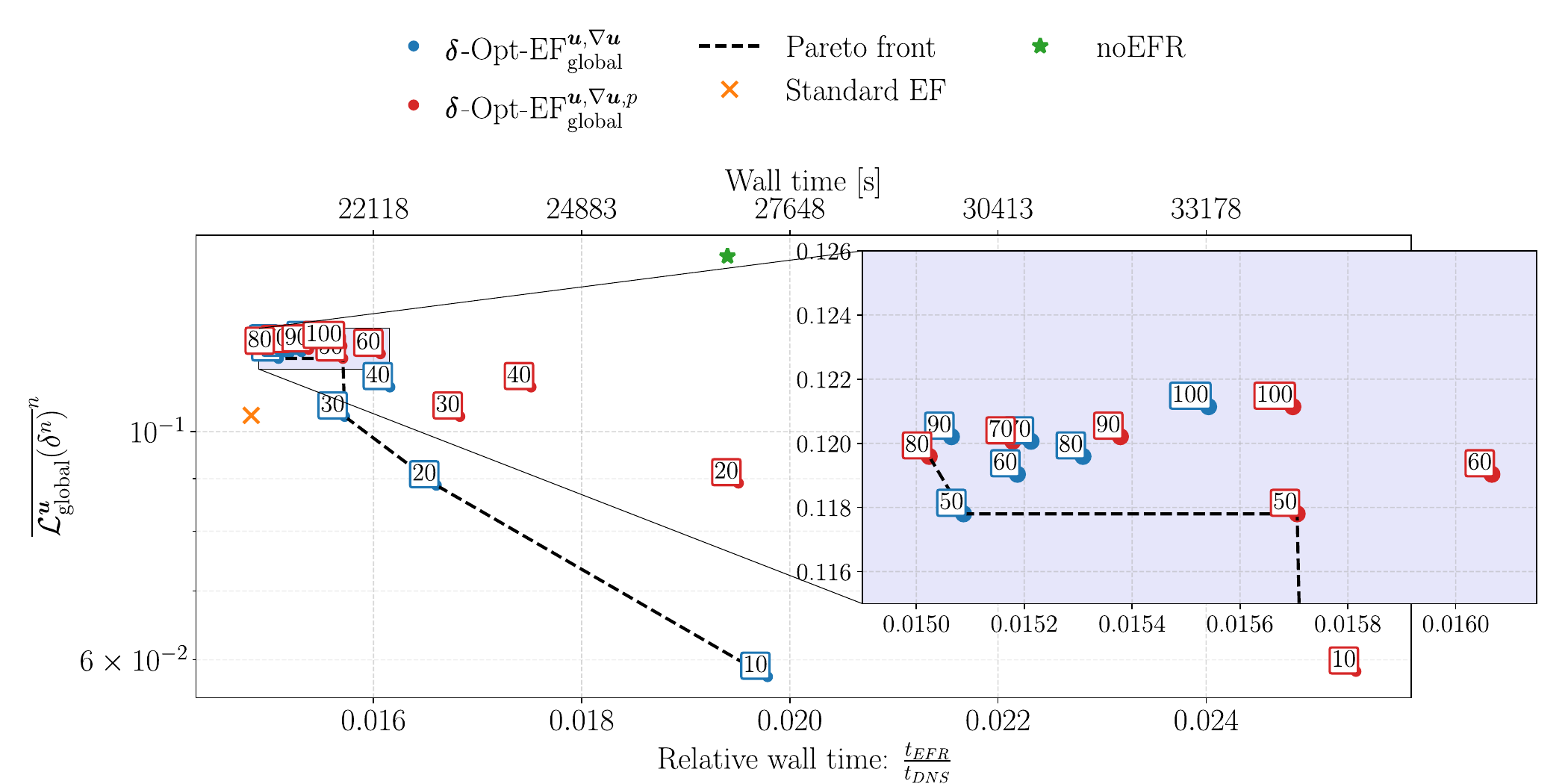}}\\
    \subfloat[Discrepancy on the pressure $L^2$ norm]{
    \includegraphics[width=0.9\textwidth]{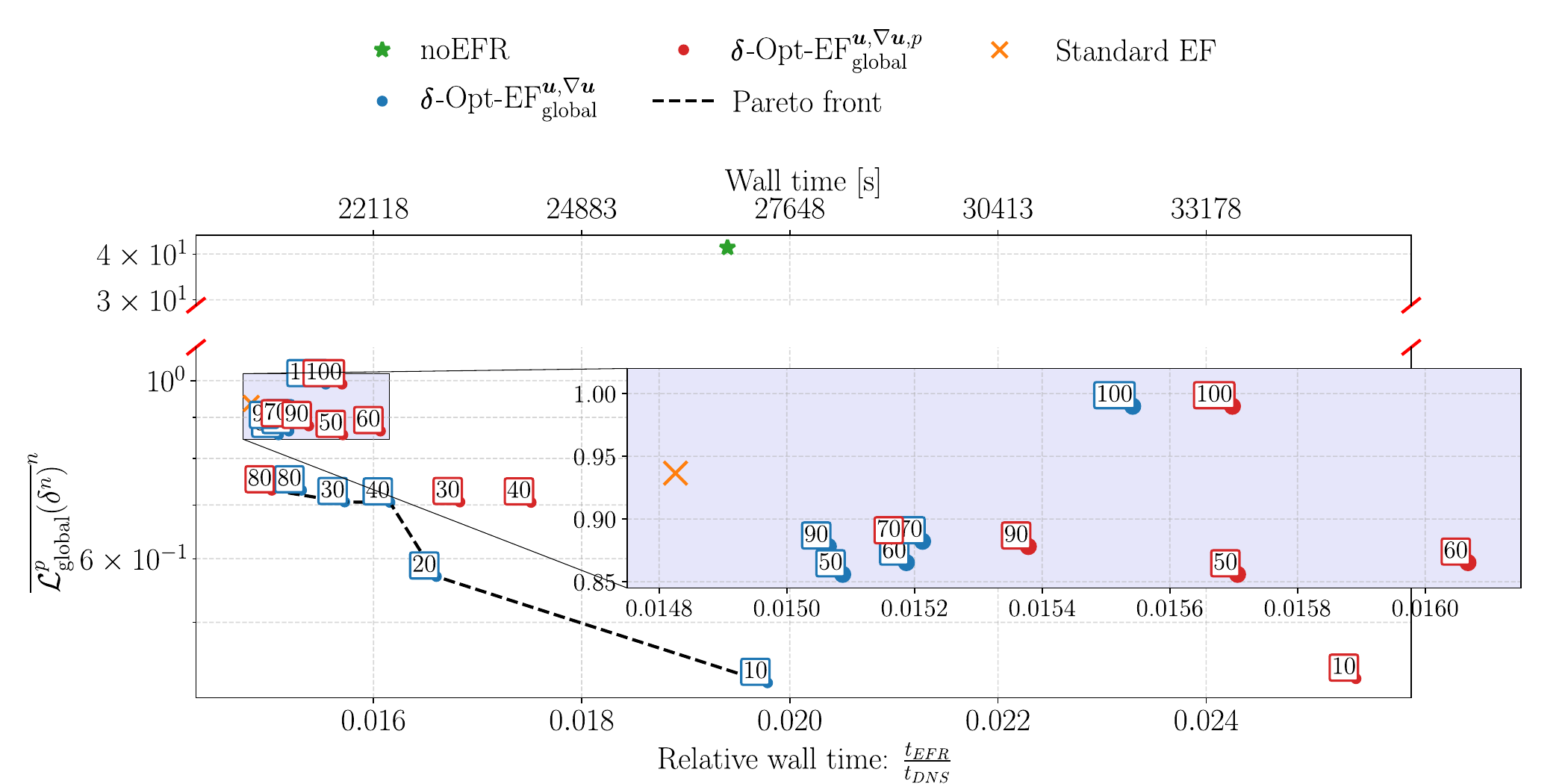}}
    \caption{Representation of how the optimization step affects the efficiency and the accuracy of the final velocity or pressure approximation in $\optEFglobGRAD{}$ and $\optEFglobPRESS{}$ simulations. The efficiency, on the $x$-axis, is measured through the relative wall clock time with respect to the DNS time. The accuracy, on the $y$-axis, is measured with respect to the $L^2$ velocity and pressure norms.}
    \label{fig:opt-ef-opt-step}
\end{figure*}

\begin{figure*}[htpb!]
    \centering
        \subfloat[Discrepancy on the velocity $H^1$ semi-norm]{\includegraphics[width=0.9\textwidth]{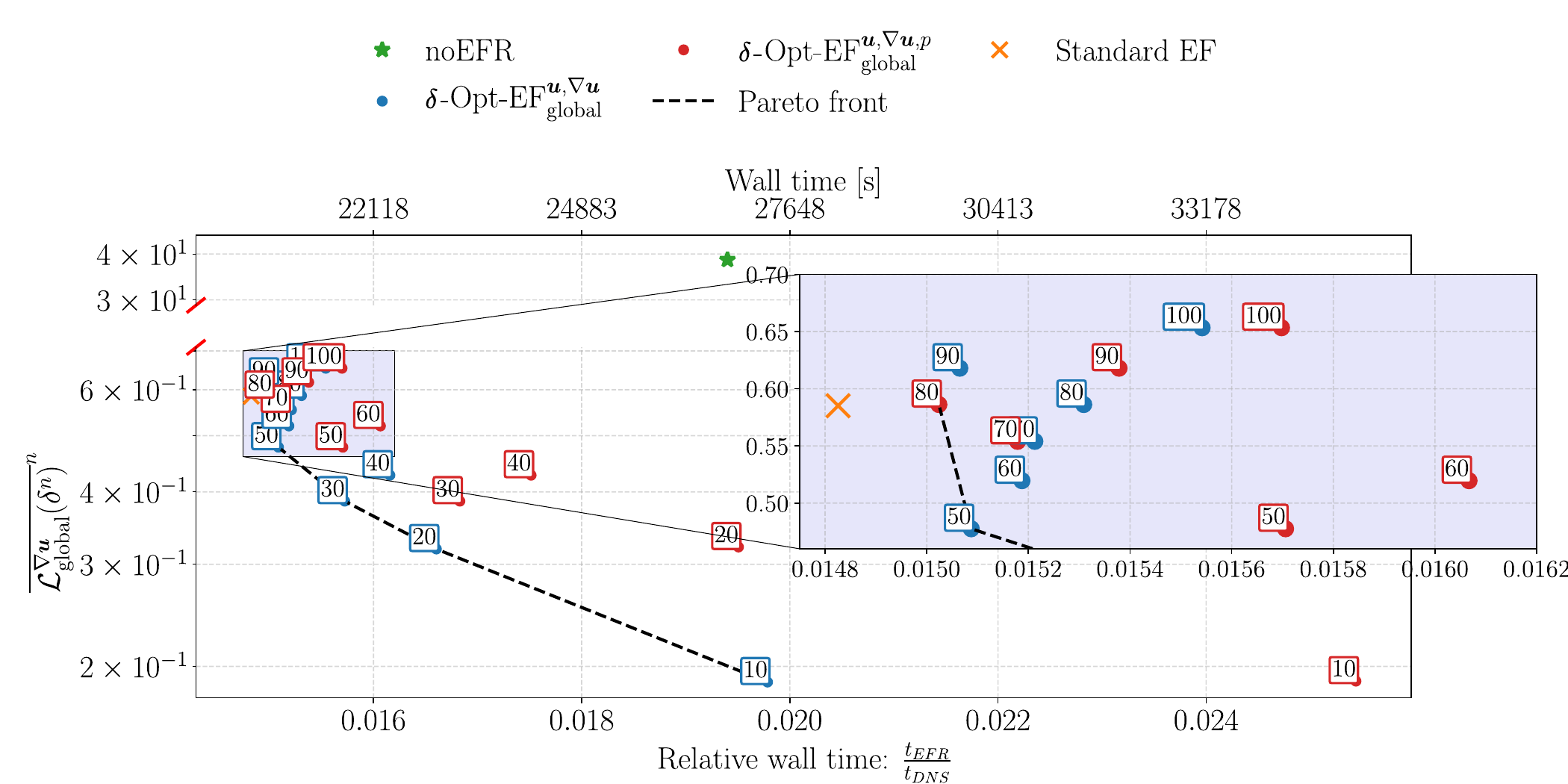}}
    \caption{Representation of how the optimization step affects the efficiency and the accuracy of the final velocity or pressure approximation in $\optEFglobGRAD{}$ and $\optEFglobPRESS{}$ simulations. The efficiency, on the $x$-axis, is measured through the relative wall clock time with respect to the DNS time.
    The accuracy, on the $y$-axis, is measured with respect to the $H^1$ velocity seminorm.}
    \label{fig:opt-ef-opt-step-h1}
\end{figure*}

The conclusions are the same as those outlined for the $\chi$-Opt-EFR algorithm.
Indeed, as can be seen from Figures \ref{fig:opt-ef-opt-step} and \ref{fig:opt-ef-opt-step-h1}, also in this case the algorithm including the pressure contribution (\optEFglobPRESS{}) is slightly more accurate, but the 
\optEFglobGRAD{} simulations are more efficient.
We can, however, 
notice some differences between the $\chi$-Opt-EFR and $\delta$-Opt-EF strategies:
\begin{itemize}
    \item[(i)] \textbf{The accuracy of the optimized EF/EFR.} Looking at the $y$-scale in Figures \ref{fig:opt-ef-opt-step} and \ref{fig:opt-efr-opt-step}, we can notice that the $\chi$-Opt-EFR algorithms provide in general more accurate results when keeping constant the filter parameter. Hence, optimizing only $\chi(t)$ and keeping the filter parameter fixed to the Kolmogorov scale produces lower errors than optimizing only $\delta(t)$ without the relaxation step.
    \item[(ii)] \textbf{The accuracy of the standard EF/EFR.} Standard EFR provides noisy results, while standard EF is over-diffusive but 
    is closer to the DNS results. It is indeed more accurate than the standard EFR with $\chi \neq 1$.
    As an additional consideration, since it is more stable and does not exhibit spurious oscillations, the Newton method converges faster to the FE solution, and its wall clock time is considerably lower than that of the standard EFR.
    \item[(iii)] \textbf{The comparison in time between standard EF and $\delta$-Opt-EF} The \optEFglobGRAD{} approach proves to be more accurate than the standard EF when considering $\Delta t_{opt}=10\Delta t$ and $\Delta t_{opt}=20\Delta t$. However, differently from what happens for $\chi$-Opt-EFR algorithms (Section \ref{sec:results-chi-opt}), 
    these methods are characterized by longer computational times than the standard EF. This is expected because EF is overdiffusive and, thus, facilitates the convergence of the Newton solver.
\end{itemize}

%% file: sections/methods-Opt-EFR.tex
This section \RA{describes} a strategy that combines standard optimization algorithms with the EFR strategy to find the optimal $\chi$, starting from a fixed value of the parameter $\delta$.

All the algorithms \RA{described} in this section predict a discretized time-dependent $\chi$, $\{\chi(t^n)\}_{n=1}^{N_T}$, minimizing an objective function.  
Specifically, we investigate several
objective functions to deduce how the 
different optimization approaches would affect the accuracy of the final EFR solution.
The features of the adopted strategies are summarized in Table \ref{tab:acronyms}.
\begin{table*}[h]
\caption{Acronyms and main features of the algorithms $\chi$-Opt-EFR. Gray cells indicate that the algorithm does not present that feature.}
\label{tab:acronyms}
\centering
{
\begin{tabular}{cccccc}
\\ \hline
\multirow{2}{*}{\textbf{Acronym}} & \multicolumn{2}{c}{\textbf{Objective function type}} & \multicolumn{3}{c}{\textbf{Contributions in objective function}}
\\ \cline{2-3}\cline{4-6}
& local & global & $\bu$ discrepancy & $\nabla \bu$ discrepancy & $p$ discrepancy \\ \hline
\optEFRloc & $\checkmark$ & \cellcolor{gray!50} & $\checkmark$ & \cellcolor{gray!50}& \cellcolor{gray!50}\\ \hline
\optEFRlocGRAD & $\checkmark$ & \cellcolor{gray!50} & $\checkmark$ & $\checkmark$ & \cellcolor{gray!50}\\ \hline
\optEFRglob & \cellcolor{gray!50} & $\checkmark$ & $\checkmark$ & \cellcolor{gray!50} & \cellcolor{gray!50}\\ \hline
\optEFRglobGRAD & \cellcolor{gray!50} & $\checkmark$ & $\checkmark$ &$\checkmark$& \cellcolor{gray!50}\\ \hline
\optEFRglobPRESS & \cellcolor{gray!50} & $\checkmark$ & $\checkmark$ &$\checkmark$ &$\checkmark$  \\ \hline
\end{tabular}}
\end{table*}

The algorithms in Table \ref{tab:acronyms} perform at each time instance Step (I), Step (II) for a given constant $\delta$ and, then, look for an optimal time-dependent $\chi(t)$ to be employed in Step (III). The process is summarized in Algorithm \ref{alg:opt-efr-chi}.


\RA{As mentioned in Section \ref{sec:methods-opt}, the objective function to minimize in order to compute the optimal $\chi^n=\chi(t^n)$ can
feature either}
\emph{local} or \emph{global} contributions.
The \RA{local} objective function \RA{is specialized here} as follows:

\begin{equation}
    \label{eq:local-loss}
\lossLOC(\chi^{n}) = w_{\bu} \lossLOCu(\chi^{n}) + w_{\nabla \bu} \lossLOCgrad(\chi^{n}), 
\end{equation}
where the individual contributions are 
\begin{equation}
\label{eq:local-contributions}
\begin{split}
    &\lossLOCu(\chi^{n}) = \text{MSE}\left( \bu^{n}(\chi^n), \bu_{ref}^{n} \right), \\
   &\lossLOCgrad(\chi^{n}) = \text{MSE}\left( \nabla \bu^{n}(\chi^n), \nabla \bu_{ref}^{n} \right).
   \end{split}
\end{equation}

The global objective function has the following expression:

\begin{equation}
    \label{eq:global-loss}
\lossGLOB(\chi^{n}) = w_{\bu} \lossGLOBu(\chi^{n}) + w_{\nabla \bu} \lossGLOBgrad(\chi^{n}) + w_{p} \lossGLOBp(\chi^{n}),
\end{equation}
where the global contributions are 
\begin{equation}
\label{eq:global-contributions}
\begin{split}
    &\lossGLOBu(\chi^n) = \left| \dfrac{\|\bu^{n}(\chi^n)\|^2_{L^2(\Omega)} - \|\bu_{ref}^{n}\|^2_{L^2(\Omega)}}{\|\bu_{ref}^{n}\|^2_{L^2(\Omega)}} \right|, \\
   & \lossGLOBgrad(\chi^{n}) =\left| \dfrac{\| \nabla \bu^{n}(\chi^n)\|^2_{L^2(\Omega)} - \| \nabla \bu_{ref}^{n}\|^2_{L^2(\Omega)}}{\| \nabla \bu_{ref}^{n}\|^2_{L^2(\Omega)}} \right|,\\
   &\lossGLOBp(\chi^{n}) = \left| \dfrac{\|p^{n+1}(\chi^n)\|^2_{L^2(\Omega)} - \|p_{ref}^{n+1}\|^2_{L^2(\Omega)}}{\|p_{ref}^{n+1}\|^2_{L^2(\Omega)}} \right|.
\end{split}
\end{equation}



\begin{algorithm*}[htpb!]
\caption{Pseudocode for algorithms $\chi$-Opt-EFR optimizing the relaxation parameter $\chi(t)$. Row 10 (underlined) is considered only for \optEFRglobPRESS{}. The optimization part is highlighted in gray.}
\label{alg:opt-efr-chi}
\begin{algorithmic}[1]
\State $\bu_{ref}^n$, $n=0, \dots, N_T - 1$; $N_h^u$; $\bu^0$, $\Delta t_{opt}$, $\chi^0 = 5\Delta t$ \Comment{$\text{\it Inputs needed}$}
\For{$n \in [0, \dots, N_T - 1]$}
\State (I)+(II) $\rightarrow \bw^{n+1}$, $\overline{\bw}^{n+1}$ \Comment{$\text{\it Perform Evolve and Filter steps}$}
\State $\bu_{ref}^{n+1} \gets \bu_{ref}^{n+1}\,  | \, \mathbb{U}^{N_h^u}$ \Comment{$\text{\it Project the reference field into the coarse mesh}$}
\State $\chi^{n+1} \gets \chi^n$ \Comment{$\text{\it Initialize } \chi$}
\highlight{\If {$\bmod(n, \Delta t_{opt}) = 0$} \Comment{$\text{\it Perform optimization every }\Delta t_{opt}$} 
    \State $iter =0$
\Repeat \Comment{$\text{\it SLSQP algorithm}$}
        \State (III) $\rightarrow \bu^{n+1}$ \Comment{$\text{\it Relax step}$}
        \State {\underline{(I) $\rightarrow p^{n+2}$} \Comment{$\text{\it \underline{(Second) Evolve step}}$}}
        \State Optimize $\mathcal{L}(\chi^{n+1})$, subject to $0 \leq \chi^{n+1} \leq 1$ 
        \State $iter \gets iter +1$
        \State Output current solution: $\chi^{n+1}$
    \Until{Convergence or maximum iterations reached}
\EndIf}
\State Solve (III) with the optimized $\chi^{n+1}\rightarrow$ Update $\bu^{n+1}$
\Comment{$\text{\it Perform Relax step (III)}$}
\EndFor 
\end{algorithmic}
\end{algorithm*}

Algorithm \ref{alg:opt-efr-chi} presents the general pseudocode 
that can be applied for all the algorithms presented in Table \ref{tab:acronyms}, by replacing the objective function $\mathcal{L}(\chi^{n+1})$ with the specific functional of the proposed algorithm, namely \eqref{eq:local-loss} or \eqref{eq:global-loss}. 

%% file: sections/methods-Opt-EF.tex
In the present section, we \RA{provide the details on the optimized} \emph{Evolve--Filter} algorithm, namely the final velocity at time step $t^{n}$ coincides with the filtered velocity: $\bu^{n} \equiv \overline{\bw}^{n}=\bu^{n}(\delta^n)$. This is also equivalent to the EFR case with $\chi=1$, and in this case, the predicted velocity would only depend on the filter parameter. 

The motivation for proceeding with this analysis is that the optimization of the relaxation parameter is influenced by the choice of the filter parameter, which is kept fixed in $\chi$-Opt-EFR.  
To overcome this limitation, in $\delta$-Opt-EF we optimize the time-varying parameter $\delta(t)$ in the EF approach, directly. Moreover, we want to compare how the optimization behaves in absence of the relaxation step, namely when only the filtered velocity is retained at each time step, so we fix $\chi=1$.

The $\delta$-Opt-EF strategy can be specialized in different versions, according to the objective function we decide to optimize, just as described in Section \ref{sec:ml-efr} for $\chi$-Opt-EFR. The cases taken into account are reported in Table \ref{tab:acronyms-EF}. As will be seen in Section \ref{sec:results-chi-opt}, the flow fields obtained considering the \emph{local} objective function are inaccurate and contain spurious oscillations. For this reason, for the $\delta$-Opt-EF objective function, we only consider the \emph{global} optimization.
%
\begin{table*}[h]
\centering
\caption{Acronyms and main features of the algorithms $\delta$-Opt-EF. Gray cells indicate that the algorithm does not present that feature.}
\label{tab:acronyms-EF}
{
\begin{tabular}{cccccc}
\\ \hline
\multirow{2}{*}{\textbf{Acronym}} & \multicolumn{2}{c}{\textbf{Objective function type}} & \multicolumn{3}{c}{\textbf{Contributions in objective function}}
\\ \cline{2-3}\cline{4-6}
& local & global & $\bu$ discrepancy & $\nabla \bu$ discrepancy & $p$ discrepancy \\ \hline
\optEFglob & \cellcolor{gray!50} & $\checkmark$ & $\checkmark$ & \cellcolor{gray!50} & \cellcolor{gray!50}\\ \hline
\optEFglobGRAD & \cellcolor{gray!50} & $\checkmark$ & $\checkmark$ &$\checkmark$& \cellcolor{gray!50}\\ \hline
\optEFglobPRESS & \cellcolor{gray!50} & $\checkmark$ & $\checkmark$ &$\checkmark$ &$\checkmark$  \\ \hline
\end{tabular}}
\end{table*}

The loss type considered is the one expressed in \eqref{eq:global-loss}, with the contributions \eqref{eq:global-contributions}.
The difference is that $\bu^{n+1}=\overline{\bw}^{n+1}$ is obtained just from system \eqref{eqn:ef-rom-2}, and, hence, the contributions will be a function of $\delta$ instead of $\chi$, say:
\begin{equation}
    \label{eq:global-loss-delta}
\lossGLOB(\delta^{n}) = w_{\bu} \lossGLOBu(\delta^{n}) + w_{\nabla \bu} \lossGLOBgrad(\delta^{n}) + w_{p} \lossGLOBp(\delta^{n}),
\end{equation}
with $w_{\bu}$, $w_{\nabla \bu}$, and $w_p$ scalar quantities, and where the loss contributions are:
\begin{equation}
\label{eq:global-contributions-delta}
\begin{split}
    &\lossGLOBu(\delta^n) = \left| \dfrac{\|\bu^{n}(\delta^n)\|^2_{L^2(\Omega)} - \|\bu_{ref}^{n}\|^2_{L^2(\Omega)}}{\|\bu_{ref}^{n}\|^2_{L^2(\Omega)}} \right|, \\
   & \lossGLOBgrad(\delta^{n}) =\left| \dfrac{\| \nabla \bu^{n}(\delta^n)\|^2_{L^2(\Omega)} - \| \nabla \bu_{ref}^{n}\|^2_{L^2(\Omega)}}{\| \nabla \bu_{ref}^{n}\|^2_{L^2(\Omega)}} \right|,\\
   &\lossGLOBp(\delta^{n}) = \left| \dfrac{\|p^{n+1}(\delta^n)\|^2_{L^2(\Omega)} - \|p_{ref}^{n+1}\|^2_{L^2(\Omega)}}{\|p_{ref}^{n+1}\|^2_{L^2(\Omega)}} \right|.
\end{split}
\end{equation}

The process is outlined in Algorithm \ref{alg:opt-ef-delta}.

\begin{algorithm*}[htpb!]
\caption{Pseudocode for algorithms $\delta$-Opt-EF optimizing the filter parameter $\delta(t)$. Row 10 (underlined) is considered only for \optEFglobPRESS{}.The optimization part is highlighted in gray.} 
\label{alg:opt-ef-delta}
\begin{algorithmic}[1]
\State $\bu_{ref}^n$, $n=0, \dots, N_T - 1$; $N_h^u$; $\bu^0$, $\Delta t_{opt}$, $\delta^0 = \eta$ \Comment{$\text{\it Inputs needed}$}
\For{$n \in [0, \dots, N_T - 1]$}
\State (I) $\rightarrow \bw^{n+1}$ \Comment{$\text{\it Perform Evolve step}$}
\State $\bu_{ref}^{n+1} \gets \bu_{ref}^{n+1}\,  | \, \mathbb{U}^{N_h^u}$ \Comment{$\text{\it Project the reference field into the coarse mesh}$}
\State $\delta^{n+1} \gets \delta^n$ \Comment{$\text{\it Initialize}$ $\delta$}
\highlight{
\If {$\bmod(n, \Delta t_{opt}) = 0$} \Comment{$\text{\it Perform optimization every}$ $\Delta t_{opt}$}
\State $iter=0$
\Repeat \Comment{$\text{\it SLSQP algorithm}$}
        \State (II) $\rightarrow \bu^{n+1}$ \Comment{$\text{\it Filter step}$}
        \State \underline{(I) $\rightarrow p^{n+2}$} \Comment{$\text{\it \underline{(Second) Evolve step}}$}
        \State Optimize $\mathcal{L}(\delta^{n+1})$, subject to $\num{1e-5} \leq \delta^{n+1} \leq \num{1e-3}$
        \State $iter \gets iter+1$
        \State Output current solution: $\delta^{n+1}$
    \Until{Convergence or maximum iterations reached}
    \EndIf
    }
    \State Solve (II) with the optimized $\delta^{n+1} \rightarrow$ Update $\bu^{n+1}$
\Comment{$\text{\it Perform Filter step (II)}$}
\EndFor 
\end{algorithmic}
\end{algorithm*}

In this optimization, we consider the Kolmogorov length scale $\eta$ as the initial guess at time $t^0$, namely, the literature common choice \cite{bertagna2016deconvolution}. The optimization is performed every $\Delta t_{opt}$, as in Opt-EFR, and at each time step, the initial guess is the optimal solution at the previous time instance.
The results related to this type of optimization are presented in Section \ref{sec:results-delta-opt}.

%% file: sections/methods-D-Opt-EFR.tex
This final methodological section is dedicated to the \emph{double} optimization algorithm\RA{, the} $\delta \chi$-Opt-EFR.

Also in this case we include Table \ref{tab:acronyms-double} describing the cases investigated. Since the results obtained with the local optimizations $\chi$-Opt-EFR$_{\text{local}}$ and $\delta$-Opt-EF$_{\text{local}}$ are not accurate, we decide to exclude these analyses here and to consider only the global optimizations in Table \ref{tab:acronyms-double}.

\begin{table*}[h]
\centering
\caption{Acronyms and main features of the algorithms $\delta \chi$-Opt-EFR. Gray cells indicate that the algorithm does not present that feature.}
\label{tab:acronyms-double}
{
\begin{tabular}{cccccc}
\\ \hline
\multirow{2}{*}{\textbf{Acronym}} & \multicolumn{2}{c}{\textbf{Objective function type}} & \multicolumn{3}{c}{\textbf{Contributions in objective function}}
\\ \cline{2-3}\cline{4-6}
& local & global& $\bu$ discrepancy & $\nabla \bu$ discrepancy & $p$ discrepancy \\ \hline
\DoptEFRglobGRAD & \cellcolor{gray!50} & $\checkmark$ & $\checkmark$ &$\checkmark$& \cellcolor{gray!50}\\ \hline
\DoptEFRglobPRESS & \cellcolor{gray!50} & $\checkmark$ & $\checkmark$ &$\checkmark$ &$\checkmark$  \\ \hline
\end{tabular}}
\end{table*}
The double optimization is described
in Algorithm \ref{alg:opt-double}, where the losses have the same form as in \eqref{eq:global-loss} and \eqref{eq:global-contributions}, but with double dependency $\mathcal{L}(\delta, \chi)$. 
In particular, the expression of the loss can be written as:

\begin{equation}
\begin{split}
    \label{eq:global-loss-double}
&\lossGLOB(\delta^n, \chi^{n}) = w_{\bu} \lossGLOBu(\delta^n, \chi^{n}) + w_{\nabla \bu} \lossGLOBgrad(\delta^n, \chi^{n}) +\\
&+w_{p} \lossGLOBp(\delta^n, \chi^{n}),
\end{split}
\end{equation}

where the global contributions are
\begin{equation}
\label{eq:global-contributions-double}
\begin{split}
    &\lossGLOBu(\delta^n, \chi^n) = \left| \dfrac{\|\bu^{n}(\delta^n, \chi^n)\|^2_{L^2(\Omega)} - \|\bu_{ref}^{n}\|^2_{L^2(\Omega)}}{\|\bu_{ref}^{n}\|^2_{L^2(\Omega)}} \right|, \\
   & \lossGLOBgrad(\delta^n, \chi^{n}) =\left| \dfrac{\| \nabla \bu^{n}(\delta^n, \chi^n)\|^2_{L^2(\Omega)} - \| \nabla \bu_{ref}^{n}\|^2_{L^2(\Omega)}}{\| \nabla \bu_{ref}^{n}\|^2_{L^2(\Omega)}} \right|,\\
   &\lossGLOBp(\delta^n, \chi^{n}) = \left| \dfrac{\|p^{n+1}(\delta^n, \chi^n)\|^2_{L^2(\Omega)} - \|p_{ref}^{n+1}\|^2_{L^2(\Omega)}}{\|p_{ref}^{n+1}\|^2_{L^2(\Omega)}} \right|.
\end{split}
\end{equation}

In this section we will only analyze cases $w_{\nabla \bu}= w_{\bu} =1$, $w_{p}=0$, and $w_{\nabla \bu}= w_{\bu} =w_p=1$, as those were proved to be the algorithms with the best performances in the $\chi$-Opt-EFR and the $\delta$-Opt-EF setting, as shown in Sections \ref{sec:results-chi-opt} and \ref{sec:results-delta-opt}.

\begin{algorithm*}[htpb!]
\caption{Pseudo-code for algorithms $\delta \chi$-Opt-EFR optimizing both parameters. Row 9 (underlined) is considered only for \DoptEFRglobPRESS{}. The optimization part is highlighted in gray.}
\label{alg:opt-double}
\begin{algorithmic}[1]
\State $\bu_{ref}^n$, $n=0, \dots, N_T - 1$; $N_h^u$; $\bu^0$, $\Delta t_{opt}$, $(\delta^0, \chi^0) = (\eta, 5\Delta t)$ \Comment{$\text{\it Inputs needed}$}
\For{$n \in [0, \dots, N_T - 1]$}
\State (I) $\rightarrow \bw^{n+1}$ \Comment{$\text{\it Perform Evolve step}$}
\State $\bu_{ref}^{n+1} \gets \bu_{ref}^{n+1}\,  | \, \mathbb{U}^{N_h^u}$ \Comment{$\text{\it Project the reference field into the coarse mesh}$}
$(\chi^{n+1}, \delta^{n+1}) \gets (\chi^n, \delta^n)$  \Comment{$\text{\it Initialize both parameters}$}
\highlight{\If {$\bmod(n, \Delta t_{opt}) = 0$} \Comment{$\text{\it Perform optimization every }\Delta t_{opt}$} 
\State $iter=0$
\Repeat 
        \State (II)+(III) $\rightarrow {\bu}^{n+1}$ \Comment{$\text{\it Filter and Relax step}$}
        \State \underline{(I) $\rightarrow p^{n+2}$} \Comment{$\text{\it \underline{(Second) Evolve step}}$}
        \State Optimize $\mathcal{L}(\delta^{n+1}, \chi^{n+1})$, subject to $\num{1e-5} \leq \delta^{n+1} \leq \num{1e-3}$, $0 \leq \chi^{n+1} \leq 1$ 
        \State $iter \gets iter + 1$
        \State Output current solution: $(\delta^{n+1}, \chi^{n+1})$
    \Until{Convergence or maximum iterations reached}
    \EndIf}
    \State Solve (II)+(III) with the optimized $\delta^{n+1}$ and $\chi^{n+1} \rightarrow$ Update $\bu^{n+1}$  \Comment{$\text{\it Perform Filter and Relax}$}
\EndFor 
\end{algorithmic}
\end{algorithm*}

In the $\delta \chi$-Opt-EFR optimization, we consider as ranges for the parameters' search the same intervals taken into account in strategies $\chi$-Opt-EFR and $\delta$-Opt-EF, namely $\delta(t) \in [\num{1e-5}, \num{1e-3}]$, $\chi(t) \in [0, 1]$.
The results related to this type of optimization are presented in Section \ref{sec:results-double-opt}. We remark that line 9 of Algorithm \ref{alg:opt-double}, is performed only if $w_p \neq 0$.